**The Symmetric Traveling Salesman Problem**

**by Howard Kleiman**

## I. INTRODUCTION

Let M be an $n \ X \ n$ symmetric cost matrix where $n$ is even. We present an algorithm that extends the concept of admissible permutation and the modified Floyd-Warshall algorithm given in math. CO/0212305 that was used to obtain to obtain near-optimal or optimal solutions to the asymmetric traveling salesman problem. Using the modified F-W, we obtain a derangement, $D_{minimal}$, that we patch into a tour, $T_{UPPERBOUND}$. Once we have obtained $T_{UPPERBOUND}$, we consider it to be a circuit consisting of edges. From it, we extract two sets of alternate edges. Each set forms a perfect matching. We denote by $\sigma_{T_{UPPERBOUND}}$ the perfect matching whose edges have a smaller value. $\sigma_{T_{UPPERBOUND}}$ can also be expressed as a product of $\frac{n}{2}$ pair-wise disjoint 2-cycles. We then construct $\sigma^{-1}_{T_{UPPERBOUND}} M^{-}$.

Applying the modified F-W algorithm to it, using sieve criteria, we can obtain paths (called either acceptable or 2-circuit paths) that can form circuits. We prove that every tour whose value is less than $|T_{UPPERBOUND}|$ and is constructed from circuits obtained by using F-W on $\sigma^{-1}_{T_{UPPERBOUND}} M^{-}$ can be obtained by patching acceptable and 2-circuit cycles. Let $S$ be a set of cycles that can be patched to form a tour. Formulae are derived that give the number of points contained in such a set of cycles as well as a bound for $|S|$. To be more precise, we use Phases 1 and 2 given in math. CO/0212305 with the following modification: Since a 2-cycle in a derangement is an edge, we must check to see if an arc extending a path in the algorithm leads to an arc symmetric to an arc in the current derangement. We disallow such paths. In Phase 2, we use the modified Floyd-Warshall algorithm to obtain an approximation to a minimally-valued derangement. During this process, we *delete* all arcs that are symmetric to the current derangement, $D_i$, when we apply the modified F-W algorithm to $D_i^{-1}M^{-}$. Using a transformation of Jonker and Volgenant [2] given in paper of Helsgaun [1], the algorithm can be applied to $n \ X \ n$ asymmetric cost matrices. If $n$ is odd, we can add a new column and new row to M. Each of the entries in the new row and column has the value – N where N is the entry of greatest value in M. Also, the algorithm in [1] could provide a good option (in place of Phases 1 and 2 and the construction of a tour) for $T_{UPPERBOUND}$.

## II. DEFINITIONS.

An *edge* is an undirected line segment connecting two points, say $a$ and $b$. A *circuit* is a set of edges and points such that precisely two edges are incident to each point. A *tour* in an $n \ X \ n$ symmetric cost matrix M is a circuit containing precisely $n$ points and $n$ edges where each edge has a value or cost. The *value of a tour* is the sum of the values of its edges. An *arc* is an edge that has been given an orientation or direction. Arc $(b\ a)$ is *symmetric* to $(a\ b)$. A *cycle* is a directed circuit. The *value of a cycle* is the sum of the values of its arcs. A *permutation* is a rearrangement of some or all of a fixed set of points, say $V = \{1,2, \ldots ,n\}$. A *derangement* is a permutation that moves all of the points in $V$. A *perfect matching or PM* is a set of $\frac{n}{2}$ edges that are pair-wise disjoint. *Henceforth, it is equivalent to a set of pair-wise disjoint 2-cycles where each edge is equivalent to a pair of symmetric entries in M, i.e. to a derangement consisting of 2-cycles* . $\sigma^{-1}M$ is a permutation of the columns of M by the derangement $\sigma^{-1}$. $\sigma^{-1}M^{-}$ is the matrix obtained from $\sigma^{-1}M$ by subtracting the value in each diagonal entry $(a\ a)$ from all entries in row $a$. A *path* in $\sigma^{-1}M^{-}$ is a sequence of arcs of the form $[a_1\ a_2 \ldots a_r]$ where $a_i \neq a_j, i = 1,2, \ldots ,r;\ j = 1,2, \ldots ,r.$ A *cycle* in $\sigma^{-1}M^{-}$ is a path that becomes a cycle. If $C$ is a cycle in $\sigma^{-1}M^{-}$, then $\sigma C$ always yields a derangement. In the symmetric case, this isn't satisfactory since a 2-cycle of $\sigma C$ is an edge. Given a symmetric cost matrix M, our purpose is to obtain a *circuit* consisting of edges. We thus define acceptable and 2-circuit paths. An *acceptable path* is a path satisfying the following condition: It yields an *acceptable path of edges* of the form $[a_1\ \sigma(a_2)\ a_2\ \sigma(a_3)\ a_3\ \ldots\ \sigma(a_r) a_r]$ where no two points belong to the same 2-cycle of $\sigma$. Due to symmetry in M, every non-diagonal entry in M defines an edge having double the value of the entry. Thus, by construction, since both $[a_i \sigma(a_{i+1})]$ and $[\sigma(a_{i+1}) a_{i+1}]$ lie in M, they may be thought of as edges if we wish to do so. An *unlinked 2-circuit path* is a path in which we allow precisely *one* pair of points to belong to the same 2-cycle of $\sigma$. In a *linked 2-circuit path*, we allow precisely *two* pairs of points to belong to the same 2-cycles provided that *the points of the 2-cycles interlace*. i.e., as we traverse the path, we don't have a point from a 2-cycle followed by the other point of the 2-cycle. An *acceptable cycle* is an acceptable path that forms a cycle. An *unlinked 2-circuit cycle* is an unlinked 2-circuit path that forms a cycle. A *linked 2-circuit cycle* is a linked 2-circuit path that forms a cycle. A *tour* is a circuit in M containing $n$ points and edges. Since we consider a perfect matching as a set of 2-cycles, we here consider a tour as a *directed* set of edges. At the same time, when applying our algorithm, we use the fact that it *is* constructed of edges. $T_{FWTSPOPT}$ denotes a minimally-valued tour that can be obtained by applying the modified F-W algorithm to M. An *optimal tour* in M which may require the use of paths not constructed by using minimal paths obtained from F-W is denoted by $T_{TSPOPT}$. Let $C$ be a cycle obtained from $\sigma^{-1}M^{-}$ where $\sigma$ is a product of $\frac{n}{2}$ pair-wise disjoint 2-cycles. Then $|C|$ is the sum of the values of the arcs of $C$. Let $a_i$ be a point of a cycle, $C$, constructed in $\sigma^{-1}M^{-}$ such that as we traverse $C$ in a clock-wise direction, the partial sum of the arcs of the arcs obtained so far is no greater than $|C|$. Then $a_i$ is a *determining vertex* of $C$.

A set of pair-wise disjoint acceptable cycles such that no two cycles have a point in the same 2-cycle of $\sigma_T^{-1}M^-$ is called an *acceptable permutation*. $pt(C)$ is the number of points in the permutation cycle $C$. $\sigma_{ABSOLUTE}$ is the minimally-valued PM in M.

## III. USEFUL THEOREMS

**Theorem 1** *Let* $C = (a_1 a_2 \ldots a_n)$ *be a cycle of length n. Assume that the weight* $w(a_i, C(a_i))$ *(i=1,2,...,n) corresponds to the arc* $(a_i, C(a_i))$ *of C. Then if*

$$W = \sum_{i=1}^{i=n} w(a_i, C(a_i)) \leq 0$$

*there exists at least one vertex* $a_{i*}$ *with* $1 \leq a_{i*} \leq n$ *such that*

$$S_m = \sum_{j=0}^{j=m} w(a_{i*+j}, C(a_{i*+j})) \leq 0 \qquad (A)$$

*where* $m = 0,1,2,\ldots,n\text{-}1$ *and* $i*+j$ is modulo $n$.

**Proof.** We prove the theorem by induction. Let k = 2. We thus have a 2-cycle. If both arcs have non-positive value, then the theorem is proved. If the non-positive arc, $(a_1\ a_2)$, has a smaller weight than a positive one, then the sum of the weights of the two arcs is positive. This can't be the case. Thus, the sum of the two weights is non-positive. Now let the theorem always be true when our cycle has k arcs. Suppose the cycle has k+1 arcs. In what follows, assume that if a value is 0, its sign is negative. Then one of the following is true: (a) there exists a pair of consecutive arcs both of whose values have the same sign or (b) the signs of values of the arcs consecutively alternate in sign. First, we consider (a). Without loss of generality, let the two arcs be $(a_1\ a_2)$ and $(a_2\ a_3)$. Assume that each one has a non-positive value . We now define the arc $(a_1\ a_3)$ where $w(a_1\ a_3) = w(a_1\ a_2) + w(a_2\ a_3)$. Now replace arcs $(a_1\ a_2)$ and $(a_2\ a_3)$ by $(a_1\ a_3)$. The result is a cycle $C'$ containing $k$ arcs. By induction, the theorem holds for $C'$. Now replace $(a_1\ a_3)$ by $(a_1\ a_2)$ and $(a_2\ a_3)$. Let $a_{i*}$ be a determining vertex of $C'$. Then the path, $P_{i*}$, from $a_{i*}$ to $a_1$ is non-positive. If both $w(a_1\ a_2)$ and $w(a_2\ a_3)$ are non-negative, then both $P_{i*} \cup w(a_1\ a_2)$ *and* $P_{i*} \cup w(a_1\ a_2) \cup w(a_2\ a_3)$ are non-positive. Thus, the theorem is valid in this case. Now assume that both $w(a_1\ a_2)$ and $w(a_2\ a_3)$ are positive. But, since the theorem is valid for $C'$, $P_{i*} \cup w(a_1\ a_3)$ is non-positive. This assures us that each of $P_{i*} \cup w(a_1\ a_2)$ and $P_{i*} \cup w(a_1\ a_2) \cup w(a_2\ a_3)$ is non-positive. Therefore, the theorem holds when the values of $(a_1\ a_2)$ and $(a_2\ a_3)$ both have the same sign. We now consider case (b). Here, the signs alternate between positive and negative. Since the arcs lie on a cycle, we may assume that the first sign is negative. If k+1 is odd, it is always true that at least one pair of consecutive arcs has the same sign. Thus, assume that k+1 is even. Then, starting with a non-positively-valued arc, we can arrange the arcs of the cycle in pairs where the first arc has a non-positive value, while the second one has a positive value. Suppose the sum of each pair of arcs is positive. Then the sum of the values of the cycle is positive. Therefore, there exists at least one pair of arcs the sum of whose values is non-positive. Remembering that the first arc of each pair is non-positive, we follow the same procedure as in (a): $C'$ is a non-positive cycle containing k arcs. Therefore, $C$ is also positive.

**Corollary 1**. *Suppose that C is a cycle such that*

$$W = \sum_{i=1}^{i=n} w(a_i, C(a_i)) \leq N \qquad (B)$$

*Then there exists a determining vertex* $a_{i*}$ *such that each partial sum,* $S_m$*, has the property that*

$$S_m = \sum_{j=0}^{j=m} w(a_{i*+j}, C(a_{i*+j})) \leq N \qquad (C)$$

*always holds . Here m =0,1,2, ... , n-1 while* $i*+j$ *is modulo n.*

**Proof**. Subtract $N$ from both sides of (B). Now let the weight of arc $(a_n \; C(a_n))$ become $w(a_n, C(a_n)) - N$. From theorem 3.1,

$$W* = W - N = \sum_{i=1}^{i=n} w(a_i, C(a_i)) \leq 0 \quad \text{(D)}$$

Therefore, we can obtain a determining vertex $a_{i*}$ having the property that every partial sum having $a_{i*}$ as it initial vertex is non-positive . It follows that if we restore $w(a_n, C(a_n))$ to its original value, every partial sum with initial vertex $a_{i*}$ is less than or equal to $N$ .

**Corollary 2**. *Let* $C$ *be a positively-valued cycle of length* $n$ *obtained from a cost matrix* $W$ *whose entries may be positive, negative or* $0$ *. Then there exists at least one determining vertex, say* $a_{i*}$*, of* $C$ *such that each subpath* $S_m$ *having initial vertex* $a_{i*}$ *has the property that*

$$S_m = \sum_{j=0}^{j=m} w(a_{i*+j}, C(a_{i*+j})) \geq 0 \qquad (E)$$

*where m = 0,1,2, ... ,n-1 and* $i*+j$ *is modulo n.*

**Proof.** Let $N = 0$ in ($B$). Then multiply each term of $(B)$ by $-1$. It follows that we obtain $.-W = \sum_{i=1}^{i=n} -w(a_i, C(a_i)) \leq 0$. The rest of the proof is similar to that of Corollary 1.

**Example 1** We now give an example of how to obtain $i = i'$.

Let $C = (a_1\ a_2\ \ldots\ a_n)$

$a_1 = -7,\ a_2 = -10,\ a_3 = +1,\ a_4 = +2,\ a_5 = -7,\ a_6 = +4,\ a_7 = -9,$

$a_8 = +11,\ a_9 = -2,\ a_{10} = -1,\ a_{11} = -4,\ a_{12} = -4,\ a_{13} = -8,\ a_{14} = +9,$

$a_{15} = +9,\ a_{16} = +21,\ a_{17} = +1,\ a_{18} = -2,\ a_{19} = -1,\ a_{20} = -3,$

$a_{21} = -3,\ a_{22} = -12,\ a_{23} = +6,\ a_{24} = +2,\ a_{25} = +3$

-7 –10 +1 +2 –7 +4 –9 +11 –2 –1 –4 –4 –8 +9 +9 + 21 +1

-2 -1 -3 -3 -12 +6 +2 +3

We now add terms with like signs going from left to right. We place the ordinal number of the first number in each sum above it.

1 3 5 6 7 8 9 14 18 23

-17 +3 -7 +4 -9 +11 -19 +40 -21 +11

We next add the positive number to the right of each negative number to the negative number.

1 5 18

-14 -3 +2 +21 -10

We now add terms with like signs going from left to right.

1 18

-17 +23 -10

Finally, assuming that all points lie on a circle, we add like terms going from left to right. We thus obtain

18 18

-27 +23 = -4

This tells us that $i'$ is the eighteenth ordinal number; its value is –2.

Thus, the partial sums are: -2, -3, -6, -9, -21, -15, -13, -10, -17, -27,

-26, -24, -31, -27, -36, -25, -27, -28, -32 –36, -44, -35, -26, -5, -4

**Theorem 1.2 The Floyd-Warshall Algorithm**

*If we perform a triangle operation for successive values j = 1,2,...,n, each entry $d_{ik}$ of an n X n cost matrix M becomes equal to the value of the shortest path from i to k provided that M contains no negative cycles.*

The version given here is modeled on theorem 6.4 in [4].

**Proof.** We shall show by induction that that after the triangle operation for $j = j_0$ is executed, $d_{ik}$ is the value of the shortest path with intermediate vertices $v \le j_0$, for all *i* and *k*. The theorem holds for $j_0 = 1$ since $v = 0$. Assume that the inductive hypothesis is true for $j = j_0 - 1$ and consider the triangle operation for $j = j_0$:

$$d_{ik} = min\{d_{ik}, d_{ij_0} + d_{j_0k}\}.$$

If the shortest path from *i* through *k* with $v \le j_0$ doesn't pass through $j_0$, $d_{ik}$ will be unchanged by this operation, the first argument in the min-operation will be selected, and $d_{ik}$ will still satisfy the inductive hypothesis. On the other hand if the shortest path from *i* to *k* with intermediate vertices $v \le j_0$ does pass through $j_0$, $d_{ik}$ will be replaced by $d_{ij_0} + d_{j_0k}$. By the inductive hypothesis, $d_{ij_0}$ and $d_{j_0k}$ are both optimal values with intermediate vertices $v \le j_0 - 1$. Therefore, $d_{ij_0} + d_{j_0k}$ is optimal with intermediate vertices $v \le j_0$.

We now give an example of how F-W works.

**Example 2** Let d(1, 3) = 5, d(3, 7) = -2, d(1, 7) = 25. Then

d(1, 3) + d(3, 7) < d(1, 7). Note however, that the intermediate vertex 3 comes from the fact that we have reached column j = 3 in the algorithm. We now substitute d(1, 3) + d(3, 7) = 3 for the entry in

(1, 7). Suppose now that d(1, 10) = 7 while d(7, 10) = -5.

$$d(1, 3) + d(3, 7) + d(7,10) = -2 < d(1, 10) = 7.$$

In what follows, Roman numerals represent numbers of iterations.

*P*

**I.** [1 3] + 3

[1 7] + 4

[1 13] + 6

[1 15] + 2

[1 19] + 4

[1 20] + 1

------------------------

**II.** [1 18] + 18

------------------------

**III.** [1 14] + 16

------------------------

*C*

**I.** [20 18] + 18

------------------------

**II.** [20 14] + 16

------------------------

**III.** [20 6] + 12

[20 7] + 1

[20 13] + 6

[20 15] + 2

[20 19] + 4

[20 20] + 1

**IV.** [1 6] + 12 -------------------------

[1 7] + 1 + **60**

------------------------

+ **67**

**Theorem 3** *Let c be any real number,* $S_a = \{ a_i \mid i = 1,2, \dots ,n\}$, *a set of real numbers in increasing order of value. For i = 1, 2, ..., n, let* $b_i = a_i + c$. *Then* $S_b = \{ b_i \mid i = 1,2, \dots ,n\}$ *preserves the ordering of* $S_a$.

**Proof**. The theorem merely states that adding a fixed number to a set ordered according to the value of its elements retains the same ordering as that of the original set.

*Comment*. This is very useful when we are dealing with entries of the value matrix which have not been changed during the algorithm. However, as the algorithm goes on, if we have entry ( i, j ) and
d(i, j), we must go through all values of j = j' where the (i, j')-th entry's value in the current $\sigma_\alpha^{-1} M^-(k)$ is less than the value of the (i, j')-th entry in $D^{-1}M^-$ (The terminology will become clearer in the examples given.)

**Theorem 4**. *Let* $D = \prod_{i=1}^{i=r} C_i$ *be a derangement in* $S_n$. *Assume that* $a_i$ $(i = 1, 2, \dots , r)$ *is a point on* $C_i$. *Then* $s = (a_1\ a_2\ \dots\ a_r)$ *has the property that* $Ds$ *is an* $n$*-cycle,* $H$.

**Proof.** Consider the arcs of $H$ obtained by the action of $D$ on $s$:
$(a_1\ a_2) \to (a_1\ D(a_2)),\ (a_2\ a_3) \to (a_2\ D(a_3)), \dots ,\ (a_r\ a_1) \to (a_r\ D(a_1))$. Thus,

$$H = (a_1D(a_2) \dots a_2D(a_3) \dots a_3D(a_4) \dots\dots a_{r-1}D(a_r) \dots a_rD(a_1) \dots a_1)$$

. **Theorem 5** Let $T$ be a tour containing an even number of edges. Let $\sigma_T$ be the set of smaller-valued alternating edges of $T$ that may be considered a set of $\frac{n}{2}$ pair-wise disjoint 2-cycles. Then there always exists an acceptable cycle, $s$, in $\sigma_T^{-1}M^-$ containing $\frac{n}{2}$ points such that $\sigma_T s = T^{-1} = T$.

**Proof**. Let $T = (a_1\ a_2\ a_3\ \dots\ a_{2n})$. $\sigma_T s = T \Rightarrow s = \sigma_T T$. Since $\sigma_T$ consisting of alternating edges of $T$, assume it is $S = \{[a_i\ a_{i+1}] \mid i = 1,3,5,\dots,2n-1\}$. Consider the cycle $s = (a_1\ a_{2n-1}\ a_{2n-3}\ \dots\ a_3)$ in $\sigma_T^{-1}M^-$. Applying $\sigma_T$ to it, where the arc $(a_i\ a_j)$ is mapped into $[a_i\ \sigma_T(a_j)][\sigma_T(a_j)\ a_j]$ we obtain $[a_1\ \underline{a_{2n}}\ a_{2n-1}\ \underline{a_{2n-2}}\ a_{2n-3}\ \dots\ a_3\ \underline{a_2}\ ]$. Thus, $s$ is an acceptable cycle containing $\frac{n}{2}$ points that yields the tour $T^{-1}$ where – due to symmetry - $|T^{-1}| = |T|$.
Furthermore, since $T^{-1}$ is a circuit consisting of edges, $T^{-1} = T$.

**Theorem 6** *Let n = 2m. A perfect matching always exists if a finite-valued hamilton cycle exists.*

**Theorem 7** *Let n = 2m. Assume that* $\sigma_{ABSOLUTE}$ *is the smallest-valued PM that we've obtained. Then*

*a necessary and sufficient condition for an acceptable cycle* $C = (a_1\ a_2\ \ldots\ a_m)$ *in* $\sigma^{-1}_{ABSOLUTE}M^{-1}$ *to yield a perfect matching is that the cycle* $C' = (\sigma_{ABSOLUTE}(a_2)\ \sigma_{ABSOLUTE}(a_1)\ \sigma_{ABSOLUTE}(a_m)\ \ldots\ \sigma_{ABSOLUTE}(a_3))$ *is disjoint from* $C$. *Furthermore, the arcs of the permutation obtained from* C *and* $C = (a_1\ a_2\ \ldots\ a_m)$ *yield a perfect matching, while*

$$T = (a_1\, \sigma_{ABSOLUTE}(a_2)\, a_2\, \sigma_{ABSOLUTE}(a_3)\, a_3\ \ldots\ \sigma_{ABSOLUTE}(a_m)\, a_m\, \sigma_{ABSOLUTE}(a_1))$$

*yields a tour.*

**Corollary 1** *If the length of C is* $m' < m$, *then T is a circuit of edges of length 2m'.*

**Proof.** From theorem 5, we know that if $n = 2m$, then we obtain a tour. It follows from theorems 5 and 7, if $m' < m$, we obtain a cycle, $C$, consisting of edges. The number of edges cannot be $n$ since $2m' < 2m = n$. Thus, $C$ cannot be a tour.

*Note 1.* Although it might be time consuming, if we can obtain a set of disjoint cycles, $C_i$, covering all $n$ vertices, we could obtain a minimally-valued derangement. This would be the best possible lower bound for $\sigma_{TSPOPT}$, i.e., a minimal derangement obtained using the modified F-W algorithm, is the smallest possible derangement obtainable: It is a lower bound for both $\sigma_{TSPOPT}$ as well as for $\sigma_{FWTSPOPT}$.

**Theorem 8.** Every acceptable cycle in $\sigma^{-1}M^{-}$ yields a unique perfect matching.

**Proof**. We first note that $\sigma(\sigma(a)) = a$. Let $C = (a_1\ a_2\ \ldots\ a_m)$. Define $C' = (\sigma(a_2)\ \sigma(a_1)\ \sigma(a_m)\ \ldots\ \sigma(a_3))$ $C'$ is an acceptable cycle since it has the same number of points as $C$ and its points are in the same 2-cycles as those of $C$. Consider the directed edges obtained from $C$: $\{(a_1\ \sigma(a_2)), (a_2\ \sigma(a_3)),\ \ldots\ , (a_m\ \sigma(a_1))\}$. They are pair-wise disjoint and contain precisely the points in the 2-cycles containing points of $C$. On the other hand, the directed edges from $C'$ are $\{(\sigma(a_2)\ a_1), (\sigma(a_1)\ a_m),\ \ldots\ , ((\sigma_3)\ a_2)\}$. The second set contains precisely all directed edges that are symmetric to those in the first set. Thus, we have obtained a new set of 2-cycles that replace those in which $C$ had points. By the definition of acceptable permutation, we have changed the edges only in the 2-cycles of $\sigma$ containing a point of $C$. Thus, $\sigma C C' = \sigma^*$, a unique perfect matching.

**Theorem 9** *Let* $\sigma$ *be a perfect matching. Then a necessary and sufficient condition that it is a minimal-valued perfect matching is that* $\sigma^{-1}M^{-}$ *contains no negatively-valued acceptable cycle.*

**Proof.** Suppose $C$ is a negatively-valued acceptable cycle. Then $\sigma C$ yields a perfect matching of smaller value than $\sigma$. On the other hand, if no acceptable cycle exists, then $\sigma$ is a minimally-valued perfect matching.

**Theorem 10** *Let* $\sigma$ *be a PM. Then the following hold:*

*(1) An acceptable circuit containing an even number of vertices is obtainable from an acceptable cycle in* $\sigma^{-1}M^{-}$.

*(2) A 2- circuit cycle in* $\sigma^{-1}M^{-}$ *contains in its alternating path precisely two circuits each of which has an odd number of edges.*

*(3) Tthe first vertex of each of the two circuits along the alternating path belongs to the same 2-cycle of* $\sigma$.

*Furthermore, if* $s$ *is a permutation such that* $\sigma s$ *is a tour, then each of the cycles of* $s$ *is either an acceptable or a 2-circuit cycle.*

**Proof.** An acceptable cycle always yields an acceptable circuit containing twice the number of edges as the cycle contains. To insure that a path is a 2 circuit path, it must have distinct vertices in $P_{20i}$, while its alternating path contains precisely one circuit of the following kind: $P = (a\ q) = [a\ \underline{a}'\ b\ \underline{b}' ...\ \underline{a}\ ...\ d\ \underline{d}'\ e\ \underline{e}'\ f\ \underline{f}'\ ...\ p\ \underline{p}'\ q]$. Our circuit consists of the edges $[a\ \underline{a}'\ b\ \underline{b}'\ ...\ \underline{a}]$. The arc $(a\ b)$ yields the edge $[a\ \underline{a}']$. But $(b\ a')$ is an arc of $\sigma$ when we consider it as a directed derangement. Thus, $(a'\ b)$ is an arc symmetric to the arc $(b\ a')$ of $\sigma$. It follows that any underlined vertex in $P$ is followed by its companion in a 2-cycle of $\sigma$. Now consider the cycle $C = (a\ b\ ...\ ?\ ...\ d\ e\ f\ ...\ p\ q\ ...\ r)$ that defines the alternating path of edges $P' = [a\ \underline{a}'\ b\ \underline{b}' ...\ \underline{a}\ ?\ ...\ d\ \underline{d}'\ e\ \underline{e}'\ f\ \underline{f}'\ ...\ p\ \underline{p}'\ q\ ...\ r\ \underline{??}\,a]$. What are the vertices represented by "?" and "??" ? Let $s$ be the companion of $a$ in a unique 2-cycle of $\sigma$. Then "??" represents $s$. But what represents "?" ? It could be $s$ since that would allow $C$ to consist of distinct vertices, while allowing us to obtain two circuits of edges from $P'$. Could it be another vertex? But if we had three or more circuits, the last vertex of the final circuit would have to be $s$ which is impossible. Thus, there can be only two circuits – each containing an odd number of vertices – in $P$. One of the edges of each circuit isn't actually obtained from a cycle. However, since it always is an edge of $\sigma$ and thus has a value of 0 in $\sigma^{-1}M^{-}$, we can assume that it is part of the circuit. Actually, the missing edge tells us precisely which edge of another circuit can be deleted. On the other hand, by including it in theory, the formula for the number of edges in all circuits - $(n + 2r - 2)$ - can be used without modification.

*Note.* In theorem 10, we have shown that if a cycle of $s$ is not acceptable, then it must be a 2 circuit cycle: No cycle can yield three or more circuits. We obtain a tour by constructing a tree of circuits in which circuits are linked by a common edge. (In the case of a 2 circuit cycle, we delete the "phantom edge" of a circuit.).

In what follows, if $C$ is an acceptable cycle, $C_{circuit}$ is the circuit of edges obtainable from $C$.

What follows is another – more detailed - version of theorem 10.

**Theorem 11** *Let* $\sigma$ *be a PM such that* $|\sigma_{ABSOLUTE}| \leq |\sigma| < |T_{UPPERBOUND}|$. *Furthermore, suppose that*

*there exists a tour T such that* $|T| < |T_{UPPERBOUND}|$. *Then*

*(a) there exists an acceptable cycle* $t$ *such that* $\sigma\, t = T$

*(b) there exists a permutation* $t'$ *such that* $\sigma_{ABSOLUTE}\, t' = T$ *where* $t'$ *consists of a set of disjoint cycles each of which is either acceptable and yields a circuit containing an even number of vertices, or else is a 2 circuit cycle yielding a set of two or more disjoint circuits each of which contains an odd number of points. In the latter case, the circuits can be obtained directly from 2-circuit paths.* .

*(c) Each circuit,* $C_{icircuit}$, *consists of edges with at most one identical to the edge of a fixed cycle* $C_{j_{circuit}}$ *where* $C_{icircuit} \neq C_{jcircuit}$, $j \neq i$, *while every circuit has an edge identical to at least one circuit.*

*(d) Defining a common edge as a "link" between two circuits, the circuits form a tree.*

*(e) The total number of edges (and vertices) in all of the circuits obtained is* $n + 2r - 2$ *where* $r$ *is the number of circuits.*

**Proof**. (a) From theorem 5, there exists a acceptable cycle $t$ whose value is non-positive such that $T_{UPPERBOUND}\, t^{-1} = \sigma$. Furthermore, each point of $t^{-1}$ is contained in precisely one 2-cycle of $\sigma$.

(b) If a cycle, $C'$, of $t'$ is acceptable, then, from theorem 9, it yields a circuit containing an even number of points. We know that $t'$ exists and contains disjoint cycles. Any cycle $C''$ that is not acceptable must still yield circuits that contain all of its vertices of the cycle. Otherwise, we couldn't obtain a tour. This can only be the case if it yields two or more circuits each of which (from theorem 10) contains an odd number of points.

(c) In order to obtain a tour from a set of circuits, $\{C_{icircuit}\}$, containing all $n$ vertices, each circuit, $C_{i_{circuit}}$, must contain a path $[a \ \ldots \ b]$ such that the edge $[a\ b]$ is identical to precisely one edge $[a\ b]$ of some circuit $C_{j_{circuit}}$ . It can't contain two or edges identical to those in $C_{j_{circuit}}$, since then we obtain at least two circuits. On the other hand, if a circuit doesn't have at least one edge in common with another circuit, then we again must obtain at least two circuits.

(d) If the circuits form a tree whose "links" are identical edges, then we can obtain a tour by deleting common edges. Otherwise, if they don't form a tree, we would obtain more than one circuit..

(e) Let $r$ be the number of circuits, $C_{i_{circuit}}$. We thus delete $2r-2$ common edges (and vertices) from the $r$ circuits forming the tree. After deleting the $2r-2$ common edges, we obtain a tour containing $n$ edges. Thus, the total number of edges (and vertices) in all circuits is $n + 2r - 2$.

**Theorem 11** Let $T$ be a tour where $\sigma_T$ is its set of alternating edges the sum of whose values is no greater than $|\frac{T}{2}|$. Suppose that $S$ is a set of acceptable and 2-circuit cycles obtained from $\sigma_T^{-1}M^-$ .Let $a$ be the number of acceptable cycles and $t$ the number of 2-circuit cycles in $S$. Assume that

$$C_{PATH} = [a\ \underline{f_1}\ f_1\ \ldots\ \underline{l_1}\ l_1\ \underline{a}\ b\ \underline{f_2}\ f_2\ \ldots\ \underline{l_2}\ l_2\ \underline{b}\ a]$$

is a 2-circuit path obtainable from a 2-circuit cycle. For the purpose of discussion, from the *2-circuit permutation cycle* $C$, define each of $(a\ f_1 \ldots l_1)$ and $(b\ f_2 \ldots l_2)$ as a cycle rather than as a circuit. Assume that two cycles of $S$ are *linked* if precisely one point in the first cycle belongs to the same 2-cycle of $\sigma_T$ to which a point in the second cycle belongs. We now add the condition that the points $a$ and $b$ of $C$ cannot be linked. Then necessary and sufficient conditions for the cycles in S to yield a tour $T'$ with $|T'| < |T|$ are the following:

(1) The number of points moved by the cycles in $S$ is $\frac{n}{2} + 3t + a - 1$ .

(2) Each cycle in $S$ is linked to at least one cycle of $S$.

(3) Each pair of cycles in $S$ has at most one point in common

(4) The cycles of $S$ form a tree by linking.

(5) The sum of the values of the cycles in $S$ is smaller than $|\ T\ |$.

**Proof.** When two cycles are linked, their points of linkage yield a common edge of circuits that can be constructed from the cycles. Depending upon whether or not the two cycles can be linked by going in a clockwise or counterclockwise direction of the circuit formed from the second cycle, two cycles may or may not have a point in common. The important thing is that only one *common* edge can be deleted from two linked cycles. Since the linked cycles form a tree, after deleting common edges from the circuits obtained by the cycles, we obtain a circuit containing all of the points of $S$. Let $p$ be the number of points moved by a 2-circuit cycle. Then the number of points in the two circuits obtained from the cycle is $2p-2$. The reason is as follows:

The circuits we obtain are ($a\ \underline{f_1}\ f_1\ \ldots\ \underline{l_1}\ l_1$), ( $b\ \underline{f_2}\ f_2\ \ldots\ \underline{l_2}\ l_2$). We thus omit the following arcs:

$\underline{a}\ b$, $\underline{b}\ a$. Let $A_i$, $i = 1, 2, \ldots, a$, be the subset of acceptable permutation cycles in $S$, while $T_j$, $j = 1, 2, \ldots, s$, is the subset of permutation 2-circuit cycles. We are assuming that we can obtain a tour by deleting arcs. Thus, the number of edges contained in all circuits constructed from the cycles in $S$ is $\sum_{i=1}^{a} 2|A_i| + \sum_{j=1}^{s}(2|T_j| - 2)$. In order to construct a tour from these cycles, we must construct $a + 2t - 1$ linkages of $a + 2t$ circuits. Each linkage deletes two edges. Thus, we must delete $2(a + 2t - 1)$ edges. Since a tour contains $n$ edges, we obtain the equation

$$\sum_{i=1}^{a} 2|A_1| + \sum_{j=1}^{s}(2|T_1| - 2) - 2(a + 2t - 1)$$

$$= 2(\sum_{i=1}^{a}|A_i| + \sum_{j=1}^{s}(|T_j|) - 2t - 2(a + 2t - 1) = n$$
$$\Rightarrow 2(\sum_{i=1}^{a}|A_i| + \sum_{j=1}^{s}(|T_j|)) = n + 2t + 2(a + 2t - 1) = n + 6t + 2a - 2$$

Therefore, the number of points in the corresponding cycles in $S$ is $\frac{n}{2} + 3t + a - 1$.

Next, each cycle yields one or more circuits composed of edges. The edges are of two kinds: (1) Those that we have obtained from arcs of the cycles. (2) Those we obtained from the edges (2-cycles) of $\sigma_T$. By definitions, each of the edges obtained from $\sigma_T$ has a value of zero. Furthermore, all of the edges that are deleted belong to $\sigma_T$. It follows that the value of one (or two) *directed* circuits corresponds to the value of a permutation cycle lying in $S$ implying that the sum of the values of the permutation cycles in $S$ corresponds to the sum of the values of the *directed* circuits obtained from corresponding permutation cycles.

**Corollary 1** Let $S$ be a set of cycles that yields a tour. Denote by $S'$ the set obtained by deleting $X$ acceptable and $Y$ 2-circuit cycles from $S$. Define $t' = t - Y$, $a' = a - X$. Then if $p'$ is the total number of points moved by the cycles in $S'$, $p' < \frac{n}{2} + 3t' + a' - 1$.

**Proof.** Each 2-circuit cycle moves at least four points, while each acceptable cycle moves at least two points. Thus, $(3t - 3Y) < (3t - 4Y)$, $(a - X) < (a - 2X)$.

**Theorem 12** Let $T$ be a tour and $C$ an acceptable cycle containing $\frac{n}{2}$ points that was obtained in $T^{-1}M^{-}$. Suppose |C| < $|T|$. Then $C$ yields a tour $T'$ such that $|T'| < |T|$.

**Proof**. By construction $C = (a_1\ T(a_2)\ a_2\ T(a_3)\ a_3\ \ldots\ a_{n-1}\ T(a_n)\ a_n\ T(a_1)) = T'$.

**Theorem 13** Let $S_{FW}$ be the set of all acceptable permutations of value no greater than $|T_{UPPERBOUND}|$ obtained from $\sigma_T^{-1}M^{-}$ using the modified F-W algorithm, while $S_{ALL}$ is the set of *all* acceptable cycles of value no greater than $|T_{UPPERBOUND}|$. Then the following holds:

Using the modified FW-algorithm, if we can obtain all elements of $S_{FW}$ ($S_{ALL}$) in polynomial time we can always obtain $\sigma_{FWTSPOPT}$ ($\sigma_{TSPOPT}$) in polynomial time.

**Proof**. Each acceptable permutation $P$ yields a perfect matching $\sigma_P$. Given each $P$, there exists a minimally-valued tour $T_{\sigma_P}$ of which it is a set of alternating edges. Since it is a minimally-valued tour, it can be obtained as the product $\sigma_P s$ where $s$ is an acceptable permutation containing precisely $\frac{n}{2}$ points obtainable from $\sigma_P^{-1}M^-$. This can always be done in polynomial time. On the other hand, $T^* = \sigma_{FWTSPOPT}$, has a value no greater than $|T_{UPPERBOUND}|$. Therefore, it contains a set of alternating edges of value no greater than $|T_{UPPERBOUND}|$ that is a perfect matching, $\sigma_{T^*}$. This perfect matching must be an element of $S_{FW}$. It follows that we can obtain $s$ in polynomial time by applying the modified F-W algorithm to $\sigma_{T^*}^{-1}M^-$.

An analogous proof holds for (2).

**Theorem 14** Let $\sigma$ be a product of $\frac{n}{2}$ disjoint 2-cycles that is obtained from a perfect matching. Suppose that $T$ is a tour. Let $s$ be a permutation such that $\sigma s = T$. Then $s$ contains at least one point in each 2-cycle of $T$ such that each cycle of $s$ is either an acceptable or 2-circuit cycle.

**Proof.** We consider alternate possibilities. First, $s$ must have at least one point in each 2-cycle in order to obtain a tour. Each arc, $(a\ b)$, in a cycle can represent a sequence of two edges $[a\ \sigma(b)][\ \sigma(b)\ b]$. In this way, each acceptable cycle of $j$ arcs represents a circuit containing $2j$ edges. Suppose a cycle is *not* acceptable. Using this procedure, there is only one way in which it can contribute cycles: Suppose the cycle is of one of the following form:

(1) $(a_1\ a_2\ \dots\ a_r\ b_1\ a_{r+1}\ a_{r+2}\ \dots\ a_s)$ where all of the $a_i's$ belong to distinct 2-cycle, while $b_1$ is the other point in the 2-cycle $(a_1\ b_1)$ of $\sigma$. Then we obtain

$$\begin{aligned}
&\left[a_1\ \underline{\sigma(a_2)}\ a_2\ \underline{\sigma(a_3)}\ \dots\ a_r\ \underline{\sigma(b_1)}\ b_1\ \underline{\sigma(a_{r+1})}\ a_{r+1}\ \underline{\sigma(a_{r+2})}\ a_{r+2}\ \underline{\sigma(a_{r+3})}\dots\ a_s\ \underline{\sigma(a_1)}\right]\\
&=\left[a_1\ \underline{\sigma(a_2)}\ \dots\ a_r\ \underline{a_1}\ b_1\ \underline{\sigma(a_{r+1})}\ a_{r+1}\ \underline{\sigma(a_{r+2})}\ a_{r+2}\ \underline{\sigma(a_{r+3})}\dots\ a_s\ \underline{b_1}\right]\\
&=\left[a_1\ \underline{\sigma(a_2)}\ \dots\ a_r\right]\left[b_1\ \underline{\sigma(a_{r+1})}\ a_{r+1}\ \underline{\sigma(a_{r+2})}\ a_{r+2}\ \underline{\sigma(a_{r+3})}\ \dots\ a_s\right]
\end{aligned}$$

Thus, we obtain two circuits. Furthermore, we cannot obtain *more than* two circuits. If the two circuits have no edge in common they are said to be unlinked. If they have precisely one edge in common, they are linked. If they are linked, they can become one circuit by deleting the common edge. Suppose we have two cases $(a_i\ b_i),\ i = 1,2$ where the points *are not interlaced* (say $a_1, b_2, a_2, b_1$)as they are with linked circuits. Then we have

$$\begin{aligned}
&\left[a_1\ \dots\ \underline{a_1}\, b_1\ \sigma(a_{r+1})\ a_{r+1}\ \dots\ \underline{a_{r+1}}\dots\ \underline{b_1}\ a_1\right]\\
&=\left[a_1\ \dots\ \underline{a_1}\right]\left[b_1\ \dots\ a_{r+1}\ \dots\ \underline{a_{r+1}}\ \dots\ \underline{b_1}\right]
\end{aligned}.$$

But $\left[b_1 \ \dots \ a_{r+1} \ \dots \ \underline{a_{r+1}} \ \dots \ \underline{b_1}\right]$ is *not* a circuit because we cannot express it as a sequence of pair-wise disjoint edges such that the first repetition of a point occurs only when we reach the initial point $b_1$.

*Note*. Using the determining vertices of cycles obtained using the modified F-W algorithm, we can employ rooted trees together with the end-points of the respective paths used to obtain cycles, to obtain *all* acceptable cycles. They are then added to the list of acceptable cycles obtained using the F-W algorithm. We would thus construct any further acceptable permutations to obtain $S_{ALL}$. Another approach is to obtain *all* paths using the modified F-W algorithm. This consists of the following: Every time we wish to replace a path by a shorter one by changing an entry in some $P_i$, we copy the original path from $P_i$. After we have obtained all acceptable cycles using the modified F-W algorithm, we delete those paths that don't have an initial vertex that is the determining vertex of an acceptable cycle. We then try to obtain an acceptable cycle from each of the remaining paths by using a rooted tree for each one. This saves us the trouble of backtracking at least some of the paths that eventually lead to acceptable cycles. A couple of final thoughts. It may make a difference if the tour is obtainable by patching together only acceptable cycles or a linked cycle is included. In the latter case, it may be more likely that our acceptable permutation in theorem 12 contains more than one cycle. Secondly, if a tour of smaller value exists, it may be more likely that the acceptable permutation contains only one cycle. Lastly, we might use our algorithm that requires obtaining not only acceptable but also unlinked and linked cycles. Using a larger cycle and a few smaller ones, we try to obtain a reasonably small-valued tour. Once we have such a tour, say $T_{large\ cycle}$, we obtain $\sigma_{T_{large\ cycle}}$. Using the modified F-W algorithm, we then obtain all acceptable permutations, $p_i$, of value less than $|\sigma_{T_{large\ cycle}}|$. Using each $p_i^{-1}M^-$, we try to obtain a still smaller tour by constructing an acceptable cycle containing precisely $\frac{n}{2}$ points.

To illustrate Phase 1, we give the following example. The symmetric cost matrix $M$ has been chosen randomly from positive integers no greater than 99.

**Example 3**

M

| | 1 | 2 | 3 | 4 | 5 | 6 | 7 | 8 | 9 | 10 | 11 | 12 | 13 | 14 | 15 | 16 | 17 | 18 | 19 | 20 | |
|---|---|---|---|---|---|---|---|---|---|---|---|---|---|---|---|---|---|---|---|---|---|
| 1 | ∞ | 26 | 4 | 30 | 74 | 5 | 4 | 38 | 28 | 78 | 81 | 7 | 97 | 10 | 94 | 40 | 98 | 49 | 40 | 70 | 1 |
| 2 | 26 | ∞ | 69 | 30 | 41 | 80 | 50 | 74 | 1 | 60 | 9 | 9 | 31 | 87 | 89 | 91 | 6 | 82 | 23 | 85 | 2 |
| 3 | 4 | 69 | ∞ | 23 | 7 | 61 | 60 | 98 | 99 | 90 | 84 | 57 | 4 | 56 | 66 | 30 | 51 | 3 | 25 | 47 | 3 |
| 4 | 30 | 30 | 23 | ∞ | 33 | 59 | 5 | 92 | 26 | 48 | 84 | 18 | 57 | 28 | 47 | 25 | 81 | 48 | 70 | 17 | 4 |
| 5 | 74 | 41 | 7 | 33 | ∞ | 82 | 29 | 80 | 5 | 87 | 87 | 97 | 55 | 45 | 72 | 94 | 20 | 9 | 90 | 20 | 5 |
| 6 | 5 | 80 | 61 | 59 | 82 | ∞ | 1 | 6 | 43 | 9 | 39 | 41 | 56 | 45 | 62 | 38 | 50 | 52 | 41 | 50 | 6 |
| 7 | 4 | 30 | 60 | 5 | 29 | 1 | ∞ | 34 | 78 | 49 | 73 | 10 | 56 | 36 | 87 | 31 | 45 | 59 | 88 | 42 | 7 |
| 8 | 38 | 74 | 98 | 92 | 80 | 6 | 34 | ∞ | 14 | 55 | 43 | 91 | 85 | 93 | 75 | 17 | 64 | 78 | 60 | 41 | 8 |
| 9 | 28 | 1 | 99 | 26 | 5 | 43 | 78 | 14 | ∞ | 50 | 28 | 81 | 98 | 95 | 25 | 31 | 73 | 63 | 87 | 56 | 9 |
| 10 | 78 | 60 | 90 | 48 | 87 | 9 | 49 | 55 | 50 | ∞ | 37 | 24 | 95 | 59 | 30 | 25 | 8 | 90 | 64 | 36 | 10 |
| 11 | 81 | 9 | 84 | 84 | 87 | 39 | 73 | 43 | 28 | 37 | ∞ | 34 | 61 | 14 | 11 | 3 | 3 | 74 | 22 | 26 | 11 |
| 12 | 7 | 9 | 57 | 18 | 97 | 41 | 10 | 91 | 81 | 24 | 34 | ∞ | 84 | 62 | 56 | 34 | 95 | 17 | 71 | 30 | 12 |
| 13 | 97 | 31 | 4 | 57 | 55 | 56 | 56 | 85 | 98 | 93 | 61 | 84 | ∞ | 66 | 66 | 30 | 49 | 82 | 23 | 86 | 13 |
| 14 | 10 | 87 | 56 | 28 | 45 | 45 | 36 | 93 | 95 | 59 | 14 | 62 | 66 | ∞ | 25 | 17 | 47 | 47 | 55 | 5 | 14 |
| 15 | 94 | 89 | 66 | 47 | 72 | 62 | 87 | 75 | 25 | 30 | 11 | 56 | 66 | 25 | ∞ | 89 | 85 | 87 | 8 | 2 | 15 |
| 16 | 40 | 91 | 30 | 25 | 94 | 38 | 31 | 17 | 31 | 25 | 3 | 34 | 30 | 17 | 89 | ∞ | 54 | 18 | 92 | 90 | 16 |
| 17 | 98 | 6 | 51 | 81 | 20 | 50 | 45 | 64 | 73 | 8 | 3 | 95 | 49 | 47 | 85 | 54 | ∞ | 66 | 48 | 90 | 17 |
| 18 | 49 | 82 | 3 | 48 | 9 | 52 | 59 | 78 | 63 | 90 | 74 | 17 | 82 | 47 | 87 | 18 | 66 | ∞ | 87 | 16 | 18 |
| 19 | 40 | 23 | 25 | 70 | 90 | 41 | 88 | 60 | 87 | 64 | 22 | 71 | 23 | 55 | 8 | 92 | 48 | 87 | ∞ | 22 | 19 |
| 20 | 70 | 85 | 47 | 17 | 20 | 50 | 42 | 41 | 56 | 36 | 26 | 30 | 86 | 5 | 2 | 90 | 90 | 16 | 22 | ∞ | 20 |
| | 1 | 2 | 3 | 4 | 5 | 6 | 7 | 8 | 9 | 10 | 11 | 12 | 13 | 14 | 15 | 16 | 17 | 18 | 19 | 20 | |

MIN(M)

| | 1 | 2 | 3 | 4 | 5 | 6 | 7 | 8 | 9 | 10 | 11 | 12 | 13 | 14 | 15 | 16 | 17 | 18 | 19 | 20 | |
|---|---|---|---|---|---|---|---|---|---|---|---|---|---|---|---|---|---|---|---|---|---|
| 1 | 3 | 7 | 6 | 12 | 14 | 2 | 9 | 4 | 8 | 16 | 19 | 18 | 20 | 5 | 10 | 11 | 15 | 13 | 17 | 1 | 1 |
| 2 | 9 | 17 | 11 | 12 | 19 | 1 | 4 | 7 | 13 | 5 | 10 | 3 | 8 | 6 | 20 | 18 | 14 | 15 | 16 | 2 | 2 |
| 3 | 18 | 13 | 1 | 5 | 4 | 19 | 16 | 20 | 17 | 14 | 12 | 7 | 6 | 15 | 2 | 11 | 10 | 8 | 9 | 3 | 3 |
| 4 | 7 | 20 | 12 | 3 | 16 | 9 | 14 | 1 | 2 | 5 | 15 | 10 | 18 | 13 | 6 | 19 | 17 | 11 | 8 | 4 | 4 |
| 5 | 9 | 3 | 18 | 17 | 20 | 7 | 4 | 2 | 14 | 13 | 15 | 1 | 8 | 6 | 10 | 11 | 19 | 16 | 12 | 5 | 5 |
| 6 | 7 | 1 | 8 | 10 | 16 | 11 | 12 | 19 | 9 | 14 | 17 | 20 | 18 | 13 | 4 | 3 | 15 | 2 | 5 | 6 | 6 |
| 7 | 6 | 1 | 4 | 12 | 5 | 2 | 16 | 8 | 14 | 20 | 17 | 10 | 13 | 18 | 3 | 11 | 9 | 15 | 19 | 17 | 7 |
| 8 | 6 | 9 | 16 | 7 | 1 | 20 | 11 | 10 | 19 | 17 | 2 | 15 | 18 | 5 | 13 | 12 | 4 | 14 | 3 | 8 | 8 |
| 9 | 2 | 5 | 8 | 15 | 4 | 1 | 11 | 16 | 6 | 10 | 20 | 18 | 17 | 7 | 12 | 19 | 14 | 13 | 3 | 9 | 9 |
| 10 | 17 | 6 | 12 | 16 | 15 | 20 | 11 | 4 | 7 | 9 | 8 | 14 | 2 | 19 | 1 | 5 | 3 | 18 | 13 | 10 | 10 |
| 11 | 16 | 17 | 2 | 15 | 14 | 19 | 20 | 9 | 12 | 10 | 6 | 8 | 13 | 7 | 18 | 1 | 3 | 4 | 5 | 11 | 11 |
| 12 | 1 | 2 | 7 | 18 | 4 | 10 | 20 | 11 | 16 | 6 | 15 | 3 | 14 | 19 | 9 | 13 | 8 | 17 | 5 | 12 | 12 |
| 13 | 3 | 19 | 16 | 2 | 17 | 5 | 6 | 7 | 4 | 11 | 14 | 15 | 18 | 12 | 8 | 20 | 10 | 1 | 9 | 13 | 13 |
| 14 | 20 | 1 | 11 | 16 | 15 | 4 | 7 | 5 | 6 | 17 | 18 | 19 | 3 | 10 | 12 | 13 | 2 | 8 | 9 | 14 | 14 |
| 15 | 20 | 19 | 11 | 9 | 14 | 10 | 4 | 12 | 6 | 3 | 13 | 5 | 8 | 17 | 7 | 18 | 2 | 16 | 1 | 15 | 15 |
| 16 | 11 | 8 | 14 | 18 | 4 | 10 | 13 | 3 | 7 | 9 | 12 | 6 | 1 | 17 | 15 | 20 | 2 | 19 | 5 | 16 | 16 |
| 17 | 11 | 2 | 10 | 5 | 7 | 14 | 19 | 13 | 6 | 3 | 16 | 8 | 18 | 9 | 4 | 15 | 20 | 12 | 1 | 17 | 17 |
| 18 | 3 | 5 | 20 | 12 | 16 | 14 | 4 | 1 | 6 | 7 | 9 | 17 | 11 | 8 | 2 | 13 | 15 | 19 | 10 | 18 | 18 |
| 19 | 15 | 11 | 20 | 2 | 13 | 3 | 1 | 6 | 17 | 14 | 8 | 10 | 4 | 12 | 9 | 18 | 7 | 5 | 16 | 19 | 19 |
| 20 | 15 | 14 | 18 | 4 | 5 | 19 | 11 | 12 | 10 | 8 | 7 | 3 | 6 | 9 | 1 | 2 | 13 | 16 | 17 | 20 | 20 |
| | 1 | 2 | 3 | 4 | 5 | 6 | 7 | 8 | 9 | 10 | 11 | 12 | 13 | 14 | 15 | 16 | 17 | 18 | 19 | 20 | |

$D_0 = (1^{26}\ 2^{69}\ 3^{23}\ 4^{63}\ 5^{82}\ 6^{1}\ 7^{34}\ 8^{14}\ 9^{50}\ 10^{37}\ 11^{34}\ 12^{84}\ 13^{66}\ 14^{25}\ 15^{89}\ 16^{54}\ 17^{66}\ 18^{87}\ 19^{22}\ 20^{20}\ )$

We now change $D_0$ into a form that can be used in our algorithm. In place of $a^{d(a\ D_0(a))}$, we put $a^{d(a\ MIN(M)(a,1))\ -\ d(a\ D_0(b))}$ yielding

$D_0 = (1^{-22}\ 2^{-68}\ 3^{-20}\ 4^{-26}\ 5^{-77}\ 6^{0}\ 7^{-33}\ 8^{-8}\ 9^{-49}\ 10^{-29}\ 11^{-31}\ 12^{-77}\ 13^{-62}\ 14^{-20}\ 15^{-87}\ 16^{-51}\ 17^{-63}\ 18^{-84}\ 19^{-14}\ 20^{-68}\ )$

Placing $D_0$ in row form, we obtain

| | | | | | | | | | | | | | | | | | | | | |
|---|---|---|---|---|---|---|---|---|---|---|---|---|---|---|---|---|---|---|---|---|
| | **-22** | **-68** | **-20** | **-26** | **–77** | **0** | **-33** | **-8** | **-49** | **-29** | **-31** | **-77** | **-62** | **-20** | **-87** | **-51** | **-63** | **-84** | **-14** | **-68** |
| | 1 | 2 | 3 | 4 | 5 | 6 | 7 | 8 | 9 | 10 | 11 | 12 | 13 | 14 | 15 | 16 | 17 | 18 | 19 | 20 |
| $D_0 =$ | | | | | | | | | | | | | | | | | | | | |
| | 2 | 3 | 4 | 5 | 6 | 7 | 8 | 8 | 10 | 11 | 12 | 13 | 14 | 15 | 16 | 17 | 18 | 19 | 20 | 1 |
| | 1 | 2 | 3 | 4 | 5 | 6 | 7 | 8 | 9 | 10 | 11 | 12 | 13 | 14 | 15 | 16 | 17 | 18 | 19 | 20 |
| $D_0^{-1}$ | | | | | | | | | | | | | | | | | | | | |
| | 20 | 19 | 18 | 17 | 16 | 15 | 14 | 13 | 12 | 11 | 10 | 9 | 8 | 7 | 6 | 5 | 4 | 3 | 2 | 1 |

Let the maximum number of blocks in a row used in a trial be $T = [log\ n] + 1$. Our first step is to sort all vertices in an ascending manner with respect to $d(a,\ D_i(a)) - d(a,\ MIN(M)(a,1))$. Let $v$ be the vertex with the smallest difference value. Our initial vertex in the first trial is $v$. We choose $(v\ MIN(M)(v,1))$ as our first arc. When we reach a vertex previously chosen, *if the vertex is the initial vertex of the cycle*, we conclude the trial using the cycle obtained. Otherwise, if each element of the block following it has a negative value, we can choose vertices in it until we reach a vertex that hasn't been chosen. If all of the vertices in the block have been previously chosen, we can try the next block. We may continue this procedure until we have gone through at most $T$ blocks. Each vertex chosen must have a negative value associated with it. If we reach an arc $A = (a\ D_i(a)),\ i = 0,1,2,\ldots$, we may choose any other vertex that lies in the same block. If no other vertex lies in that block, we may choose any vertex in the next block. If we still cannot find a vertex that has not been repeated, we check each possible cycle to see which one has the smallest value, say $s_i$. This completes one trial. The first arc chosen in the second trial is $(v\ MIN(M)(v,2))$. The maximum number of trials is $T$ Place the edges obtained from $s_i$ into the row form for $D_i$ in place of the edges with corresponding initial vertices. The new derangement obtained is $D_{i+1}$ Now sort the vertices of $D_{i+1}$ The first vertex obtained in the sorting is vertex used in the second trial. These comments roughly sketch the algorithm in the asymmetric case.

We now go over changes we must make to the algorithm used in the asymmetric case.

PHASE 1.

(1) Before choosing an arc $A = (a\ b)$ in a permutation to be applied to $D_i$, we check the row form of $D_i$ to see if the arc $(D_i(b)\ a)$ is in the $D_i(b)$ column. If it is, we can't use arc $A$ in the permutation.

(2) During each trial, as we place an arc $(a\ b)$ in a permutation, we place both the arc $(a\ D_i(b))$ and the arc $(D_i(b)\ a)$ in a list sorted in increasing order of magnitude of the initial integer of the arc. Before a new arc $(a\ b)$ is added to a permutation, we check to see if $(D_i(b)\ a)$ is in the list. If it is, we can't place $(a\ b)$ in the permutation.

TRIAL 1. MIN(M)(15,1) = 20.

$(15\ \ 20) \rightarrow (15\ \ 19)$ **-87** **1.** **11.**

| | | |
|---|---|---|
| (19 15) → (19 14) **-14** | **2. (2 9)** | **12.** |
| (14 20) **repeat** | **3. (3 16)** | **13. (13 3)** |
| (14 1) → (14 20) **-5** | **4.** | **14. (14 1)** |
| (20 15) **not allowed** | **5. (5 18)** | **15. (15 20)** |
| (20 14) → (20 13) **-65** | **6.** | **16. (16 8)** |
| (13 3) → (13 2) **-62** | **7. (7 4)** | **17. (17 11)** |
| (2 9) → (2 8) **-68** | **8. (8 6)** | **18.** |
| (8 6) → (8 5) **-8** | **9.** | **19. (19 15)** |
| (5 9) **repeat** | **10. (10 17)** | **20. (20 14)** |

(5 3) **repeat**
(5 18) → (5 17) **-73**
(17 11) → (17 10) **-63**
(10 17) → (10 16) **-29**
(16 10) **repeat**
(16 8) → (16 7) **-14**
(7 6) **repeat**
(7 1) **repeat**
(7 4) → (7 3) **-29**
(3 13) **not allowed**
(3 6) **repeat**
(3 5) **repeat**

We can't proceed further since *[log 20] + 1 = 3* . We thus must form a cycle using (3 16).

(3 16) → (3 15) **+7**

Since all of the differences are negative except the last one, and 7 is small, **the permutation with the smallest difference value is (15 19 14 20 13 2 8 5 17 10 16 7 3): -516**

TRIAL 2. MIN(15, 2) = 19

| | | |
|---|---|---|
| (15 19) → (15 18) **-79** | **1.** | **11..** |
| (18 3) → (18 2) **-84** | **2. (2 9)** | **12.** |
| (2 8)(2 9) → (2 8) **-68** | **3.** | **13.** |
| (8 6) → (8 5) **-8** | **4.** | **14. (14 16)** |
| (5 18) → (5 17) **-62** | **5. (5 18)** | **15. (15 19)** |
| (10 17) → (10 16) **-29** | **6.** | **16. (16 8)** |
| (16 8) → (16 7) **-37** | **7. (7 1)** | **17. (17 11)** |
| (7 1) → (7 20) **-30** | **8.** | **18. (18 3)** |
| (20 15) → (20 14) **-68** | **9.** | **19.** |
| (14 16) → (14 15) **-8** | **10. (10 17)** | **20. (20 15)** |

**Our cycle is (15 18 2 8 5 17 10 16 7 20 14): -536**

| | | |
|---|---|---|
| TRIAL 3. MIN(M)(15, 3) = 11 | **1..** | **11.** |
| (15 11) → (15 10) **-78** | **2.** | **12.** |
| (10 17) $\rightarrow$ (10 16) **-29** | **3.** | **13.** |
| (16 11) **repeat** | **4.** | **14.** |
| (16 8) $\rightarrow$ (16 7) **-14** | **5. (5 9)** | **15. (15 11)** |
| (7 6) $\rightarrow$ (7 5) **-33** | **6.** | **16. (16 8)** |
| (5 9) $\rightarrow$ (5 8) **-77** | **7. (7 6)** | **17.** |
| (8 6) **repeat** | **8. (8 16)** | **18.** |
| (8 9) **arc of $D_0$** | **9.** | **19.** |
| (8 16) $\rightarrow$ (18 15) +**3** | **10. (10 17)** | **20.** |

**Our cycle is (15 10 16 7 5 8): -228**

$s$ is the cycle obtained in Trial 2: **(15 18 2 8 5 17 10 16 7 20 14).** Our new edges are (15 19), (18 3), (2 9), (8 6), (5 18), (17 11), (10 17), (7 1), (20 15). $D_1 = D_0 s$ . We thus obtain:

| | | | | | | | | | | | | | | | | | | | | |
|---|---|---|---|---|---|---|---|---|---|---|---|---|---|---|---|---|---|---|---|---|
| | **-22** | **0** | **-20** | **-28** | **-4** | **0** | **-3** | **0** | **-49** | **0** | **-31** | **-77** | **-62** | **-12** | **-6** | **-14** | **0** | **0** | **-14** | **0** |
| | 1 | 2 | 3 | 4 | 5 | 6 | 7 | 8 | 9 | 10 | 11 | 12 | 13 | 14 | 15 | 16 | 17 | 18 | 19 | 20 |
| $D_1$ = | | | | | | | | | | | | | | | | | | | | |
| | 2 | 9 | 4 | 5 | 18 | 7 | 1 | 6 | 10 | 17 | 12 | 13 | 14 | 16 | 19 | 8 | 11 | 3 | 20 | 15 |
| | 1 | 2 | 3 | 4 | 5 | 6 | 7 | 8 | 9 | 10 | 11 | 12 | 13 | 14 | 15 | 16 | 17 | 18 | 19 | 20 |
| $D_1^{-1}$ = | | | | | | | | | | | | | | | | | | | | |
| | 7 | 1 | 18 | 3 | 4 | 8 | 6 | 18 | 2 | 9 | 17 | 11 | 12 | 13 | 20 | 14 | 10 | 5 | 15 | 19 |

Choose 12.

| | | |
|---|---|---|
| TRIAL 1. MIN(M)(12,1) = 1 | **1.** | **11.** |
| (12 1) $\rightarrow$ (12 7) **-77** | **2.** | **12. (12 1)** |
| (7 6) **not allowed** | **3. (3 13)** | **13.** |
| (7 1) **arc repeat** | **4.** | **14.** |
| (7 4) $\rightarrow$ (7 3) **+1** | **5.** | **15.** |
| (3 18) **not allowed** | **6.** | **16.** |
| (3 13) $\rightarrow$ (3 12) **-19** | **7. (7 4)** | **17.** |
| | **8.** | **18.** |
| | **9.** | **19.** |
| . | **10.** | **20.** |

Our cycle is **(12 7 3): -95.**

| | | |
|---|---|---|
| Trial 2. MIN(M)(12,2) = 2 | **1. (1 3)** | **11.** |
| (12 2) $\rightarrow$ (12 1) **-78** | **2.** | **12. (12 2)** |
| (1 3) $\rightarrow$ (1 18) **-22** | **3.** | **13.** |
| (18 3) **arc** | **4.** | **14.** |
| (18 5) **not allowed** | **5.** | **15.** |
| (18 20) → (18 19) **+13** | **6.** | **16.** |
| (19 15) **not allowed** | **7.** | **17.** |
| (19 11) → (19 17) **0** | **8.** | **18. (18 20)** |
| (17 11) **arc repeat** | **9.** | **19. (19 13)** |
| (17 2) **repeat** | **10.** | **20.** |

We can't go further, d(17 13) - d(17 11) = +46. d(19 13) - d(19 20) = 23 - 22 = 1.

d(18 13) - d(18 3) = 82 - 3 = 79. d(1 13) - d(1 2) = 97 - 26 = 51. Thus, our smallest valued cycle is

| | |
|---|---|
| **(12 1 18 19): -86**: | **1**. |
| Trial 3. MIN(M)(12,3) = 7 | **2**. |
| (12 7) → (12 6) **-74** | **3. (3 13)** |
| (6 7) **arc repeat** | **4.** |
| (6 1) → (6 7) **+4** | **5**. |
| (7 6) **not allowed** | **6.** |
| (7 1) **arc** | **7. (7 4)** |
| (7 4) → (7 3) **+1** | **8**.. |
| (3 18) **not allowed** | **9.** |
| (3 13) → (3 12) **-19** | **10**. |

Our smallest-valued cycle is obtained from Trial 1: **(12 7 3): -95.** Our new edges are (12 1), (7 4), (3 13).:

| | | | | | | | | | | | | | | | | | | | |
|---|---|---|---|---|---|---|---|---|---|---|---|---|---|---|---|---|---|---|---|
| **-22** | **0** | **-1** | **-28** | **-4** | **0** | **-4** | **0** | **-49** | **0** | **-31** | **-3** | **-62** | **-12** | **-6** | **-14** | **0** | **0** | **-14** | **0** |
| 1 | 2 | 3 | 4 | 5 | 6 | 7 | 8 | 9 | 10 | 11 | 12 | 13 | 14 | 15 | 16 | 17 | 18 | 19 | 20 |

$D_2$ =

| | | | | | | | | | | | | | | | | | | | |
|---|---|---|---|---|---|---|---|---|---|---|---|---|---|---|---|---|---|---|---|
| 2 | 9 | 13 | 5 | 18 | 7 | 4 | 6 | 10 | 17 | 12 | 1 | 14 | 16 | 19 | 8 | 11 | 3 | 20 | 15 |
| 1 | 2 | 3 | 4 | 5 | 6 | 7 | 8 | 9 | 10 | 11 | 12 | 13 | 14 | 15 | 16 | 17 | 18 | 19 | 20 |

$D_2^{-1}$

| | | | | | | | | | | | | | | | | | | | |
|---|---|---|---|---|---|---|---|---|---|---|---|---|---|---|---|---|---|---|---|
| 12 | 1 | 18 | 7 | 4 | 8 | 6 | 16 | 2 | 9 | 17 | 11 | 3 | 13 | 20 | 14 | 10 | 5 | 15 | 19 |

Choose 13.

| | | |
|---|---|---|
| Trial 1. MIN(M)(13,1) = 3 | **1. (1 3)** | **11.** |
| (13 3) **not allowed** | **2**. | **12.** |
| Trial 2. MIN(M)(13,2) = 19 | **3**. | **13. (13 19)** |
| (13 19) $\rightarrow$ (13 15) **-43** | **4.** | **14.** |
| (15 20) **not allowed** | **5**. | **15. (15 11)** |
| (15 19) **arc** | **6.** | **16.** |
| (15 11) $\rightarrow$ (15 17) **+3** | **7.** | **17. (17 2)** |
| (17 11) **repeat** | **8.** | **18. (18 20)** |
| (17 2) $\rightarrow$ (17 1) **+3** | **9.** | **19.** |
| (1 3) $\rightarrow$ (1 18) **-22** | **10.** | **20.** |

(18 3) **arc**
(18 5) **not allowed**
(18 20) $\rightarrow$ (18 19) **+13**
(19 15) **not allowed**
(19 11) **repeat**
(19 20) **arc**
(19 2) **repeat positive**

We can go no further d(19 14) - d(19 11) = +33. d(18 14) - d(18 3) = +44. . d(1 14) - d(1 2) = -16.. d(17 14) - d(17 11) = +44. d(15 14) - d(15 19) = +17. Our smallest-valued cycle is **(13 15 17 1): -62**

| | | |
|---|---|---|
| Trial 3. MIN(M)(13,3) = 16 | **1. (1 3)** | **11**. |
| (13 16) $\rightarrow$ (13 14) **-36** | **2**. | **12**. |
| (14 20) $\rightarrow$ (14 19) **-12** | **3.** | **13. (13 15)** |
| (19 15) **not allowed** | **4.** | **14. (14 20)** |
| (19 11) $\rightarrow$ (19 17) **0** | **5.** | **15.** |
| (17 11) **arc** | **6.** | **16.** |
| (17 2) $\rightarrow$ (17 1) **+3** | **7.** | **17. (17 2)** |
| (1 3) $\rightarrow$ (1 18) **-22** | **8.** | **18.** |
| (18 3) **arc** | **9.** | **19. (19 11)** |
| (18 5) **not allowed** | **10.** | **20.** |

(18 20) **repeat positive**

We can go no further. d(18 14) - d(18 3) = 44. d(1 14) - d(1 2) = -16. d(17 14) - d(17 11) = 44. d(19 14) - d(19 20) = 55 - 22 = 33. Our best cycle is **(13 14 19 17 1): -61**

Our smallest-valued cycle occurs in Trial 3: **(13 15 17 1): -62.** Our new edges are: (13 19), (15 11), (17 2), (1 14).

| | | | | | | | | | | | | | | | | | | | |
|---|---|---|---|---|---|---|---|---|---|---|---|---|---|---|---|---|---|---|---|
| **-6** | **0** | **-1** | **-28** | **-4** | **0** | **-4** | **0** | **-49** | **0** | **-31** | **0** | **-19** | **-12** | **-9** | **-14** | **-3** | **0** | **-14** | **0** |
| 1 | 2 | 3 | 4 | 5 | 6 | 7 | 8 | 9 | 10 | 11 | 12 | 13 | 14 | 15 | 16 | 17 | 18 | 19 | 20 |

$D_3$ =

| | | | | | | | | | | | | | | | | | | | |
|---|---|---|---|---|---|---|---|---|---|---|---|---|---|---|---|---|---|---|---|
| 14 | 9 | 13 | 5 | 18 | 7 | 4 | 6 | 10 | 17 | 12 | 1 | 19 | 16 | 11 | 8 | 2 | 3 | 20 | 15 |
| 1 | 2 | 3 | 4 | 5 | 6 | 7 | 8 | 9 | 10 | 11 | 12 | 13 | 14 | 15 | 16 | 17 | 18 | 19 | 20 |

$D_3^{-1}$

| | | | | | | | | | | | | | | | | | | | |
|---|---|---|---|---|---|---|---|---|---|---|---|---|---|---|---|---|---|---|---|
| 12 | 17 | 18 | 7 | 4 | 8 | 6 | 16 | 2 | 9 | 15 | 11 | 3 | 1 | 20 | 14 | 10 | 5 | 13 | 19 |

Choose 9.

TRIAL 1. MIN(M)(9,1) = 2.

(9 2) **not allowed**

TRIAL 2. MIN(M)(9 2) = 5

(9 5) → (9 4) **-45**

(4 7) **not allowed**

(4 20) → (4 19) **-16**

(19 15) → (19 20) **-14**

(20 15) **arc repeat**

(20 14) → (20 1) **+3**

(1 3) → (1 18) **-6**

(18 3) **arc**

(18 5) **not allowed**

(18 20) **repeat positive**

We can go no further.

d(20, 10) - d(20 15) = **+34**

Our best cycle is **(9 4 19 20): -41**

| | |
|---|---|
| **1. (1 3)** | **11.** |
| **2.** | **12.** |
| **3.** | **13.** |
| **4. (4 20)** | **14.** |
| **5.** | **15.** |
| **6.** | **16.** |
| **7.** | **17.** |
| **8.** | **18.** |
| **9. (9 5)** | **19. (19 15)** |
| **10.** | **20. (20 14)** |

TRIAL 3. MIN(M)(9,3) = 8.

(9 8) → (9 16) **-35**

(16 11) → (16 15) **-14**

(15 20) **not allowed**

(15 19) → (15 13) **-3**

(13 3) **not allowed**

(13 19) **arc repeat**

(13 16) → (13 14) **+7**

(14 20) **repeat**

(14 1) **not allowed**

(14 11) **repeat positive**

| | |
|---|---|
| **1.** | **11.** |
| **2.** | **12.** |
| **3.** | **13. (13 16)** |
| **4.** | **14.** |
| **5.** | **15. (15 19)** |
| **6.** | **16. (16 11)** |
| **7.** | **17.** |
| **8.** | **18.** |
| **9. (9 8)** | **19.** |
| **10.** | **20.** |

We can go no further**.** d(15 10) - d(15 11) = +19. Our smallest-valued cycle is **(9 16 15): -30.**

It follows that $s$ is obtained from TRIAL 2: **(9 4 19 20): -41.** Our new edges are (9 5), (4 20), (19 15), (20 10) .

| | | | | | | | | | | | | | | | | | | | |
|---|---|---|---|---|---|---|---|---|---|---|---|---|---|---|---|---|---|---|---|
| **-6** | **0** | **-1** | **-12** | **-4** | **0** | **-4** | **0** | **-4** | **0** | **-31** | **0** | **-19** | **-12** | **-9** | **-14** | **-3** | **0** | **0** | **-34** |

$D_4$ =

| 1 | 2 | 3 | 4 | 5 | 6 | 7 | 8 | 9 | 10 | 11 | 12 | 13 | 14 | 15 | 16 | 17 | 18 | 19 | 20 |
|---|---|---|---|---|---|---|---|---|---|---|---|---|---|---|---|---|---|---|---|
| 14 | 9 | 13 | 20 | 18 | 7 | 4 | 6 | 5 | 17 | 12 | 1 | 19 | 16 | 11 | 8 | 2 | 3 | 15 | 10 |

$D_4^{-1}$ =

| 1 | 2 | 3 | 4 | 5 | 6 | 7 | 8 | 9 | 10 | 11 | 12 | 13 | 14 | 15 | 16 | 17 | 18 | 19 | 20 |
|---|---|---|---|---|---|---|---|---|---|---|---|---|---|---|---|---|---|---|---|
| 12 | 17 | 18 | 7 | 9 | 8 | 6 | 16 | 2 | 20 | 15 | 11 | 3 | 1 | 19 | 14 | 10 | 5 | 13 | 4 |

Choose 20.

TRIAL 1. MIN(M)(20,1) = 15

(20 15) $\rightarrow$ (20 19) **-34**

(19 15) **not allowed**

(19 11) $\rightarrow$ (19 15) **+14**

(15 20) **not allowed**

(15 19) **not allowed**

(15 11) **arc repeat**

Can't go further. Have gone through *[log n] + 1* columns of MIN((M). Furthermore, path is no longer negative. No negative cycles can be obtained from any subpath.

TRIAL 2. MIN(M)(20,2) = 14

(20 14) $\rightarrow$ (20 1) **-31**

(1 3) $\rightarrow$ (1 18) **-6**

(18 3) **arc**

(18 5) **not allowed**

(18 20) $\rightarrow$ (18 4) **+13**

(4 7) **not allowed**

(4 20) **arc repeat**

(4 12) $\rightarrow$ (4 11) **+1**

(11 16) $\rightarrow$ (11 14) **-31**

(14 20) **not allowed**

(14 1) **not allowed**

(14 11) $\rightarrow$ (14 15) **-3**

(15 20) **repeat**

(15  19)  **not allowed**

(15  11)  **arc repeat**

Can't go further. d(15  10)  -  d(15  11)  =  **+19**

Our smallest-valued cycle in this tiral is **(20  1  18  4  11  14  15):  -38**

TRIAL 3.  MIN(M)(20,3) = 18.

(20  18)  → (20  5)  **-20**

(5  9)  **not allowed**

(5  3)  → (5  18)  **-2**

(18  3)  **arc**

(18  5)  **not allowed**

(18  20)  → (18  4)  **+13**

(4  7)  **not allowed**

(4  20)  **arc repeat**

(4  12)  → (4  11)  **+1**

(11  16)  → (11  14)  **-31**

(14  20)  **repeat**

(14  1)  **not allowed**

(14  11)  → (14  15)  **-3**

(15  20)  **repeat**

(15  19)  → **not allowed**

(15  11)  **arc**

Can't go further. The smallest-valued cycle obtainable in this trial is **(20  5  18  4  11  14  15):  -23**

It follows that our smallest-valued cycle is **(20  1  18  4  11  14  15):  -38.** Our new edges are: (20  14),  (1  3), (18  20),  (4  12),  (11  16),  (14  11),  (15  10).

| | | | | | | | | | | | | | | | | | | | |
|---|---|---|---|---|---|---|---|---|---|---|---|---|---|---|---|---|---|---|---|
| **0** | **0** | **-1** | **-13** | **-4** | **0** | **-4** | **0** | **-4** | **0** | **0** | **0** | **-19** | **-9** | **-28** | **-14** | **-3** | **-13** | **0** | **-3** |
| 1 | 2 | 3 | 4 | 5 | 6 | 7 | 8 | 9 | 10 | 11 | 12 | 13 | 14 | 15 | 16 | 17 | 18 | 19 | 20 |

$D_5$ =

| | | | | | | | | | | | | | | | | | | | |
|---|---|---|---|---|---|---|---|---|---|---|---|---|---|---|---|---|---|---|---|
| 3 | 9 | 13 | 12 | 18 | 7 | 4 | 6 | 5 | 17 | 16 | 1 | 19 | 11 | 10 | 8 | 2 | 20 | 15 | 14 |
| 1 | 2 | 3 | 4 | 5 | 6 | 7 | 8 | 9 | 10 | 11 | 12 | 13 | 14 | 15 | 16 | 17 | 18 | 19 | 20 |

$D_5^{-1}$ =

| | | | | | | | | | | | | | | | | | | | |
|---|---|---|---|---|---|---|---|---|---|---|---|---|---|---|---|---|---|---|---|
| 12 | 17 | 1 | 7 | 9 | 8 | 6 | 16 | 2 | 15 | 14 | 4 | 3 | 20 | 19 | 14 | 10 | 5 | 13 | 18 |

Choose 15.

TRIAL 1. MIN(M)(15,1) = 20

(15 20) → (15 18) **-28**

(18 3) → (18 1) **-13**

(1 3) **arc**

(1 7) → (1 6) **0**

(6 7) **arc repeat**

(6 20) **repeat**

(6 8) **not allowed**

Can't go further. d(6 10) - d(6 7) = **+8**

Our smallest-valued negative cycle is **(15 18 1 7 6): -33**.

TRIAL 2. MIN(M)(15,2) = 19

(15 19) → (15 13) **-22**

(13 3) **not allowed**

(13 19) **arc repeat**

(13 16) → (13 11) **+7**

(11 16) **arc**

(11 17) → (11 10) **0**

(10 17) **arc**

(10 6) → (10 8) **+1**

(8 6) **arc repeat**

(8 9) → (8 2) **+8**

(2 9) **arc**

(2 17) **repeat positive**

Can't go further.

No negative cycle exists.

Our smallest-valued negative cycle occurs in Trial 1: **(15 18 1 6): -33.**

New edges are: (15 20), (18 3), (1 7), (6 10).

| | | | | | | | | | | | | | | | | | | | | |
|---|---|---|---|---|---|---|---|---|---|---|---|---|---|---|---|---|---|---|---|---|
| | **0** | **0** | **-1** | **-13** | **-4** | **-8** | **-4** | **0** | **-4** | **0** | **0** | **0** | **-19** | **-9** | **0** | **-14** | **-3** | **0** | **0** | **-3** |
| | 1 | 2 | 3 | 4 | 5 | 6 | 7 | 8 | 9 | 10 | 11 | 12 | 13 | 14 | 15 | 16 | 17 | 18 | 19 | 20 |
| $D_6$ = | | | | | | | | | | | | | | | | | | | | |
| | 7 | 9 | 13 | 12 | 18 | 10 | 4 | 6 | 5 | 17 | 16 | 1 | 19 | 11 | 20 | 8 | 2 | 3 | 15 | 14 |
| | 1 | 2 | 3 | 4 | 5 | 6 | 7 | 8 | 9 | 10 | 11 | 12 | 13 | 14 | 15 | 16 | 17 | 18 | 19 | 20 |
| $D_6^{-1}$ = | | | | | | | | | | | | | | | | | | | | |
| | 12 | 17 | 18 | 7 | 9 | 8 | 1 | 16 | 2 | 6 | 14 | 4 | 3 | 20 | 19 | 11 | 10 | 5 | 13 | 15 |

**(1 7 4 12)(2 9 5 18 3 13 19 15 20 14 11 16 8 6 10 17): 157**

Choose 13.

Trial 1. MIN(M)(13,1) = 3

(13 3) **not allowed**

(13 19) **arc**

Can't go further. Thus, using Phase 1, we cannot obtain a derangement containing $n$ edges that has a value less than that of $D_6$.

**Example 4.**

In this example, we chose a non-random matrix 3-Cycle M containing a derangement consisting of five cycles. That derangement is (1 3 7)(2 5 11 14 19)(4 16 8)(6 18 13)(10 17 12). We assume that this derangement was obtained at the end of Phase 1. In Phase 2, we don't use arcs that are symmetric to arcs in the derangement. We thus always obtain a derangement. We now give Phases 2 and 3.

3-cycle M

| | 1 | 2 | 3 | 4 | 5 | 6 | 7 | 8 | 9 | 10 | 11 | 12 | 13 | 14 | 15 | 16 | 17 | 18 | 19 | 20 | |
|---|---|---|---|---|---|---|---|---|---|---|---|---|---|---|---|---|---|---|---|---|---|
| 1 | ∞ | 26 | 1 | 30 | 74 | 5 | 3 | 38 | 28 | 78 | 81 | 7 | 97 | 10 | 94 | 40 | 98 | 49 | 40 | 70 | 1 |
| 2 | 26 | ∞ | 69 | 30 | 1 | 80 | 50 | 74 | 1 | 60 | 3 | 9 | 31 | 87 | 89 | 91 | 6 | 82 | 23 | 85 | 2 |
| 3 | 1 | 69 | ∞ | 23 | 7 | 61 | 2 | 98 | 99 | 90 | 84 | 57 | 4 | 56 | 66 | 30 | 51 | 3 | 25 | 47 | 3 |
| 4 | 30 | 30 | 23 | ∞ | 33 | 59 | 5 | 3 | 26 | 48 | 84 | 18 | 57 | 28 | 47 | 1 | 81 | 48 | 70 | 17 | 4 |
| 5 | 74 | 1 | 7 | 33 | ∞ | 82 | 29 | 80 | 5 | 87 | 2 | 97 | 3 | 45 | 72 | 94 | 20 | 9 | 90 | 20 | 5 |
| 6 | 5 | 80 | 61 | 59 | 82 | ∞ | 1 | 6 | 43 | 9 | 39 | 41 | 3 | 45 | 62 | 38 | 50 | 1 | 41 | 50 | 6 |
| 7 | 3 | 30 | 2 | 5 | 29 | 1 | ∞ | 34 | 78 | 49 | 73 | 10 | 56 | 36 | 87 | 31 | 45 | 59 | 88 | 42 | 7 |
| 8 | 38 | 74 | 98 | 3 | 80 | 6 | 34 | ∞ | 14 | 55 | 43 | 91 | 85 | 93 | 75 | 2 | 64 | 78 | 60 | 1 | 8 |
| 9 | 28 | 1 | 99 | 26 | 5 | 43 | 78 | 14 | ∞ | 50 | 28 | 81 | 98 | 95 | 3 | 31 | 73 | 63 | 87 | 2 | 9 |
| 10 | 78 | 60 | 90 | 48 | 87 | 9 | 49 | 55 | 50 | ∞ | 37 | 1 | 95 | 59 | 30 | 25 | 3 | 90 | 64 | 36 | 10 |
| 11 | 81 | 3 | 84 | 84 | 2 | 39 | 73 | 43 | 28 | 37 | ∞ | 3 | 61 | 14 | 11 | 3 | 3 | 74 | 22 | 26 | 11 |
| 12 | 7 | 9 | 57 | 18 | 97 | 41 | 10 | 91 | 81 | 1 | 2 | ∞ | 84 | 62 | 56 | 34 | 2 | 17 | 71 | 30 | 12 |
| 13 | 97 | 31 | 4 | 57 | 5 | 3 | 56 | 85 | 98 | 93 | 61 | 84 | ∞ | 66 | 66 | 30 | 49 | 2 | 23 | 86 | 13 |
| 14 | 10 | 87 | 56 | 28 | 45 | 45 | 36 | 93 | 95 | 59 | 14 | 62 | 66 | ∞ | 25 | 14 | 47 | 47 | 2 | 5 | 14 |
| 15 | 94 | 89 | 66 | 47 | 72 | 62 | 87 | 75 | 3 | 30 | 11 | 56 | 66 | 25 | ∞ | 89 | 85 | 87 | 8 | 2 | 15 |
| 16 | 40 | 91 | 30 | 1 | 94 | 38 | 31 | 2 | 31 | 25 | 3 | 34 | 30 | 14 | 89 | ∞ | 54 | 18 | 92 | 90 | 16 |
| 17 | 98 | 6 | 51 | 81 | 20 | 50 | 45 | 64 | 73 | 3 | 3 | 2 | 49 | 47 | 85 | 54 | ∞ | 66 | 48 | 90 | 17 |
| 18 | 49 | 82 | 3 | 48 | 9 | 1 | 59 | 78 | 63 | 90 | 74 | 17 | 2 | 47 | 87 | 18 | 66 | ∞ | 87 | 16 | 18 |
| 19 | 40 | 23 | 25 | 70 | 90 | 41 | 88 | 60 | 87 | 64 | 22 | 71 | 23 | 2 | 8 | 92 | 48 | 87 | ∞ | 22 | 19 |
| 20 | 70 | 85 | 47 | 17 | 20 | 50 | 42 | 1 | 2 | 36 | 26 | 30 | 86 | 5 | 2 | 90 | 90 | 16 | 22 | ∞ | 20 |
| | 1 | 2 | 3 | 4 | 5 | 6 | 7 | 8 | 9 | 10 | 11 | 12 | 13 | 14 | 15 | 16 | 17 | 18 | 19 | 20 | |

3-cycle MIN(M)

| | 1 | 2 | 3 | 4 | 5 | 6 | 7 | 8 | 9 | 10 | 11 | 12 | 13 | 14 | 15 | 16 | 17 | 18 | 19 | 20 | |
|---|---|---|---|---|---|---|---|---|---|---|---|---|---|---|---|---|---|---|---|---|---|
| 1 | 3 | 7 | 6 | 12 | 14 | 2 | 9 | 4 | 8 | 16 | 19 | 18 | 20 | 5 | 10 | 11 | 15 | 13 | 17 | 1 | 1 |
| 2 | 9 | 5 | 11 | 17 | 12 | 19 | 1 | 4 | 13 | 7 | 10 | 3 | 8 | 6 | 20 | 18 | 14 | 15 | 16 | 2 | 2 |
| 3 | 1 | 7 | 18 | 13 | 5 | 4 | 19 | 16 | 20 | 17 | 14 | 12 | 6 | 15 | 2 | 11 | 10 | 8 | 9 | 3 | |
| 4 | 16 | 8 | 7 | 20 | 12 | 3 | 9 | 14 | 1 | 2 | 5 | 15 | 10 | 18 | 13 | 6 | 19 | 17 | 11 | 4 | 4 |
| 5 | 2 | 11 | 13 | 9 | 3 | 18 | 17 | 20 | 7 | 4 | 14 | 15 | 1 | 8 | 6 | 10 | 19 | 16 | 12 | 5 | |
| 6 | 7 | 18 | 13 | 1 | 8 | 10 | 16 | 11 | 19 | 9 | 14 | 17 | 20 | 12 | 4 | 3 | 15 | 2 | 5 | 6 | 6 |
| 7 | 6 | 3 | 1 | 4 | 12 | 5 | 2 | 16 | 8 | 14 | 20 | 17 | 10 | 18 | 3 | 11 | 9 | 15 | 19 | 17 | 7 |
| 8 | 20 | 16 | 4 | 6 | 9 | 7 | 1 | 11 | 10 | 19 | 17 | 2 | 15 | 18 | 5 | 13 | 12 | 14 | 3 | 8 | 8 |
| 9 | 2 | 29 | 15 | 5 | 8 | 4 | 1 | 11 | 16 | 6 | 10 | 18 | 17 | 7 | 12 | 19 | 14 | 13 | 3 | 9 | 9 |
| 10 | 12 | 17 | 6 | 16 | 15 | 20 | 11 | 4 | 7 | 9 | 8 | 14 | 2 | 19 | 1 | 5 | 3 | 18 | 13 | 10 | 10 |
| 11 | 5 | 2 | 12 | 16 | 17 | 15 | 14 | 19 | 20 | 9 | 10 | 6 | 8 | 13 | 7 | 19 | 1 | 3 | 4 | 11 | 11 |
| 12 | 10 | 11 | 17 | 1 | 2 | 7 | 18 | 4 | 20 | 16 | 6 | 15 | 3 | 14 | 19 | 9 | 13 | 8 | 5 | 12 | 12 |
| 13 | 18 | 6 | 3 | 5 | 19 | 16 | 2 | 17 | 7 | 4 | 11 | 14 | 15 | 8 | 12 | 20 | 10 | 1 | 9 | 13 | 13 |
| 14 | 19 | 20 | 1 | 11 | 16 | 15 | 4 | 7 | 5 | 6 | 17 | 18 | 3 | 10 | 12 | 13 | 2 | 8 | 9 | 14 | 14 |
| 15 | 20 | 9 | 19 | 11 | 14 | 10 | 4 | 12 | 6 | 3 | 13 | 5 | 8 | 17 | 7 | 18 | 2 | 16 | 1 | 15 | 15 |
| 16 | 8 | 11 | 14 | 18 | 4 | 10 | 13 | 3 | 7 | 9 | 12 | 6 | 1 | 17 | 15 | 20 | 2 | 19 | 5 | 16 | 16 |
| 17 | 12 | 10 | 11 | 2 | 5 | 7 | 14 | 19 | 13 | 6 | 3 | 16 | 8 | 18 | 9 | 4 | 15 | 20 | 1 | 17 | 17 |
| 18 | 6 | 13 | 3 | 5 | 20 | 12 | 16 | 14 | 4 | 1 | 7 | 9 | 17 | 11 | 8 | 2 | 15 | 19 | 10 | 18 | 18 |
| 19 | 14 | 15 | 11 | 20 | 2 | 13 | 3 | 1 | 6 | 17 | 8 | 10 | 4 | 12 | 9 | 18 | 7 | 5 | 16 | 19 | 19 |
| 20 | 15 | 14 | 18 | 4 | 5 | 19 | 11 | 12 | 10 | 8 | 7 | 3 | 6 | 9 | 1 | 2 | 13 | 16 | 17 | 20 | 20 |
| | 1 | 2 | 3 | 4 | 5 | 6 | 7 | 8 | 9 | 10 | 11 | 12 | 13 | 14 | 15 | 16 | 17 | 18 | 19 | 20 | |

$D^{-1}_{3\text{-}cycle\,2}M^{-}$

| | | | | | | | | | | | | | | | | | | | | | |
|---|---|---|---|---|---|---|---|---|---|---|---|---|---|---|---|---|---|---|---|---|---|
| | 3 | 5 | 7 | 16 | 11 | 18 | 1 | 4 | 2 | 17 | 14 | 10 | 6 | 19 | 20 | 8 | 12 | 13 | 15 | 9 | |
| | 1 | 2 | 3 | 4 | 5 | 6 | 7 | 8 | 9 | 10 | 11 | 12 | 13 | 14 | 15 | 16 | 17 | 18 | 19 | 20 | |
| 1 | 0 | 73 | 2 | 39 | 80 | 48 | ∞ | 29 | 26 | 97 | 10 | 77 | 4 | 39 | 69 | 37 | 6 | 96 | 93 | 27 | 1 |
| 2 | 68 | 0 | 49 | 90 | 2 | 81 | 25 | 29 | ∞ | 5 | 87 | 59 | 79 | 22 | 84 | 73 | 8 | 30 | 88 | 0 | 2 |
| 3 | ∞ | 5 | 0 | 28 | 82 | 1 | -1 | 21 | 69 | 49 | 56 | 88 | 59 | 23 | 45 | 96 | 55 | 2 | 64 | 97 | 3 |
| 4 | 22 | 32 | 4 | 0 | 83 | 47 | 29 | ∞ | 30 | 80 | 28 | 47 | 58 | 69 | 16 | 2 | 17 | 56 | 46 | 25 | 4 |
| 5 | 5 | ∞ | 27 | 92 | 0 | 7 | 72 | 31 | 1 | 18 | 45 | 85 | 80 | 78 | 18 | 78 | 95 | 1 | 79 | 3 | 5 |
| 6 | 60 | 81 | 0 | 37 | 38 | 0 | 4 | 58 | 80 | 49 | 45 | 8 | ∞ | 40 | 49 | 5 | 40 | 2 | 61 | 42 | 6 |
| 7 | -1 | 26 | ∞ | 28 | 70 | 56 | 0 | 2 | 30 | 42 | 36 | 46 | -2 | 85 | 39 | 31 | 7 | 53 | 84 | 75 | 7 |
| 8 | 95 | 77 | 31 | -1 | 40 | 75 | 35 | 0 | 74 | 61 | 93 | 52 | 3 | 57 | -2 | ∞ | 88 | 82 | 72 | 11 | 8 |
| 9 | 95 | 1 | 74 | 27 | 24 | 59 | 24 | 22 | 0 | 69 | 94 | 46 | 39 | 83 | -2 | 10 | 77 | 94 | -1 | ∞ | 9 |
| 10 | 87 | 84 | 46 | 22 | 34 | 87 | 75 | 45 | 60 | 0 | 59 | ∞ | 6 | 61 | 33 | 52 | -2 | 92 | 27 | 47 | 10 |
| 11 | 67 | -15 | 56 | -14 | ∞ | 57 | 64 | 67 | -11 | -11 | 14 | 20 | 22 | 5 | 9 | 26 | -14 | 44 | -6 | 11 | 11 |
| 12 | 56 | 96 | 9 | 33 | 1 | 16 | 6 | 17 | 9 | 1 | 62 | 0 | 40 | 70 | 29 | 90 | ∞ | 83 | 55 | 80 | 12 |
| 13 | 1 | 2 | 53 | 27 | 58 | -1 | 94 | 54 | 31 | 46 | 66 | 90 | 0 | 20 | 83 | 82 | 81 | ∞ | 63 | 95 | 13 |
| 14 | 54 | 43 | 34 | 15 | 12 | 45 | 8 | 26 | 87 | 45 | ∞ | 57 | 43 | 0 | 3 | 91 | 60 | 64 | 23 | 93 | 14 |
| 15 | 64 | 70 | 85 | 87 | 9 | 85 | 92 | 45 | 89 | 83 | 25 | 28 | 60 | 6 | 0 | 73 | 54 | 64 | ∞ | 1 | 15 |
| 16 | 28 | 92 | 29 | ∞ | 1 | 16 | 38 | -1 | 91 | 52 | 17 | 23 | 36 | 90 | 88 | 0 | 32 | 28 | 87 | 29 | 16 |
| 17 | 49 | 18 | 43 | 52 | 1 | 64 | 96 | 79 | 6 | ∞ | 47 | 1 | 48 | 46 | 88 | 62 | 0 | 47 | 83 | 71 | 17 |
| 18 | 1 | 7 | 57 | 16 | 72 | ∞ | 47 | 46 | 82 | 64 | 47 | 88 | -1 | 85 | 14 | 76 | 15 | 0 | 85 | 61 | 18 |
| 19 | 17 | 82 | 80 | 84 | 14 | 79 | 32 | 62 | 17 | 40 | 4 | 54 | 33 | ∞ | 14 | 32 | 63 | 15 | 0 | 79 | 19 |
| 20 | 45 | 18 | 40 | 88 | 24 | 14 | 68 | 15 | 85 | 88 | 5 | 34 | 48 | 20 | ∞ | -1 | 28 | 84 | 0 | 0 | 20 |
| | 1 | 2 | 3 | 4 | 5 | 6 | 7 | 8 | 9 | 10 | 11 | 12 | 13 | 14 | 15 | 16 | 17 | 18 | 19 | 20 | |

We now try to obtain a negative permutation, $s$, such that $D_{3\text{-}cycle\ 2}s$ is a derangement containing $n$ edges.
j = 4

(11 4)(4 3) = (11 3): -10

j = 9

(11 9)(9 2) = )11 2): -10

$P_{20}$

| | 1 | 2 | 3 | 4 | 5 | 6 | 7 | 8 | 9 | 10 | 11 | 12 | 13 | 14 | 15 | 16 | 17 | 18 | 19 | 20 | |
|---|---|---|---|---|---|---|---|---|---|---|---|---|---|---|---|---|---|---|---|---|---|
| 1 | | | | | | | | | | | | | | | | | | | | | 1 |
| 2 | | | | | | | | | | | | | | | | | | | | | 2 |
| 3 | | | | | | | | | | | | | | | | | | | | | 3 |
| 4 | | | | | | | | | | | | | | | | | | | | | 4 |
| 5 | | | | | | | | | | | | | | | | | | | | | 5 |
| 6 | | | | | | | | | | | | | | | | | | | | | 6 |
| 7 | | | | | | | | | | | | | | | | | | | | | 7 |
| 8 | | | | | | | | | | | | | | | | | | | | | 8 |
| 9 | | | | | | | | | | | | | | | | | | | | | 9 |
| 10 | | | | | | | | | | | | | | | | | | | | | 10 |
| 11 | | 9 | 4 | 11 | | | | | 11 | | | | | | | | | | | | 11 |
| 12 | | | | | | | | | | | | | | | | | | | | | 12 |
| 13 | | | | | | | | | | | | | | | | | | | | | 13 |
| 14 | | | | | | | | | | | | | | | | | | | | | 14 |
| 15 | | | | | | | | | | | | | | | | | | | | | 15 |
| 16 | | | | | | | | | | | | | | | | | | | | | 16 |
| 17 | | | | | | | | | | | | | | | | | | | | | 17 |
| 18 | | | | | | | | | | | | | | | | | | | | | 18 |
| 19 | | | | | | | | | | | | | | | | | | | | | 19 |
| 20 | | | | | | | | | | | | | | | | | | | | | 20 |
| | 1 | 2 | 3 | 4 | 5 | 6 | 7 | 8 | 9 | 10 | 11 | 12 | 13 | 14 | 15 | 16 | 17 | 18 | 19 | 20 | |

$D^{-1}_{3\text{-cycle }2}M^{-}(20)$

| | 3 | 5 | 7 | 16 | 11 | 18 | 1 | 4 | 2 | 17 | 14 | 10 | 6 | 19 | 20 | 8 | 12 | 13 | 15 | 9 | |
|---|---|---|---|---|---|---|---|---|---|---|---|---|---|---|---|---|---|---|---|---|---|
| | 1 | 2 | 3 | 4 | 5 | 6 | 7 | 8 | 9 | 10 | 11 | 12 | 13 | 14 | 15 | 16 | 17 | 18 | 19 | 20 | |
| 1 | 0 | 73 | 2 | 39 | 80 | 48 | ∞ | 29 | 26 | 97 | 10 | 77 | 4 | 39 | 69 | 37 | 6 | 96 | 93 | 27 | 1 |
| 2 | 68 | 0 | 49 | 90 | 2 | 81 | 25 | 29 | ∞ | 5 | 87 | 59 | 79 | 22 | 84 | 73 | 8 | 30 | 88 | 0 | 2 |
| 3 | ∞ | 5 | 0 | 28 | 82 | 1 | -1 | 21 | 69 | 49 | 56 | 88 | 59 | 23 | 45 | 96 | 55 | 2 | 64 | 97 | 3 |
| 4 | 22 | 32 | 4 | 0 | 83 | 47 | 29 | ∞ | 30 | 80 | 28 | 47 | 58 | 69 | 16 | 2 | 17 | 56 | 46 | 25 | 4 |
| 5 | 5 | ∞ | 27 | 92 | 0 | 7 | 72 | 31 | 1 | 18 | 45 | 85 | 80 | 78 | 18 | 78 | 95 | 1 | 79 | 3 | 5 |
| 6 | 60 | 81 | 0 | 37 | 38 | 0 | 4 | 58 | 80 | 49 | 45 | 8 | ∞ | 40 | 49 | 5 | 40 | 2 | 61 | 42 | 6 |
| 7 | -1 | 26 | ∞ | 28 | 70 | 56 | 0 | 2 | 30 | 42 | 36 | 46 | -2 | 85 | 39 | 31 | 7 | 53 | 84 | 75 | 7 |
| 8 | 95 | 77 | 31 | -1 | 40 | 75 | 35 | 0 | 74 | 61 | 93 | 52 | 3 | 57 | -2 | ∞ | 88 | 82 | 72 | 11 | 8 |
| 9 | 95 | 1 | 74 | 27 | 24 | 59 | 24 | 22 | 0 | 69 | 94 | 46 | 39 | 83 | -2 | 10 | 77 | 94 | -1 | ∞ | 9 |
| 10 | 87 | 84 | 46 | 22 | 34 | 87 | 75 | 45 | 60 | 0 | 59 | ∞ | 6 | 61 | 33 | 52 | -2 | 92 | 27 | 47 | 10 |
| 11 | 67 | -10 | -10 | -14 | ∞ | 57 | 64 | 67 | -11 | -14 | 14 | 20 | 22 | 5 | 9 | 26 | -14 | 44 | -6 | 11 | 11 |
| 12 | 56 | 96 | 9 | 33 | 1 | 16 | 6 | 17 | 9 | 1 | 62 | 0 | 40 | 70 | 29 | 90 | ∞ | 83 | 55 | 80 | 12 |
| 13 | 1 | 2 | 53 | 27 | 58 | **-1** | 94 | 54 | 31 | 46 | 66 | 90 | 0 | 20 | 83 | 82 | 81 | ∞ | 63 | 95 | 13 |
| 14 | 54 | 43 | 34 | 15 | 12 | 45 | 8 | 26 | 87 | 45 | ∞ | 57 | 43 | 0 | 3 | 91 | 60 | 64 | 23 | 93 | 14 |
| 15 | 64 | 70 | 85 | 87 | 9 | 85 | 92 | 45 | 89 | 83 | 25 | 28 | 60 | 6 | 0 | 73 | 54 | 64 | ∞ | 1 | 15 |
| 16 | 28 | 92 | 29 | ∞ | 1 | 16 | 38 | -1 | 91 | 52 | 17 | 23 | 36 | 90 | 88 | 0 | 32 | 28 | 87 | 29 | 16 |
| 17 | 49 | 18 | 43 | 52 | 1 | 64 | 96 | 79 | 6 | ∞ | 47 | 1 | 48 | 46 | 88 | 62 | 0 | 47 | 83 | 71 | 17 |
| 18 | 1 | 7 | 57 | 16 | 72 | ∞ | 47 | 46 | 82 | 64 | 47 | 88 | -1 | 85 | 14 | 76 | 15 | 0 | 85 | 61 | 18 |
| 19 | 17 | 82 | 80 | 84 | 14 | 79 | 32 | 62 | 17 | 40 | 4 | 54 | 33 | ∞ | 14 | 32 | 63 | 15 | 0 | 79 | 19 |
| 20 | 45 | 18 | 40 | 88 | 24 | 14 | 68 | 15 | 85 | 88 | 5 | 34 | 48 | 20 | ∞ | -1 | 28 | 84 | 0 | 0 | 20 |
| | 1 | 2 | 3 | 4 | 5 | 6 | 7 | 8 | 9 | 10 | 11 | 12 | 13 | 14 | 15 | 16 | 17 | 18 | 19 | 20 | |

$P_{40}$

| | 1 | 2 | 3 | 4 | 5 | 6 | 7 | 8 | 9 | 10 | 11 | 12 | 13 | 14 | 15 | 16 | 17 | 18 | 19 | 20 | |
|---|---|---|---|---|---|---|---|---|---|---|---|---|---|---|---|---|---|---|---|---|---|
| 1 | | | | | | | | | | | | | | | | | | | | | 1 |
| 2 | | | | | | | | | | | | | | | | | | | | | 2 |
| 3 | | | | | | | | | | | | | | | | | | | | | 3 |
| 4 | | | | | | | | | | | | | | | | | | | | | 4 |
| 5 | | | | | | | | | | | | | | | | | | | | | 5 |
| 6 | | | | | | | | | | | | | | | | | | | | | 6 |
| 7 | | | | | | | | | | | | | | | | | | | | | 7 |
| 8 | | | | | | | | | | | | | | | | | | | | | 8 |
| 9 | | | | | | | | | | | | | | | | | | | | | 9 |
| 10 | | | | | | | | | | | | | | | | | | | | | 10 |
| 11 | 18 | 9 | 4 | 11 | 16 | 3 | 6 | 7 | 11 | | | | | | | 11 | | 3 | | | 11 |
| 12 | | | | | | | | | | | | | | | | | | | | | 12 |
| 13 | | | | | | | | | | | | | | | | | | | | | 13 |
| 14 | | | | | | | | | | | | | | | | | | | | | 14 |
| 15 | | | | | | | | | | | | | | | | | | | | | 15 |
| 16 | | | | | | | | | | | | | | | | | | | | | 16 |
| 17 | | | | | | | | | | | | | | | | | | | | | 17 |
| 18 | | | | | | | | | | | | | | | | | | | | | 18 |
| 19 | | | | | | | | | | | | | | | | | | | | | 19 |
| 20 | | | | | | | | | | | | | | | | | | | | | 20 |
| | 1 | 2 | 3 | 4 | 5 | 6 | 7 | 8 | 9 | 10 | 11 | 12 | 13 | 14 | 15 | 16 | 17 | 18 | 19 | 20 | |

$D^{-1}_{3\text{-}cycle\,2}M^{-}(40)$

| | 3 | 5 | 7 | 16 | 11 | 18 | 1 | 4 | 2 | 17 | 14 | 10 | 6 | 19 | 20 | 8 | 12 | 13 | 15 | 9 | |
|---|---|---|---|---|---|---|---|---|---|---|---|---|---|---|---|---|---|---|---|---|---|
| | 1 | 2 | 3 | 4 | 5 | 6 | 7 | 8 | 9 | 10 | 11 | 12 | 13 | 14 | 15 | 16 | 17 | 18 | 19 | 20 | |
| 1 | 0 | 73 | 2 | 39 | 80 | 48 | ∞ | 29 | 26 | 97 | 10 | 77 | 4 | 39 | 69 | 37 | 6 | 96 | 93 | 27 | 1 |
| 2 | 68 | 0 | 49 | 90 | 2 | 81 | 25 | 29 | ∞ | 5 | 87 | 59 | 79 | 22 | 84 | 73 | 8 | 30 | 88 | 0 | 2 |
| 3 | ∞ | 5 | 0 | 28 | 82 | 1 | -1 | 21 | 69 | 49 | 56 | 88 | 59 | 23 | 45 | 96 | 55 | 2 | 64 | 97 | 3 |
| 4 | 22 | 32 | 4 | 0 | 83 | 47 | 29 | ∞ | 30 | 80 | 28 | 47 | 58 | 69 | 16 | 2 | 17 | 56 | 46 | 25 | 4 |
| 5 | 5 | ∞ | 27 | 92 | 0 | 7 | 72 | 31 | 1 | 18 | 45 | 85 | 80 | 78 | 18 | 78 | 95 | 1 | 79 | 3 | 5 |
| 6 | 60 | 81 | 0 | 37 | 38 | 0 | 4 | 58 | 80 | 49 | 45 | 8 | ∞ | 40 | 49 | 5 | 40 | 2 | 61 | 42 | 6 |
| 7 | -1 | 26 | ∞ | 28 | 70 | 56 | 0 | 2 | 30 | 42 | 36 | 46 | -2 | 85 | 39 | 31 | 7 | 53 | 84 | 75 | 7 |
| 8 | 95 | 77 | 31 | -1 | 40 | 75 | 35 | 0 | 74 | 61 | 93 | 52 | 3 | 57 | -2 | ∞ | 88 | 82 | 72 | 11 | 8 |
| 9 | 95 | 1 | 74 | 27 | 24 | 59 | 24 | 22 | 0 | 69 | 94 | 46 | 39 | 83 | -2 | 10 | 77 | 94 | -1 | ∞ | 9 |
| 10 | 87 | 84 | 46 | 22 | 34 | 87 | 75 | 45 | 60 | 0 | 59 | ∞ | 6 | 61 | 33 | 52 | -2 | 92 | 27 | 47 | 10 |
| 11 | -7 | -10 | *-10* | -14 | -1 | *-9* | *-5* | *-3* | -11 | -14 | 14 | 20 | 22 | 5 | 9 | 26 | -14 | *-8* | -6 | 11 | 11 |
| 12 | 56 | 96 | 9 | 33 | 1 | 16 | 6 | 17 | 9 | 1 | 62 | 0 | 40 | 70 | 29 | 90 | ∞ | 83 | 55 | 80 | 12 |
| 13 | 1 | 2 | 53 | 27 | 58 | -1 | 94 | 54 | 31 | 46 | 66 | 90 | 0 | 20 | 83 | 82 | 81 | ∞ | 63 | 95 | 13 |
| 14 | 54 | 43 | 34 | 15 | 12 | 45 | 8 | 26 | 87 | 45 | ∞ | 57 | 43 | 0 | 3 | 91 | 60 | 64 | 23 | 93 | 14 |
| 15 | 64 | 70 | 85 | 87 | 9 | 85 | 92 | 45 | 89 | 83 | 25 | 28 | 60 | 6 | 0 | 73 | 54 | 64 | ∞ | 1 | 15 |
| 16 | 28 | 92 | 29 | ∞ | 1 | 16 | 38 | -1 | 91 | 52 | 17 | 23 | 36 | 90 | 88 | 0 | 32 | 28 | 87 | 29 | 16 |
| 17 | 49 | 18 | 43 | 52 | 1 | 64 | 96 | 79 | 6 | ∞ | 47 | 1 | 48 | 46 | 88 | 62 | 0 | 47 | 83 | 71 | 17 |
| 18 | 1 | 7 | 57 | 16 | 72 | ∞ | 47 | 46 | 82 | 64 | 47 | 88 | -1 | 85 | 14 | 76 | 15 | 0 | 85 | 61 | 18 |
| 19 | 17 | 82 | 80 | 84 | 14 | 79 | 32 | 62 | 17 | 40 | 4 | 54 | 33 | ∞ | 14 | 32 | 63 | 15 | 0 | 79 | 19 |
| 20 | 45 | 18 | 40 | 88 | 24 | 14 | 68 | 15 | 85 | 88 | 5 | 34 | 48 | 20 | ∞ | -1 | 28 | 84 | 0 | 0 | 20 |
| | 1 | 2 | 3 | 4 | 5 | 6 | 7 | 8 | 9 | 10 | 11 | 12 | 13 | 14 | 15 | 16 | 17 | 18 | 19 | 20 | |

We cannot extend the negative paths (11 1) or (11 5) to negative paths smaller than they are. Before constructing a tour by patching together cycles of $D_{3\text{-}cycle\ 2}$, we prove the following theorem that may reduce the running time of the patching algorithm.

**Theorem 15** *Let two branches obtained from our algorithm go through the same cycles and have the following properties: (1) the terminal vertex of both branches is* $a$. *(2) If the arcs deleted are of the respective forms*

$$\begin{array}{ccc} \Gamma: & \qquad & \Delta: \\ \mathbf{A} \downarrow \mathbf{B} & & \mathbf{C} \downarrow \mathbf{D} \\ b \rightarrow a \leftarrow c & & b \rightarrow a \leftarrow c \end{array}$$

*with* $\max\{\mathbf{A}, \mathbf{B}\} \le \min\{\mathbf{C}, \mathbf{D}\}$, *then we can continue our algorithm using only* $\Gamma$.

**Proof**. *Our paths from* $b, a$ *of* $\Gamma$ *will always be no greater than if we had used* $b, a$ *of* $\Delta$.

A tour of smallest value occurs when our initial vertex is 11.

$C_3$

11

0 ↓

10→→17←←12

**6** ↓ $C_{3,5}$ ↓

13→ 6←18 7 →1←3

↙ **8** $C_{3,4,5}$ ↙ ↘ $C_{1,3,5}$ ↘

1→3→ 7 1→3←7 13→ 6 ←18 6→18→13

↙ $C_{1,3,4,5}$ ↓ **12** $C_{1,3,4,5}$ ↓ $C_{1,3,4,5}$ ↓ ↘ $C_{1,3,4,5}$ ↘

8→4←16 8→4←16 4→16←8 4→16→8 16→8←4 4→16←8

$C_{1,2,3,4,5}$ ↓ $C_{1,2,3,4,5}$ ↓ **12** $C_{1,2,3,4,5}$ $C_{1,2,3,4,5}$ $C_{1,2,3,4,5}$ $C_{1,2,3,4,5}$

5 11→14

$T_1 = (11\ \ 17\ \ 12\ \ 10\ \ 6\ \ 18\ \ 13\ \ 3\ \ 1\ \ 7\ \ 4\ \ 8\ \ 16\ \ 14\ \ 19\ \ 15\ \ 20\ \ 9\ \ 2\ \ 5\ ) = |D_{n\ edges}| + 12 = 68$.

At the start of Phase 3, when constructing acceptable paths, $P_{ni}$ and $\sigma_1^{-1}M^-(20)$ are blank tables. We choose an initial arc of a path from $\sigma_1^{-1}M^-$. We then add one arc at a time from $\sigma_1^{-1}M^-$ and place it in $P_{ni}$ (at the same time placing the value of the new path in $\sigma_1^{-1}M^-(20)$) to extend a path. When constructing acceptable paths, we can avoid backtracking in the following case: Let $c$ belong to the 2-cycle $(c\ c')$. If $(b\ c)$ is a new proposed arc to extend an acceptable path $(a\ b)$, we go to $P_{ni}$ and check whether $(a\ c')$ is blank or contains $b$. If so, then $(a\ b)(b\ c) = (a\ c)$ must be an acceptable path. Otherwise, we must backtrack in $P_{ni}$ from $c$ to $a$ to guarantee that no two points lie in the same 2-cycle. In the process of constructing acceptable paths and cycles, we can construct unlinked 2-circuit paths using the following strategy: If $(a\ c)$ is an acceptable path and $(c\ d)$ has the property that $a \neq d$ while $a$ and $d$ lie in the same 2-cycle of $\sigma_1$, then $(a\ d)$ is an unlinked 2-circuit path.

*Note.* When a path, $P$, gets reasonably long, we can place after it the set of 2-cycles *none of which* contains a point of $P$ in it. We then choose only points from 2-cycles in the set, deleting one 2-cycle at a time.

$T_1^{-1}M^-$

| | 7 | 5 | 1 | 8 | 11 | 18 | 4 | 16 | 2 | 6 | 17 | 10 | 3 | 19 | 20 | 14 | 12 | 13 | 15 | 9 | |
|---|---|---|---|---|---|---|---|---|---|---|---|---|---|---|---|---|---|---|---|---|---|
| | 1 | 2 | 3 | 4 | 5 | 6 | 7 | 8 | 9 | 10 | 11 | 12 | 13 | 14 | 15 | 16 | 17 | 18 | 19 | 20 | |
| 1 | **0** | 71 | ∞ | 35 | 78 | 46 | 27 | 40 | 23 | 2 | 95 | 75 | -2 | 37 | 67 | 7 | 4 | 94 | 91 | 25 | 1 |
| 2 | 49 | **0** | 25 | 73 | 2 | 81 | 29 | 91 | ∞ | 79 | 5 | 59 | 68 | 22 | 84 | 86 | 8 | 30 | 88 | 0 | 2 |
| 3 | 1 | 6 | **0** | 97 | 83 | 2 | 22 | 30 | 68 | 60 | 50 | 89 | ∞ | 24 | 46 | 55 | 56 | 3 | 65 | 98 | 3 |
| 4 | 2 | 30 | 27 | **0** | 81 | 45 | ∞ | 1 | 27 | 56 | 78 | 45 | 20 | 67 | 14 | 25 | 15 | 54 | 44 | 23 | 4 |
| 5 | 27 | ∞ | 72 | 78 | **0** | 7 | 31 | 94 | -1 | 80 | 18 | 85 | 5 | 88 | 18 | 43 | 95 | 1 | 70 | 3 | 5 |
| 6 | 0 | 81 | 4 | 5 | 38 | **0** | 58 | 38 | 79 | ∞ | 49 | 8 | 60 | 40 | 49 | 44 | 40 | 2 | 61 | 42 | 6 |
| 7 | ∞ | 24 | -2 | 29 | 68 | 54 | **0** | 26 | 25 | -4 | 40 | 44 | -3 | 83 | 37 | 31 | 5 | 51 | 82 | 73 | 7 |
| 8 | 32 | 78 | 36 | ∞ | 41 | 76 | 1 | **0** | 72 | 4 | 62 | 53 | 96 | 58 | -1 | 91 | 89 | 83 | 73 | 12 | 8 |
| 9 | 77 | 4 | 27 | 13 | 27 | 62 | 25 | 30 | **0** | 42 | 72 | 49 | 98 | 86 | 1 | 94 | 80 | 97 | 2 | ∞ | 9 |
| 10 | 40 | 78 | 69 | 46 | 28 | 81 | 39 | 16 | 51 | **0** | -6 | ∞ | 81 | 55 | 27 | 50 | -8 | 86 | 21 | 41 | 10 |
| 11 | 70 | -1 | 78 | 40 | ∞ | 71 | 81 | 0 | 0 | 36 | **0** | 34 | 81 | 19 | 23 | 11 | 0 | 58 | 8 | 25 | 11 |
| 12 | 9 | 96 | 6 | 90 | 1 | 16 | 17 | 33 | 8 | 40 | 1 | **0** | 56 | 70 | 29 | 61 | ∞ | 83 | 55 | 80 | 12 |
| 13 | -2 | 1 | 14 | 81 | 57 | -2 | 53 | 26 | 27 | -1 | 45 | 89 | **0** | 19 | 82 | 62 | 80 | ∞ | 62 | 94 | 13 |
| 14 | 2 | 43 | 17 | 91 | 12 | 45 | 26 | 12 | 85 | 43 | 45 | 57 | 54 | **0** | 3 | ∞ | 59 | 64 | 23 | 93 | 14 |
| 15 | 2 | 70 | 18 | 73 | 9 | 85 | 45 | 87 | 87 | 60 | 83 | 28 | 64 | 6 | **0** | 23 | 54 | 64 | ∞ | 1 | 15 |
| 16 | -1 | 80 | -6 | -12 | -11 | 4 | -13 | ∞ | 77 | 24 | 40 | 11 | 16 | 78 | 76 | **0** | 20 | 16 | 75 | 17 | 16 |
| 17 | 12 | 18 | 10 | 62 | 1 | 64 | 79 | 52 | 4 | 48 | ∞ | 1 | 49 | 46 | 88 | 45 | **0** | 47 | 83 | 71 | 17 |
| 18 | 14 | 7 | 4 | 76 | 72 | ∞ | 46 | 16 | 80 | -1 | 64 | 88 | 1 | 85 | 14 | 45 | 15 | **0** | 85 | 61 | 18 |
| 19 | -5 | 82 | 6 | 54 | 14 | 79 | 62 | 84 | 15 | 33 | 40 | 56 | 17 | ∞ | 14 | -6 | 63 | 15 | **0** | 79 | 19 |
| 20 | 9 | 18 | 13 | -1 | 24 | 14 | 15 | 88 | 83 | 48 | 88 | 34 | 45 | 20 | ∞ | 3 | 28 | 84 | 0 | **2** | 20 |
| | 1 | 20 | 3 | 4 | 5 | 6 | 7 | 8 | 9 | 10 | 11 | 12 | 13 | 14 | 15 | 16 | 17 | 18 | 19 | 50 | |

In general, by using MIN(M) and giving the inverse of a derangement, we can obtain the values of any row in ascending order of magnitude**.**

| | | | | | | | | | | | | | | | | | | | | |
|---|---|---|---|---|---|---|---|---|---|---|---|---|---|---|---|---|---|---|---|---|
| | 1 | 2 | 3 | 4 | 5 | 6 | 7 | 8 | 9 | 10 | 11 | 12 | 13 | 14 | 15 | 16 | 17 | 18 | 19 | 20 |
| $T_1^{-1}$ = | | | | | | | | | | | | | | | | | | | | |
| | 3 | 9 | 13 | 7 | 2 | 10 | 1 | 4 | 20 | 12 | 5 | 17 | 18 | 16 | 19 | 8 | 11 | 6 | 14 | 15 |

In this case, suppose I want to find the third smallest value in row 5. MIN(M)(5, 3) = 13. . $T_1^{-1}(13) = 18$

This, (5 18) has the third smallest value in row 5 of $T_1^{-1}M^-$. Since $n$ is even, by choosing those arcs in $T_1^{-1}M^-$ that represent the inverse, $T_1^{-1}$, of $T_1$, we obtain two disjoint cycles - $s_-$ and $s_+$ - each of which moves $\frac{n}{2}$ points.

Furthermore, the sum of their respective values have the same absolute value but different signs. This makes sense since $T_1 s_- s_2 = T_1^{-1}$. On the other hand, $T_1 s_-$ is a product of $\frac{n}{2}$ pair-wise disjoint cycles, i. e., a PM. We now construct $T_1 s_-$.= $\sigma_1$. We give the inverse of $\sigma_1$ in order to use it and MIN(M) to obtain the smallest entry values in each row.

| | | | | | | | | | | | | | | | | | | | | |
|---|---|---|---|---|---|---|---|---|---|---|---|---|---|---|---|---|---|---|---|---|
| | 1 | 2 | 3 | 4 | 5 | 6 | 7 | 8 | 9 | 10 | 11 | 12 | 13 | 14 | 15 | 16 | 17 | 18 | 19 | 20 |
| $\sigma_1^{-1}$ = | | | | | | | | | | | | | | | | | | | | |
| | 3 | 9 | 1 | 7 | 11 | 10 | 4 | 16 | 2 | 6 | 5 | 17 | 18 | 19 | 20 | 8 | 12 | 13 | 14 | 15 |

*Note*. We can also obtain $\sigma_1$ as the smaller-valued set of alternating edges of $T_1$. When represented as a set of edges, $\sigma_1$ is a perfect matching .

Before proceeding further, we denote the 2-cycles of $\sigma$ in the following way:

**1** (1 3); **2** (2 9); **3** (4 7); **4** (5 11); **5** (6 10); **6** (8 16); **7** (12 17); **8** (13 18); **9** (14 19); **10** (15 20) .

A cycle may be obtained more than once using different determining vertices. In that case, we include only its first occurrence. The remaining times that it occurs are included in our list of determining vertices. Finally, every acceptable cycle has a companion cycle that consists of arcs that are symmetric to those of the acceptable cycle. Together, when they multiply a derangement consisting only of 2-cyles, they form another derangement of the same form. In other words, they transform a perfect matching into another perfect matching. We don't include companion cycles in our list of acceptable cycles.

$\sigma_1^{-1} M^-$

| | 3 | 9 | 1 | 7 | 11 | 10 | 4 | 16 | 2 | 6 | 5 | 17 | 18 | 19 | 20 | 8 | 12 | 13 | 14 | 15 | |
|---|---|---|---|---|---|---|---|---|---|---|---|---|---|---|---|---|---|---|---|---|---|
| | 1 | 2 | 3 | 4 | 5 | 6 | 7 | 8 | 9 | 10 | 11 | 12 | 13 | 14 | 15 | 16 | 17 | 18 | 19 | 20 | |
| 1 | **0** | 27 | ∞ | 2 | 80 | 77 | 29 | 39 | 25 | 4 | 73 | 97 | 48 | 39 | 69 | 37 | 6 | 96 | 9 | 93 | 1 |
| 2 | 68 | **0** | 25 | 49 | 2 | 59 | 29 | 90 | ∞ | 79 | 0 | 5 | 81 | 22 | 84 | 73 | 8 | 30 | 86 | 88 | 2 |
| 3 | ∞ | 98 | **0** | 1 | 83 | 89 | 22 | 29 | 68 | 60 | 6 | 50 | 2 | 24 | 46 | 97 | 56 | 3 | 55 | 65 | 3 |
| 4 | 18 | 21 | 25 | **0** | 79 | 43 | ∞ | -4 | 25 | 53 | 28 | 76 | 43 | 65 | 12 | -2 | 13 | 52 | 23 | 42 | 4 |
| 5 | 5 | 3 | 72 | 27 | **0** | 85 | 31 | 92 | -1 | 80 | ∞ | 18 | 7 | 88 | 18 | 78 | 95 | 1 | 43 | 70 | 5 |
| 6 | 52 | 34 | -4 | -8 | 30 | **0** | 50 | 29 | 71 | ∞ | 73 | 41 | -8 | 32 | 41 | -3 | 32 | -6 | 36 | 53 | 6 |
| 7 | **-3** | 73 | -2 | ∞ | 68 | 44 | **0** | 26 | 25 | -4 | 24 | 40 | 54 | 83 | 37 | 29 | 5 | 51 | 31 | 82 | 7 |
| 8 | 96 | 12 | 36 | 32 | 41 | 53 | 1 | **0** | 72 | 4 | 78 | 62 | 76 | 58 | -1 | ∞ | 89 | 83 | 91 | 73 | 8 |
| 9 | 98 | ∞ | 27 | 77 | 27 | 49 | 25 | 30 | **0** | 42 | 4 | 72 | 62 | 86 | 1 | 13 | 80 | 97 | 94 | 2 | 9 |
| 10 | 81 | 41 | 69 | 40 | 29 | ∞ | 39 | 16 | 51 | **0** | 78 | -6 | 81 | 55 | 45 | 46 | -8 | 86 | 50 | 21 | 10 |
| 11 | 82 | 26 | 79 | 71 | ∞ | 35 | 82 | 1 | 1 | 37 | **0** | 1 | 72 | 20 | 24 | 41 | 1 | 59 | 12 | 9 | 11 |
| 12 | 55 | 79 | 5 | 8 | 0 | -1 | 16 | 32 | 7 | 39 | 95 | **0** | 15 | 69 | 28 | 89 | ∞ | 82 | 60 | 54 | 12 |
| 13 | 2 | 96 | 95 | 54 | 59 | 91 | 55 | 28 | 29 | 1 | 3 | 47 | **0** | 21 | 84 | 83 | 82 | ∞ | 64 | 64 | 13 |
| 14 | 54 | 93 | 8 | 34 | 12 | 57 | 26 | 15 | 85 | 43 | 43 | 45 | 45 | **0** | 3 | 91 | 60 | 64 | ∞ | 23 | 14 |
| 15 | 64 | 1 | 92 | 85 | 9 | 28 | 45 | 87 | 87 | 60 | 70 | 83 | 85 | 6 | **0** | 73 | 54 | 64 | 23 | ∞ | 15 |
| 16 | 28 | 29 | 38 | 29 | 1 | 23 | -1 | ∞ | 89 | 36 | 92 | 64 | 64 | 46 | 88 | **0** | 52 | 47 | 45 | 83 | 16 |
| 17 | 49 | 71 | 96 | 43 | 1 | 1 | 79 | 52 | 4 | 48 | 18 | ∞ | 64 | 46 | 88 | 62 | **0** | 47 | 45 | 83 | 17 |
| 18 | 1 | 61 | 47 | 57 | 72 | 88 | 46 | 16 | 80 | -1 | 7 | 64 | ∞ | 85 | 14 | 76 | 15 | **0** | 45 | 85 | 18 |
| 19 | 23 | 85 | 12 | 86 | 20 | 62 | 68 | 90 | 21 | 39 | 88 | 46 | 85 | ∞ | 20 | 58 | 69 | 21 | **0** | 6 | 19 |
| 20 | 45 | 0 | 68 | 40 | 24 | 34 | 15 | 88 | 83 | 48 | 18 | 88 | 14 | 20 | ∞ | -1 | 28 | 84 | 3 | **0** | 20 |
| | 1 | 2 | 3 | 4 | 5 | 6 | 7 | 8 | 9 | 10 | 11 | 12 | 13 | 14 | 15 | 16 | 17 | 18 | 19 | 20 | |

We now construct acceptable paths to obtain acceptable and 2-circuit cycles. We are assuming that the value of a positive cycle is no greater than $|T_1| - |\sigma_1| - 1 = 11$

j = 1

(5 1)(1 4) = (5 4): 7; (5 1)(1 10) = (5 10): 9; (5 1)(1 17) = (5 17): 11; (7 1)(1 17) = (7 17): 3; (7 1)(1 19) = (7 19): 6; (13 1)(1 4) = (13 4): 4; (13 1)(1 17) = (13 17): 8;

(13 1)(1 19) = (13 19): 11; (18 1)(1 4) = (18 4): 3; (18 1)(1 17) = (18 17): 7;

(18 1)(1 19) = (18 19): 10

j = 2

**(5 2)(2 5) = (5 5): 5**

**CYCLE P = [5 9 2 11 5]: (5 2): 5**

(5 2)(2 12) = (5 12): 8; (15 2)(2 5) = (15 5): 3; (15 2)(2 11) = (15 11): 1;

(15 2)(2 12) = (15 12): 6; (15 2)(2 17) = (15 17): 9; (20 2)(2 5) = (20 5): 3;

(20 2)(2 11) = (20 11): 0; (20 2)(2 12) = (20 12): 5; (20 2)(2 17) = (20 17): 8

j = 3

(6 3)(3 11) = (6 11): 2; (7 3)(3 11) = (7 11): 4; (7 3)(3 13) = (7 13): 0; (7 3)(3 18) = (7 18): 1;

(12 3)(3 4) = (12 4): 6; (12 3)(3 13) = (12 13): 7; (12 3)(3 18) = (13 18): 8;

(14 3)(3 4) = (14 4): 9; (14 3)(3 13) = (14 13): 10; (14 3)(3 18) = (14 18): 11

j = 4

(1 4)(4 8) = (1 8): -2; (1 4)(4 16) = (1 16): 0; (3 4)(4 8) = (3 8): -3; (3 4)(4 16) = (3 16): 0;

(5 4)(4 8) = (5 8): 3; (5 4)(4 16) = (5 16): 5; (6 4)(4 8) = (6 8): -12; (6 4)(4 16) = (6 16): -10;

(12 4)(4 8) = (12 8): 2; (12 4)(4 16) = (12 16): 4; (13 4)(4 8) = (13 8): 0;

(13 4)(4 16) = (13 16): 2; (14 4)(4 8) = (14 8): 5; (14 4)(4 16) = (14 16): 7;

(18 4)(4 8) = (18 8): -1; (18 4)(4 16) = (18 16): 1

j = 5

(2 5)(5 1) = (2 1): 7;

**(2 5)(5 2) = (2 2): 5**

**CYCLE P = [5 9 2 11 5]: (2 5): 5**

(2 5)(5 13) = (2 13): 3; (2 5)(5 18) = (2 18): 3; (12 5)(5 1) = (12 1): 5;

(12 5)(5 2) = (12 2): 3; (12 5)(5 9) = (12 9): -1; (12 5)(5 13) = (12 13): 1; (12 5)(5 18) = (12 18): 1;

(15 5)(5 1) = (15 1): 8; (15 5)(5 13) = (15 13): 4; (15 5)(5 18) = (15 18): 4;

(16 5)(5 1) = (16 1): 6; (16 5)(5 2) = (16 2): 4; (16 5)(5 9) = (16 9): 0; (16 5)(5 13) = (16 13): 2;

(16 5)(5 18) = (16 18): 2; (17 5)(5 1) = (17 1): 6; (17 5)(5 2) = (17 2): 4;

(17 5)(5 9) = (17 9): 0; (17 5)(5 13) = (17 13): 2; (17 5)(5 18) = (17 18): 2;

(20 5)(5 1) = (20 1): 7; (20 5)(5 13) = (20 13): 3; (20 5)(5 18) = (20 18): 3

j = 6

(12 6)(6 3) = (12 3): -5; (12 6)(6 4) = (12 4): -9; (12 6)(6 13) = (12 13): -9;

(12 6)(6 16) = (12 16): -4; (12 6)(6 18) = (12 18): -7; (17 6)(6 3) = (17 3): -3;

(17 6)(6 4) = (17 4): -7; (17 6)(6 13) = (17 13): -8; (17 6)(6 16) = (17 16): -2;

(17 6)(6 18) = (17 18): -5

j = 7

(8 7)(7 1) = (8 1): -2; (8 7)(7 3) = (8 3): -1; (8 7)(7 10) = (8 10): -3; (8 7)(7 17) = (8 17): 6;

(16 7)(7 1) = (16 1) -4; (16 7)(7 3) = (16 3): -3; (16 7)(7 10) = (16 10): -5; (16 7)(7 17): 4

j = 8

(1 8)(8 10) = (1 10): 2; (1 8)(8 15) = (1 15): -3; (3 8)(8 10) = (3 10): 1; (3 8)(8 15) = (3 15): -4;

(4 8)(8 10) = (4 10): 0; (5 8)(8 7) = (5 7): 4; (5 8)(8 10) = (5 10): 7; (5 8)(8 15) = (5 15): 2;

(6 8)(8 15) = (6 15): -13; (11 8)(8 7) = (11 7): 2; (11 8)(8 10) = (11 10): 5;

(11 8)(8 15) = (11 15): 0; (12 8)(8 10) = (12 10): 6; (12 8)(8 15) = (12 15): 1;

(13 8)(8 15) = (13 15): -1; (14 8)(8 10) = (14 10): 9

j = 9

(5 9)(9 15) = (5 15): 0; (5 9)(9 20) = (5 20): 1;

**(11 9)(9 11) = (11 11): 5**

**CYCLE P = [11 2 9 5 11]: (11 9): 5**

(11 9)(9 20) = (11 20): 3; (12 9)(9 20) = (12 20): 1; (17 9)(9 15) = (17 15): 1;

(17 9)(9 20) = (17 20): 2

j = 10

(1 10)(10 12) = (1 12): -4; (1 10)(10 17) = (1 17): -6; (3 10)(10 12) = (3 12): -5;

(3 10)(10 17) = (3 17): -7; (4 10)(10 12) = (4 12): -6; (4 10)(10 17) = (4 17): -8;

(5 10)(10 12) = (5 12): 1; (5 10)(10 17) = (5 17): -1; (7 10)(10 12) = (7 12): -10;

(7 10)(10 17) = (7 17): -12; (8 10)(10 12) = (8 12): -9; (8 10)(10 17) = (8 17): -11;

(11 10)(10 12) = (11 12): -1; (11 10)(10 17) = (11 17): -3; (13 10)(10 12) = (13 12): -5;

(13 10)(10 17) = (13 17): -7; (14 10)(10 12) = (14 12): 3; (14 10)(10 17) = (14 17): 1;

(16 10)(10 12) = (16 12): -11; (16 10)(10 17) = (16 17): -13; (18 10)(10 12) = (18 12): -7;

(18 10)(10 17) = (18 17): -9

j = 11

(3 11)(11 9) = (3 9): 7; (6 11)(11 9) = (6 9): 3; (6 11)(11 12) = (6 12): 3;

(6 11)(11 17) = (6 17): 3; (6 11)(11 20) = (6 20): 11; (7 11)(11 8) = (7 8): 5;

(7 11)(11 9) = (7 9): 5; (9 11)(11 8) = (9 8): 5;

**(9 11)(11 9) = (9 9): 5**

**CYCLE P = [9 5 11 2 9]: (9 11): 5**

(9 11)(11 12) = (9 12): 5; (9 11)(11 17): 5; (13 11)(11 9) = (13 9): 4; (15 11)(11 8) = (15 8): 2;

(15 11)(11 12) = (15 12): 2; (15 11)(11 17) = (15 17): 2; ((18 11)(11 8) = (18 8): 8;

(18 11)(11 9) = (18 9): 8; (20 11)(11 8) = (20 8): 1; (20 11)(11 12) = (20 12): 1;

(20 11)(11 17) = (20 17): 1;

**(20 11)(11 20) = (20 20): 9**

**CYCLE P = [20 9 2 5 11 15 20]: (20 2 11): 9**

j = 12

(1 12)(12 5) = (1 5): -4; (1 12)(12 9) = (1 9): 3;

**(3 12)(12 3) = (3 3): 0**

**CYCLE P = [3 7 4 16 8 6 10 17 12 1 3]: (4 8 10 12): 2**

(3 12)(12 5) = (3 5): -5; (3 12)(12 9) = (3 9): 5; (4 12)(12 3) = (4 3): -1;

**(4 12)(12 4) = (4 4): 2**

**CYCLE P = [4 16 8 6 10 17 12 7 4]: (4 8 10 12): 2**

(4 12)(12 5) = (4 5): -6; (4 12)(12 9) = (4 9): 1;

**(5 12)(12 5) = (5 5): 1**

**CYCLE P = [5 3 1 7 4 16 8 6 10 17 12 11 5]: (5 1 4 8 10 12): 1**

(7 12)(12 3) = (7 3): -5; (7 12)(12 5) = (7 5): -10; (7 12)(12 9) = (7 9): -3;

(8 12)(12 3) = (8 3): -4; (8 12)(12 5) = (8 5): -9; (8 12)(12 9) = (8 9): -2;

(9 12)(12 5) = (9 5): 5; (10 12)(12 3) = (10 3): -1; (10 12)(12 4) = (10 4): 2;

(10 12)(12 5) = (10 5): -6; (10 12)(12 9) = (10 9): 1; (11 12)(12 4) = (11 4): 7;

(11 12)(12 9) = (11 9): 6; (13 12)(12 5) = (13 5): -2; (13 12)(12 9) = (13 9): 2;

(14 12)(12 5) = (14 5): 6; (14 12)(12 9) = (14 9): 10; (15 12)(12 3) = (15 3): 7;

(16 12)(12 3) = (16 3): -6; (16 12)(12 5) = (16 5): -11; (16 12)(12 9) = (16 9): -4;

(18 12)(12 3) = (18 3): -2; (18 12)(12 4) = (18 4): 1; (18 12)(12 5) = (18 5): -7;

(18 12)(12 9) = (18 9): 0; (20 12)(12 3) = (20 3): 10

j = 13

(2 13)(13 1) = (2 1): 5; (2 13)(13 10) = (2 10): 4; (3 13)(13 10) = (3 10): 3;

(3 13)(13 11) = (3 11): 5; (6 13)(13 1) = (6 1): -6; (6 13)(13 11) = (6 11): -5;

(7 13)(13 11) = (7 11): 3; (12 13)(13 1) = (12 1): -7; (12 13)(13 10) = (12 10): -8;

(15 13)(13 1) = (15 1): 6; (15 13)(13 10) = (15 10): 5; (17 13)(13 1) = (17 1): -6;

(17 13)(13 11) = (17 11): -5; (20 13)(13 1) = (20 1): 5; (20 13)(13 10) = (20 10): 4

j = 14

**(15 14)(14 15) = (15 15): 9**

**CYCLE P = [15 19 14 20 15]: (15 14): 9**

j = 15

(1 15)(15 2) = (1 2): -2; (3 15)(15 2) = (3 2): -3;

**(5 15)(15 5) = (5 5): 9**

**CYCLE P = [5 2 9 20 15 11 5]: (5 9 15): 9**

(5 15)(15 14) = (5 14): 6; (6 15)(15 14) = (6 14): -7; (8 15)(15 2) = (8 2): 0;

(8 15)(15 14) = (8 14): 5; (9 15)(15 14) = (9 14): 7; (11 15)(15 2) = (11 2): 1;

(12 15)(15 2) = (12 2): 2; (12 15)(15 14) = (12 14): 7; (13 15)(15 2) = (13 2): 0;

(13 15)(15 14) = (13 14): 5; (14 15)(15 2) = (14 2): 4;

**(14 15)(15 14) = (14 14): 9**

**CYCLE P = [14 20 15 19 14]: (14 15): 9**

j = 16

**(5 16)(16 5) = (5 5): 6**

**CYCLE P = [5 3 1 7 4 8 16 11 5]: (5 1 4 16): 6**

(6 16)(16 5) = (6 5): -11; (17 16)(16 5) = (17 5): -1; (17 16)(16 7) = (17 7): -3;

(20 16)(16 5) = (20 5): 0; (20 16)(16 7) = (20 7): -2

j = 17

(1 17)(17 5) = (1 5): -5; (1 17)(17 9) = (1 9): -2; (2 17)(17 6) = (2 6): 9;

(3 17)(17 5) = (3 5): -6; (3 17)(17 5) = (3 5): -6; (3 17)(17 9) = (3 9): -3;

(4 17)(17 5) = (4 5): -7; (4 17)(17 9) = (4 9): -4;

**(5 17)(17 5) = (5 5): 0**

**CYCLE P = [5 3 1 7 4 16 8 6 10 12 17 11 5]: (5 1 4 8 10 17): 0**

(7 17)(17 5) = (7 5): -11; (7 17)(17 9) = (7 9): -8; (8 17)(17 5) = (8 5): -10;

(8 17)(17 9) = (8 9): -7; (9 17)(17 6) = (9 6): 6;

**(9 17)(17 9) = (9 9): 9**

**CYCLE P = [9 5 11 12 17 2 9]: (9 11 17): 9**

(10 17)(17 5) = (10 5): -7; (13 17)(17 5) = (13 5): -6; (13 17)(17 9) = (13 9): -3;

(14 17)(17 5) = (14 5): 2; (14 17)(17 9) = (14 9): 5; (15 17)(17 6) = (15 6): 10;

(16 17)(17 5) = (16 5): -12; (16 17)(17 9) = (16 9): -9; (18 17)(17 5) = (18 5): -8;

(18 17)(17 9) = (18 9): -5; (20 17)(17 6) = (20 6): 9

j = 18

(2 18)(18 1) = (2 1): 4; (2 18)(18 10) = (2 10): 2; (5 18)(18 1) = (5 1): 2;

(5 18)(18 10) = (5 10): 0; (15 18)(18 1) = (15 1): 5; (20 18)(18 1) = (20 1): 4;

(20 18)(18 10) = (20 10): 2

j = 19

**(20 19)(19 20) = (20 20): 9**

**CYCLE P = [20 14 19 15 20]: (20 19): 9**

j = 20

(5 20)(20 16) = (5 16): 0; (5 20)(20 14) = (5 14): 4; (6 20)(20 2) = (6 2): 1;

(6 20)(20 19) = (6 19): 4; (9 20)(20 16) = (9 16): 1; (9 20)(20 19) = (9 19): 5;

(11 20)(20 19) = (11 19): 6; (12 20)(20 19) = (12 19): 4;

**(19 20)(20 19) = (19 19): 9**

**CYCLE P = [19 15 20 14 19]: (19 20): 9**

$P_{20}$

| | 1 | 2 | 3 | 4 | 5 | 6 | 7 | 8 | 9 | 10 | 11 | 12 | 13 | 14 | 15 | 16 | 17 | 18 | 19 | 20 | |
|---|---|---|---|---|---|---|---|---|---|---|---|---|---|---|---|---|---|---|---|---|---|
| 1 | | 15 | 12 | 1 | 17 | | | 4 | 17 | 8 | | 10 | | | 8 | 4 | 10 | | | | 1 |
| 2 | 18 | | | | 2 | 17 | | | | 18 | | | 5 | | | | 1 | 5 | 1 | | 2 |
| 3 | | 15 | | 3 | 17 | | | 4 | 17 | 13 | 13 | 10 | 3 | | 8 | 4 | 10 | | | | 3 |
| 4 | | | 12 | | 17 | | 8 | 4 | 17 | 8 | | 10 | | | | | 10 | | | | 4 |
| 5 | 18 | 5 | | | | | 8 | 4 | 5 | 18 | | 10 | | 20 | 9 | 20 | 10 | 5 | | 9 | 5 |
| 6 | 13 | 20 | 6 | 6 | 16 | | 8 | 4 | 11 | 8 | 13 | 11 | 6 | 15 | 8 | 4 | | | 20 | 11 | 6 |
| 7 | | | 12 | | 17 | | | 11 | 17 | 7 | 13 | 10 | 3 | | | | 10 | 3 | 1 | | 7 |
| 8 | 7 | 15 | 12 | | 17 | | 8 | | 17 | 7 | | 10 | | 15 | 8 | | 10 | | | | 8 |
| 9 | | 15 | | | 12 | 17 | | 11 | | | 9 | 11 | | 15 | 9 | 20 | 11 | | 20 | 9 | 9 |
| 10 | | | 12 | 12 | 17 | | | | 12 | | | 10 | | | | | 10 | | | | 10 |
| 11 | | 15 | | 12 | | | 8 | 11 | 12 | 8 | | 10 | | | 8 | | 10 | 5 | 20 | 9 | 11 |
| 12 | 13 | 15 | 6 | 6 | 12 | 12 | | 4 | 5 | 13 | | | 6 | 15 | 8 | 6 | | 6 | 20 | 9 | 12 |
| 13 | 13 | 15 | | 1 | 17 | | | 4 | 17 | 13 | 13 | | | 15 | 8 | 4 | 10 | | 1 | | 13 |
| 14 | | 15 | 14 | 12 | 17 | | | 4 | | 8 | | 10 | 3 | | 14 | 4 | 10 | 3 | | | 14 |
| 15 | 18 | 15 | 12 | | 2 | 17 | | 11 | 17 | 13 | 2 | 11 | 5 | | | | 11 | 5 | | 11 | 15 |
| 16 | 7 | 5 | 12 | | 17 | | 16 | | 17 | 7 | | 10 | 5 | | | | 10 | 5 | | 11 | 16 |
| 17 | 13 | 5 | 6 | 6 | 16 | 17 | 16 | | 5 | | 13 | | 6 | | 9 | 6 | | 6 | | 9 | 17 |
| 18 | 18 | | 12 | 12 | 17 | | | 4 | 17 | 18 | 18 | 10 | | | | 4 | 10 | | 1 | | 18 |
| 19 | | | | | | | | | | | | | | | | | | | | | 19 |
| 20 | 18 | 20 | 12 | | 16 | 17 | 16 | 11 | | 18 | 2 | 11 | 5 | | | 20 | 11 | 5 | | | 20 |
| | 1 | 2 | 3 | 4 | 5 | 6 | 7 | 8 | 9 | 10 | 11 | 12 | 13 | 14 | 15 | 16 | 17 | 18 | 19 | 20 | |

$\sigma_1^{-1}M^{-}(20)$

| | 3 | 9 | 1 | 7 | 11 | 10 | 4 | 16 | 2 | 6 | 5 | 17 | 18 | 19 | 20 | 8 | 12 | 13 | 14 | 15 | |
|---|---|---|---|---|---|---|---|---|---|---|---|---|---|---|---|---|---|---|---|---|---|
| | 1 | 2 | 3 | 4 | 5 | 6 | 7 | 8 | 9 | 10 | 11 | 12 | 13 | 14 | 15 | 16 | 17 | 18 | 19 | 20 | |
| 1 | **0** | -2 | ∞ | 2 | -5 | | | | -2 | 2 | | -4 | | | -3 | | -6 | | 9 | | 1 |
| 2 | 4 | **0** | | | 2 | 9 | | | ∞ | 2 | 0 | 5 | 3 | | | | 8 | 3 | | | 2 |
| 3 | ∞ | -3 | **0** | 1 | -6 | | | | -3 | 1 | 5 | -5 | 2 | | -4 | 0 | -7 | 3 | | | 3 |
| 4 | | | | **0** | -7 | | ∞ | -4 | -4 | 0 | | -6 | | | | -2 | -8 | | | | 4 |
| 5 | 2 | 3 | | 6 | **0** | | 4 | 3 | -1 | 0 | ∞ | 1 | 7 | 4 | 0 | 0 | -1 | 1 | | 1 | 5 |
| 6 | -6 | 1 | -4 | -8 | -11 | **0** | | -12 | 3 | ∞ -8 | -5 | -14 | -8 | -7 | -13 | -10 | -16 | -6 | 4 | 11 | 6 |
| 7 | -3 | | -5 | ∞ | -11 | | **0** | 5 | -8 | -4 | 3 | -10 | 0 | | | | -12 | 1 | 6 | | 7 |
| 8 | -2 | 0 | -4 | | -10 | | 1 | **0** | -7 | -3 | | -9 | | 5 | -1 | ∞ | -11 | | | | 8 |
| 9 | | ∞ | | | 5 | 6 | | 5 | **0** | | 4 | 5 | | 7 | 1 | 1 | 5 | | 5 | 2 | 9 |
| 10 | | | -1 | 2 | -7 | ∞ | | | 1 | **0** | | -6 | | | | | -8 | | | | 10 |
| 11 | | | | | ∞ | | | 1 | 1 | 5 | **0** | -1 | | | 0 | | -3 | | 6 | 3 | 11 |
| 12 | -7 | 2 | -5 | -9 | 0 | -1 | | 2 | -1 | -8 | | **0** | -9 | 7 | 1 | -4 | ∞ | -7 | 4 | 1 | 12 |
| 13 | 2 | 0 | | 4 | -6 | | | | -3 | 1 | 3 | -5 | **0** | 5 | -1 | 2 | -7 | ∞ | 11 | | 13 |
| 14 | | 4 | 8 | 9 | 2 | | | 5 | 5 | 9 | | 3 | 10 | **0** | 3 | 7 | 1 | 11 | ∞ | | 14 |
| 15 | 5 | 1 | 7 | | 3 | 3 | | 2 | | -3 | 1 | 6 | 4 | 6 | **0** | | 2 | 4 | | ∞ | 15 |
| 16 | -4 | 4 | -6 | | -12 | | -1 | ∞ | -9 | -5 | | -11 | | | | **0** | -13 | 2 | | | 16 |
| 17 | -6 | 4 | -3 | -7 | -1 | 1 | -3 | | 0 | | -5 | ∞ | -8 | | | -2 | **0** | -5 | | | 17 |
| 18 | 1 | | -2 | 1 | -8 | | | 8 | -5 | -1 | 7 | -7 | ∞ | | | | -9 | **0** | 10 | | 18 |
| 19 | | | | | | | | | | | | | | ∞ | | | | | **0** | 6 | 19 |
| 20 | 4 | 0 | 10 | | 0 | 9 | -2 | | | 2 | 0 | 5 | 3 | | ∞ | -1 | 8 | 3 | 3 | **0** | 20 |
| | 1 | 2 | 3 | 4 | 5 | 6 | 7 | 8 | 9 | 10 | 11 | 12 | 13 | 14 | 15 | 16 | 17 | 18 | 19 | 20 | |

**C followed by a set 2-cycles (boldface) represents 2-cycles each of which contains a point in a given cycle.**

**DC followed by a set of 2-cycles represents 2-cycles none of whose points occur in the given cycle.**

**The set of 2-cycles was given earlier.**

## ACCEPTABLE CYCLES

**1. C 2, 4; 5; (11 9): 5. 2. C 9, 10; 9; (14 15). 3. C 2, 4, 10; 9; (5 9 15). 4. C 2, 4, 7; 9; (9 11 17).**

**5. C 2, 4, 10; 9; (20 2 11). 6. C 3, 5, 6, 7; 2; (4 8 10 12). 7. C 1, 3, 5, 6, 7; 0; (3 4 8 10 12).**

**8. DC 2, 8, 9, 10; 1; (5 1 4 8 10 12).**

## 2-CIRCUIT UNLINKED PATHS

**1. C 5, 8; -7; [*6* 18 *10*]. 2. C 5, 7; -7; [*10* 17 *6*]. 3. C 3, 6; -3; [*4* 8 *7*]. 4. C 1, 3; -1; [*7* 1 *4*].**

**5. C 2, 4; 1; [*2* 5 *9*]. 6. C 2, 10; 2; [*9* 15 2]. 7. C 4, 7; 2; [*11* 17 *5*]. 8. C 1, 8; 3; [*3* 13 *1*].**

**9. C 2, 4; 3; [*5* 2 *11*]; 10. C 2, 4; 3; [*5* 9 *11*]. 11. C 3, 5, 6; -8; [*6* 4 8 *10*].**

**12. C 3, 5, 7; -2; [*7* 10 12 *4*]. 13. C 2, 4, 10; 10; [*15* 2 11 *20*]. 14. C 1, 3, 5, 6, 7; 1; [*1* 4 8 10 12 *3*].**

We now obtain paths used to construct $\sigma_1^{-1}M^{-}(40)$.

j = 1

(2 1)(1 4) = (2 4): 6; (5 1)(1 4) = (5 4): 4; (5 1)(1 19) = (5 19): 11; (6 1)(1 17) = (6 17): 0;

(8 1)(1 19) = (8 19): 7; (12 1)(1 19) = 2; (15 1)(1 4) = (15 4): 7; (15 1)(1 17) = (15 17): 10

(15 1)(1 19) = (15 19): 10; (16 1)(1 19) = (16 19): 5;

**(17 1)(1 17) = (17 17): 0**

**CYCLE P = [17 10 6 18 13 3 1 12 17]: (17 6 13 1): 0**

(17 1)(1 19) = (17 19): 3; (20 1)(1 4) = (20 4): 6

j = 2

(1 2)(2 11) = (1 11): -2; (3 2)(2 11) = (3 11): -3; (8 2)(2 11) = (8 11): 0;

(12 2)(2 11) = (12 11): 2;

**(12 2)(2 12) = (12 12): 7**

**CYCLE P = [12 7 4 16 8 20 15 9 2 17 12]: (12 4 8 15 2): 7**

(13 2)(2 11) = (13 11): 0; (17 2)(2 12) = (17 12): 7

**(17 2)(2 17) = (17 17): 10**

**CYCLE P = [17 10 6 8 16 11 5 9 2 12 17]: (17 6 16 5 2): 10**

j = 3

(7 3)(3 11) = (7 11): 1; (7 3)(3 13) = (7 13): -3; (7 3)(3 18) = (7 18): -2;

(8 3)(3 13) = (8 13): -2; (8 3)(3 18) = (8 18): -3; (10 3)(3 4) = (10 4): 0;

(10 3)(3 11) = (10 11): 5; (10 3)(3 13) = (10 13): 0; (10 3)(3 18) = (10 18): 2;

(12 3)(3 11) = (12 11): 1; (16 3)(3 11) = (16 11): 0; (16 3)(3 13) = (16 13): -4;

(16 3)(3 18) = (16 18): -3; (18 3)(3 4) = (18 4): -1; (18 3)(3 11) = (18 11): 4;

**(18 3)(3 18) = (18 18): 1**

**CYCLE P = [18 6 10 17 12 1 3 13 18]: (18 10 12 3): 1**

j = 4

(2 4)(4 8) = (2 8): 2; (2 4)(4 16) = (2 16): 4; (5 4)(4 8) = (5 8): 0; (10 4)(4 8) = (10 8): -4;

(10 4)(4 16) = (10 16): -2; (12 4)(4 8) = (12 8): -13; (12 4)(4 15) = (12 15): 3;

(12 4)(4 16) = (12 16): -11; (15 4)(4 16) = (15 16): 5; (17 4)(4 8) = (17 8): -11;

(17 4)(4 15) = (17 15): 5; (17 4)(4 16) = (17 16): -9;

**(17 4)(4 17) = (17 17): 6**

**CYCLE P = [17 10 6 7 4 12 17]: (17 6 4): 6**

(15 4)(4 8) = (15 8): -3; (15 4)(4 16) = (15 16): -1; (18 4)(4 8) = (18 8): -3;

(18 4)(4 16) = (18 16): -1

j = 5

**(1 5)(5 1) = (1 1): 0**

**CYCLE P = [1 7 4 16 8 6 10 12 17 11 5 3 1]: (1 4 8 10 17 5): 0**

(3 5)(5 9) = (3 9): -7; (4 5)(5 1) = (4 1): -2; (4 5)(5 2) = (4 2): -4; (4 5)(5 9) = (4 9): -8;

(4 5)(5 13) = (4 13): 0; (4 5)(5 2) = (4 2): -4; (4 5)(5 15) = (4 15): 11; (4 5)(5 18) = (4 18): -6;

(6 5)(5 2) = (6 2): -8; (6 5)(5 9) = (6 9): -12; (6 5)(5 18) = (6 18): -8; (7 5)(5 1) = (7 1): -6;

(7 5)(5 2) = (7 2): -8; (7 5)(5 9) = (7 9): -12; (7 5)(5 13) = (7 13): -4; (7 5)(5 15) = (7 15): 7;

(7 5)(5 18) = -10; (8 5)(5 1) = (8 1): -5; (8 5)(5 2) = (8 2): -7; (8 5)(5 9) = (8 9): -11;

(8 5)(5 13) = (8 13): -3; (8 5)(5 18) = (8 18): -9; (9 5)(5 1) = (9 1): 10; (10 5)(5 1) = (10 1): -2;

(10 5)(5 2) = (10 2): -4; (10 5)(5 9) = (10 9): -8; (10 5)(5 15) = (10 15): 11;

(10 5)(5 18) = (10 18): -6; (13 5)(5 2) = (13 2): -3; (13 5)(5 9) = (13 9): -7;

**(13 5)(5 13) = (13 13): 1**

**CYCLE P = [13 <u>3</u> 1 <u>7</u> 4 <u>16</u> 8 <u>6</u> 10 <u>12</u> 17 <u>11</u> 5 <u>18</u> 13]: (13 1 4 8 10 17 5): 1**

(14 5)(5 2) = (14 2): <u>4</u>; (14 5)(5 9) = (14 9): 1; (14 5)(5 13) = (14 13): 9;

(14 5)(5 18) = (14 18): 3; (16 5)(5 1) = (16 1): <u>-7</u>; (16 5)(5 2) = (16 2): <u>-9</u>;

(16 5)(5 9) = (16 9): -13; (16 5)(5 13) = (16 13): -5; (16 5)(5 15) = (16 15): 6;

(16 5)(5 18) = (16 18): -11; (17 5)(5 2) = (17 2): 2; (17 5)(5 9) = (17 9): -2;

(18 5)(5 1) = (18 1): <u>-3</u>; (18 5)(5 2) = (18 2): <u>-5</u>; (18 5)(5 9) = (18 9): -9;

(18 5)(5 15) = (18 15): 10;

**(18 5)(5 18) = (18 18): 10**

**CYCLE P = [18 <u>6</u> 10 <u>12</u> 17 <u>11</u> 5 <u>13</u> 18]: (18 10 17 5): -7**

(20 5)(5 9) = (20 9): -1; (20 5)(5 18) = (20 18): 1

j = 6

(9 6)(6 7) = (9 7): 7; (9 6)(6 13) = (9 13): 9; (9 6)(6 18) = (9 18): 7; (15 6)(6 3) = (15 3): <u>-1</u>;

(15 6)(6 4) = (15 4): <u>-5</u>; (15 6)(6 13) = (15 13): -5; (15 6)(6 16) = (15 16): 0;

(15 6)(6 18) = (15 18): -3

j = 7

(9 7)(7 1) = (9 1): <u>4</u>; (9 7)(7 3) = (9 3): <u>-2</u>; (17 7)(7 1) = (17 1): <u>-6</u>; (17 7)(7 3) = (17 3): <u>-5</u>;

**(17 7)(7 17) = (17 17): 2**

**CYCLE P = [17 <u>10</u> 6 <u>8</u> 16 <u>4</u> 7 <u>12</u> 17]: (17 6 16 7): 2**

(20 7)(7 1) = (20 1): <u>-5</u>; (20 7)(7 3) = (20 3): <u>-4</u>; (20 7)(7 10) = (20 10): -6

j = 8

(2 8)(8 15) = (2 15): 1; (10 8)(8 4) = (10 4): <u>0</u>; (10 8)(8 15) = (10 15): -5;

**(10 8)(8 10) = (10 10): 0**

**CYCLE P = [10 <u>17</u> 12 <u>7</u> 4 <u>16</u> 8 <u>6</u> 10]: (10 12 4 8): 0**

(15 8)(8 7) = (15 7): <u>3</u>;

**(15 8)(8 15) = (15 15): -4**

**CYCLE P = [15 <u>11</u> 5 <u>13</u> 18 <u>3</u> 1 <u>7</u> 4 <u>16</u> 8 <u>20</u> 15]: (15 5 18 1 4 8): -4**

(17 8)(8 15) = (17 15): -12; (18 8)(8 2) = (18 2): <u>9</u>; (18 8)(8 15) = (18 15): -4

j = 9

(1 9)(9 16) = (1 16): 11; (1 9)(9 20) = (1 20): 0; (4 9)(9 11) = (4 11): -4;

(4 9)(9 20) = (4 20): -6; (7 9)(9 11) = (7 11): -8; (7 9)(9 20) = (7 20): -10;

(7 9)(9 15) = (7 15): -11; (7 9)(9 20) = (7 20): -10; (8 9)(9 11) = (8 11): -3;

(8 9)(9 15) = (8 15): -6; (8 9)(9 20) = (8 20): -9; (10 9)(9 11) = (10 11): -4;

(10 9)(9 15) = (10 15): -7; (10 9)(9 20) = (10 20): -6; (13 9)(9 15) = (13 15): -6;

(13 9)(9 20) = (13 20): -5; (14 9)(9 11) = (14 11): 5; (14 9)(9 15) = (14 15): 2;

(14 9)(9 20) = (14 20): 3; (16 9)(9 15) = (16 15): -12; (16 9)(9 20) = (16 20): -11;

(17 9)(9 15) = (17 15): -1; (17 9)(9 20) = (17 20): 0; (18 9)(9 15) = (18 15): -8;

(18 9)(9 20) = (18 20): -7

j = 10

(1 10)(10 12) = (1 12): -4; (1 10)(10 17) = (1 17): -6; (2 10)(10 12) = (2 12): -4;

(2 10)(10 17) = (2 17): -6; (3 10)(10 17) = (3 17): -5; (5 10)(10 12) = (5 12): -6;

(5 10)(10 17) = (5 17): -8;

**(12 10)(10 12) = (12 12): 2**

**CYCLE P = [12 11 5 18 13 6 10 17 12]: (12 5 13 10): 2**

j = 11

(3 11)(11 19) = (3 19): 9; (7 11)(11 8) = (7 8): -7; (7 11)(11 19) = (7 19): 4;

**(8 11)(11 8) = (8 8): -4**

**CYCLE P = [8 4 7 6 10 12 17 2 9 5 11 16 8]: (8 7 10 17 9 11): -4**

(8 11)(11 20) = (8 20): 4;

**(12 11)(11 12) = (12 12): 2**

**CYCLE P = [12 10 6 1 3 7 4 16 8 20 15 9 2 5 11 17 12]: (12 6 3 4 8 15 2 11): 2**

(17 11)(11 9) = (17 9): -4;

**(17 11)(11 17) = (17 17): -4**

**CYCLE P = [17 10 6 18 13 5 11 12 17]: (17 6 13 11): -4**

j = 12

(2 12)(12 3) = (2 3): 1; (2 12)(12 4) = (2 4): 4; (5 12)(12 3) = (5 3): -1;

**(5 12)(12 5) = (5 5): -4**

**CYCLE P = [5 13 18 6 10 17 12 11 5]: (5 18 10 12): -4**

**(6 12)(12 6) = (6 6): -8**

**CYCLE P = [6 18 13 5 11 17 12 10 6]: (6 13 11 12): -8**

(15 12)(12 3) = (15 3): 7

j = 13

**(10 13)(13 10) = (10 10): 1**

**CYCLE P = [10 17 12 1 3 18 13 6 10]: (10 12 3 13): 1**

(15 13)(13 1) = (15 1): -3; (15 13)(13 10) = (15 10): -4; (20 13)(13 1) = (20 1): -8;

(20 13)(13 14) = (20 14): 11

j = 15

(2 15)(15 14) = (2 14): 7; (3 15)(15 14) = (3 14): 0;

**(2 15)(15 2) = (2 2): 2**

**CYCLE P = [2 11 5 13 18 3 1 7 4 16 8 20 15 9 2]: (2 5 18 1 4 8 15): 2**

(6 15)(15 14) = (6 14): -5; (7 15)(15 14) = (7 14): -5; (8 15)(15 5) = (8 5): 1;

(8 15)(15 14) = (8 14): -2; (10 15)(15 14) = (10 14): -1;

**(14 15)(15 14) = (14 14): 8**

**CYCLE P = [14 1 3 7 4 16 8 6 10 12 17 11 5 2 9 20 15 19 14]: (14 3 4 8 10 17 5 9 15): 8**

(16 15)(15 14) = (16 14): -6; (17 15)(15 2) = (17 2): -11; (17 15)(15 14) = (17 14): -6;

(18 15)(15 14) = (18 14): -2

j = 16

**(5 16)(16 5) = (5 5): 1**

**CYCLE P = [5 2 9 15 20 8 16 11 5]: (5 9 20 16): 1**

(5 16)(16 7) = (5 7): -1; (9 16)(16 5) = (9 5): 2; (9 16)(16 7) = (9 7): 0;

(10 16)(16 7) = (10 7): -8; (17 16)(16 5) = (17 5): -8

j = 17

(3 17)(17 5) = (3 5): -4;

**(5 17)(17 5) = (5 5): -7**

**CYCLE P = [5 13 18 6 10 12 17 11 5]: (5 18 10 17): -7**

(5 17)(17 9) = (5 9): -7

j = 18

(4 18)(18 1) = (4 1): -5; (4 18)(18 15) = (4 15): 8; (6 18)(18 1) = (6 1): -7;

(7 18)(18 1) = (7 1): -9; (8 18)(18 1) = (8 1): -8; (10 18)(18 1) = (10 1): -5;

**(10 18)(18 10) = (10 10): -7**

**CYCLE P = [10 12 17 11 5 13 18 6 10]: (10 17 5 18): -7**

(15 18)(18 1) = (15 1): 3; (15 18)(18 10) = (15 10): 1; (16 18)(18 1) = (16 1): -10

j = 19

(6 19)(19 20) = (6 20): 9

j = 20

(1 20)(20 16) = (1 16): -1; (1 20)(20 19) = (1 19): 3; (3 20)(20 7) = (3 7): 10;

(3 20)(20 16) = (3 16): -6; (3 20)(20 19) = (3 19): 0;

**(7 20)(20 7) = (7 7): 5**

**CYCLE P = [7 6 10 12 17 11 5 2 9 15 20 4 7]: (7 10 17 5 9 20): 5**

(7 20)(20 16) = (7 16): -11; (7 20)(20 19) = (7 19): -7; (8 20)(20 19) = (8 19): -6;

(10 20)(20 16) = (10 16): -7; (10 20)(20 19) = (10 19): -3;

**(16 20)(20 16) = (16 16): -12**

**CYCLE P = [16 4 7 6 10 12 17 11 5 2 9 15 20 8 16]: (16 7 10 17 5 9 20): -12**

(16 20)(20 19) = (16 19): -8; (17 20)(20 19) = (17 19): 0; (18 20)(20 16) = (18 16): -8;

(18 20)(20 19) = (18 19): 0

**ACCEPTABLE CYCLES**

**1. DC 8; 8; (14 3 4 8 10 17 5 9 15). 2. DC 8, 9; 2; (12 6 3 4 8 15 2 11).**

**3. DC 1, 8, 9; -12; (16 7 10 17 5 9 20). 4. DC 2, 9, 10; 1; (13 1 4 8 10 17 5).**

**5. DC 5, 7, 9; 2; (2 5 18 1 4 8 15). 6. DC 1, 8, 9, 10; -4; (8 7 10 17 9 11).**

**7. DC 2, 5, 7, 9; -4; (15 5 18 1 4 8). 8. DC 2, 8, 9, 10; 0; (1 4 8 10 17 5).**

**9. DC 1, 5, 6, 8, 9; 5; (7 10 17 5 9 20). 10. C 2, 4, 5, 6, 7; 10; (17 6 16 5 2).**

**11. C 2, 3, 6, 7, 10; 7; (12 4 8 15 2); 12. C 4, 5, 7, 8; -7; (18 10 17 5). 13. C 1, 5, 7, 8; 0; (17 16 13 1).**

**14. C 2, 4, 6, 10; 1; (5 9 20 16). 15. C 4, 5, 7, 8; 2; (12 5 13 10). 16. C 3, 5, 7, 6; (17 6 16 7).**

**17. C 3, 5, 7; 6; (17 6 4).**

**2-CIRCUIT UNLINKED PATHS**

**1. DC 1, 4, 8, 9; -5; [*12* 6 4 8 15 2 *17*]. 2. C 2, 4, 5, 7, 8; 6; [*2* 5 13 10 17 *9*].**

**3. C 2, 3, 5, 6, 7; 6; [*8* 7 10 17 9 *16*]. 4. C 2, 4, 5, 6, 7; 7; [*17* 6 16 5 2 *12*].**

**5. C 1, 5, 7, 8; -1; [*12* 6 13 1 *17*]. 6. C 1, 3, 5, 7; -4; [*7* 10 12 3 *4*]. 7. C 4, 5, 7 ,8; -3; [*17* 6 13 11 *12*];**

**8. C 4, 5, 7, 8; -1; [*18* 10 17 5 *13*]. 9. C 1, 5, 7, 8; 0; [*18* 10 12 3 *13*]. 10. C 2, 4, 6, 10; 0; [*20* 16 5 9 *15*].**

**11. C 4, 5, 7, 8; 10; [*5* 18 10 17 *11]*. 12. C 1, 6, 8; -2; [*6* 13 1 *10*].**

$P_{40}$

| | 1 | 2 | 3 | 4 | 5 | 6 | 7 | 8 | 9 | 10 | 11 | 12 | 13 | 14 | 15 | 16 | 17 | 18 | 19 | 20 | |
|---|---|---|---|---|---|---|---|---|---|---|---|---|---|---|---|---|---|---|---|---|---|
| 1 | | 15 | 12 | 1 | 17 | | | 4 | 17 | 8 | 2 | 10 | | | 8 | 20 | 10 | | 20 | 9 | 1 |
| 2 | 18 | | 12 | 12 | 2 | 17 | | 4 | | 18 | | 10 | 5 | 15 | 8 | 4 | 10 | 5 | 1 | | 2 |
| 3 | | 15 | | 3 | 17 | | | 4 | 5 | 13 | 2 | 10 | 3 | 15 | 8 | 20 | 10 | | 20 | 9 | 3 |
| 4 | 18 | 5 | 12 | | 17 | | 8 | 4 | 5 | 8 | 9 | 10 | 5 | | 18 | | 10 | 5 | | 9 | 4 |
| 5 | 18 | 5 | 12 | 1 | | | 16 | 4 | 17 | 18 | | 10 | | 20 | 9 | 20 | 10 | 5 | 1 | 9 | 5 |
| 6 | 18 | 5 | 6 | 6 | 16 | | 8 | 4 | 5 | 8 | 13 | 11 | 6 | 15 | 8 | 4 | 1 | 5 | 20 | 11 | 6 |
| 7 | 18 | 5 | 12 | | 17 | | | 11 | 5 | 7 | 9 | 10 | 5 | 15 | 9 | 20 | 10 | 5 | 20 | 9 | 7 |
| 8 | 18 | 5 | 12 | | 17 | | 8 | | 5 | 7 | 9 | 10 | 5 | 15 | 9 | | 10 | 5 | 20 | 9 | 8 |
| 9 | 18 | 15 | 7 | | 16 | 17 | 16 | 11 | | | 9 | 11 | 6 | 15 | 9 | 20 | 11 | 6 | 20 | 9 | 9 |
| 10 | 18 | 5 | 12 | 8 | 17 | | 16 | 4 | 5 | | 3 | 10 | 3 | 15 | 9 | 20 | 10 | 5 | 20 | 9 | 10 |
| 11 | | 15 | | 12 | | | 8 | 11 | 12 | 8 | | 10 | | | 8 | | 10 | 5 | 20 | 9 | 11 |
| 12 | 13 | 15 | 6 | 6 | 12 | 12 | | 4 | 5 | 13 | 3 | | 6 | 15 | 8 | 4 | | 6 | 1 | 9 | 12 |
| 13 | 13 | 5 | | 1 | 17 | | | 4 | 5 | 13 | 2 | | | 15 | 9 | 4 | 10 | | 1 | 9 | 13 |
| 14 | | 5 | 14 | 12 | 17 | | | 4 | 5 | 8 | 9 | 10 | 5 | | 9 | 4 | 10 | 5 | | 9 | 14 |
| 15 | 13 | 15 | 6 | 6 | 2 | 17 | 16 | 11 | 17 | 13 | 2 | 11 | 6 | | | 6 | 11 | 6 | | 11 | 15 |
| 16 | 18 | 5 | 12 | | 17 | | 16 | | 5 | 7 | 3 | 10 | 5 | 15 | 8 | | 10 | 5 | 1 | 9 | 16 |
| 17 | 7 | 15 | 7 | 6 | 16 | 17 | | 4 | 11 | | 13 | 2 | 6 | 15 | 9 | 4 | | 6 | 1 | 9 | 17 |
| 18 | 5 | 8 | 12 | 3 | 17 | | | 4 | 5 | 18 | 3 | 10 | | 15 | 9 | 20 | 10 | | 1 | 9 | 18 |
| 19 | | | | | | | | | | | | | | | | | | | | | 19 |
| 20 | 13 | 20 | 7 | 1 | 16 | 17 | 16 | 11 | 5 | 7 | 2 | 11 | 5 | 13 | | 20 | 11 | 5 | | | 20 |
| | 1 | 2 | 3 | 4 | 5 | 6 | 7 | 8 | 9 | 10 | 11 | 12 | 13 | 14 | 15 | 16 | 17 | 18 | 19 | 20 | |

$\sigma_1^{-1} M^{-}(40)$

| | 3 | 9 | 1 | 7 | 11 | 10 | 4 | 16 | 2 | 6 | 5 | 17 | 18 | 19 | 20 | 8 | 12 | 13 | 14 | 15 | |
|---|---|---|---|---|---|---|---|---|---|---|---|---|---|---|---|---|---|---|---|---|---|
| | 1 | 2 | 3 | 4 | 5 | 6 | 7 | 8 | 9 | 10 | 11 | 12 | 13 | 14 | 15 | 16 | 17 | 18 | 19 | 20 | |
| 1 | **0** | -2 | ∞ | 2 | -5 | | | | -2 | 2 | -2 | -4 | | | -3 | -1 | -6 | | 3 | 0 | 1 |
| 2 | 4 | **0** | 1 | 4 | 2 | 9 | | 2 | ∞ | 2 | 0 | -4 | 3 | 7 | 1 | 4 | -6 | 3 | | | 2 |
| 3 | ∞ | -3 | **0** | 1 | -4 | | | | -3 | 1 | -3 | -5 | 2 | 2 | -6 | -6 | -5 | 3 | 0 | -5 | 3 |
| 4 | -5 | | | **0** | -7 | | ∞ | -4 | -8 | 0 | -4 | -6 | 0 | | 8 | -2 | -8 | -6 | -3 | -6 | 4 |
| 5 | 2 | 3 | -1 | 4 | **0** | | -1 | 3 | -1 | 0 | ∞ | -4 | 7 | 4 | 0 | 0 | -8 | 1 | 11 | 1 | 5 |
| 6 | -7 | -12 | -4 | -8 | -11 | **0** | | -12 | -12 | ∞ | -5 | -7 | -8 | -7 | -13 | -10 | -1 | -10 | 3 | 9 | 6 |
| 7 | -9 | -8 | -5 | ∞ | -11 | | **0** | -7 | -12 | -4 | -8 | -10 | -4 | | -11 | -11 | -12 | -10 | -7 | -10 | 7 |
| 8 | -8 | -7 | -4 | | -10 | | 1 | **0** | -9 | -3 | -5 | -9 | -3 | 5 | -8 | ∞ | -11 | -9 | -6 | -9 | 8 |
| 9 | 4 | ∞ | 3 | | 2 | 6 | 0 | 5 | **0** | | 4 | 5 | 6 | 7 | 1 | 1 | 5 | 7 | 5 | 2 | 9 |
| 10 | -5 | -4 | -1 | 0 | -7 | ∞ | -8 | -4 | -8 | **0** | -4 | -6 | 0 | -1 | -7 | -7 | -8 | -6 | -8 | -6 | 10 |
| 11 | | | | | ∞ | | | 1 | 1 | 5 | **0** | -1 | | | 0 | | -3 | | 6 | 3 | 11 |
| 12 | -7 | 2 | -5 | -9 | 0 | -1 | | 2 | -1 | -8 | 1 | **0** | -9 | 7 | 1 | -4 | ∞ | -7 | 2 | 1 | 12 |
| 13 | 2 | -3 | | 4 | -6 | | | | -7 | 1 | 0 | -5 | **0** | 5 | -6 | 2 | -7 | ∞ | -2 | -5 | 13 |
| 14 | | 4 | 8 | 9 | 2 | | | 5 | 1 | 9 | 4 | 3 | 9 | **0** | 2 | 7 | 1 | 3 | ∞ | 3 | 14 |
| 15 | -3 | 1 | -1 | -5 | 3 | 3 | 3 | 2 | | -3 | 3 | -3 | -5 | 6 | **0** | 0 | 2 | -3 | | ∞ | 15 |
| 16 | -10 | -9 | -6 | | -12 | | -1 | ∞ | -13 | -5 | 0 | -11 | -5 | -6 | -12 | **0** | -13 | -11 | -8 | -11 | 16 |
| 17 | -6 | -11 | -3 | -7 | -8 | 1 | | -11 | 0 | | -5 | ∞ | -8 | -6 | -12 | -9 | **0** | -5 | 0 | -3 | 17 |
| 18 | 1 | 9 | -2 | 1 | -8 | | | -3 | -9 | -1 | 4 | -7 | ∞ | -2 | -8 | -8 | -9 | **0** | 0 | -3 | 18 |
| 19 | | | | | | | | | | | | | | ∞ | | | | | **0** | 6 | 19 |
| 20 | -5 | 0 | -4 | 1 | 0 | 9 | -2 | | -1 | -6 | 0 | 5 | 1 | 11 | ∞ | -1 | 3 | 1 | 3 | **0** | 20 |
| | 1 | 2 | 3 | 4 | 5 | 6 | 7 | 8 | 9 | 10 | 11 | 12 | 13 | 14 | 15 | 16 | 17 | 18 | 19 | 20 | |

j = 1

**(4 1)(1 4) = (4 4): -3**

**CYCLE P = [4 16 8 6 10 12 17 11 5 13 18 3 1 7 4]: (4 8 10 17 5 13 1): -3**

(6 1)(1 17) = (6 17): -3; (6 1)(1 19) = (6 19): 0; (9 1)(1 4) = (9 4): 6; (10 1)(1 4) = (10 4): -3;

**(10 1)(1 10) = (10 10): -1**

**CYCLE P = [10 12 17 11 5 13 18 3 1 6 10]: (10 17 5 18 1): -1**

(20 1)(1 17) = (20 17): 1

j = 2

(17 2)(2 5) = (17 5): -9; (17 2)(2 11) = (17 11): -11;

(17 2)(2 17) = (17 17): -3

CYCLE P = [17 10 6 7 4 16 8 20 15 9 2 12 17]: (17 6 4 8 15 2): -3

j = 3

(2 3)(3 4) = (2 4): 2; (5 3)(3 4) = (5 4): 0; (9 3)(3 13) = (9 13): 0;

(9 3)(3 18) = (9 18): 1; (20 3)(3 13) = (20 13): -2; (20 3)(3 18) = (20 18): -1

j = 4

(2 4)(4 8) = (2 8): 0; (2 4)(4 16) = (2 16): 2; (5 4)(4 8) = (5 8): -4; (10 4)(4 8) = (10 8): -7;

(15 4)(4 8) = (15 8): 1; (15 4)(4 16) = (15 16): 3

j = 5

(3 5)(5 9) = (3 9): -5;

**(9 5)(5 9) = (9 9): 1**

**CYCLE P = [9 15 20 8 16 11 5 2 9]: (9 20 16 5): 1**

(17 5)(5 18) = (17 18): -8

j = 7

(9 7)(7 1) = (9 1): -3; (9 7)(7 3) = (9 3): -2; (9 7)(7 10) = (9 10): -4; (10 7)(7 1) = (10 1): -11;

(10 7)(7 3) = (10 3): -10;

**(10 7)(7 10) = (10 10): -12**

**CYCLE P = [10 12 17 11 5 2 9 15 20 8 16 4 7 6 10]: (10 17 5 9 20 16 7): -12**

(17 7)(7 1) = (17 1): -13; (17 7)(7 3) = (17 3): -12

j = 8

(2 8)(8 15) = (2 15); -1; (5 8)(8 15) = (5 15): -5;

**(10 8)(8 10) = (10 10): -3**

**CYCLE P = [10 12 17 11 5 13 18 3 1 7 4 16 8 6 10]: (10 17 5 18 1 4 8): -3**

(10 8)(8 15) = (10 15): -8;

**(15 8)(8 15) = (15 15): 2**

**CYCLE P = [15 9 2 11 5 13 18 3 1 7 4 16 8 20 15]: (15 2 18 1 4 8): 2**

j = 10

(9 10)(10 12) = (9 12): -10; (9 10)(10 17) = (9 17): -12; (20 10)(10 12) = (20 12): -12;

(20 10)(10 17) = (20 17): -14

j = 11

**(17 11)(11 17) = (17 17): -10**

**CYCLE P = [17 10 6 7 4 16 8 20 15 9 2 5 11 12 17]: (17 6 4 8 15 2 11): -10**

j = 12

(9 12)(12 3) = (9 3): -5; (9 12)(12 5) = (9 5): -10;

**(9 12)(12 9) = (9 9): -3**

**CYCLE P = [9 15 20 8 16 4 7 6 10 17 12 2 9]: (9 20 16 7 10 12): -3**

(20 12)(12 3) = (20 3): -7; (20 12)(12 5) = (20 5): -12; (20 12)(12 9) = (20 9): -5

j = 13

(20 13)(13 11) = (20 11): -2

j = 15

**(2 15)(15 2) = (2 2): 0**

**CYCLE P = [2 11 5 13 18 6 10 17 12 1 3 7 4 16 8 20 15 9 2]: (2 5 18 10 12 3 4 8 15): 0 ***

(5 15)(15 2) = (5 2): -4;

**(5 15)(15 5) = (5 5): -4**

**CYCLE P = [5 13 18 6 10 17 12 1 3 7 4 16 8 20 15 11 5]: (5 8 10 12 3 4 8 15): -4**

(10 15)(15 2) = (10 2): -7; (10 15)(15 14) = (10 14): -2;

j = 16

(3 16)(16 7) = (3 7): -6;

**(5 16)(16 5) = (5 5): -1**

**CYCLE P = [5 <u>13</u> 18 <u>6</u> 10 <u>17</u> 12 <u>1</u> 3 <u>7</u> 4 <u>8</u> 16 <u>11</u> 5]: (5 18 10 12 3 4 16): -1**

**(7 16)(16 7) = (7 7): -12**

**CYCLE P = [7 <u>6</u> 10 <u>12</u> 17 <u>11</u> 5 <u>2</u> 9 <u>15</u> 20 <u>8</u> 16 <u>4</u> 7]: (7 10 17 5 9 20 16): -12**

(18 16)(16 7) = (18 7): <u>-9</u>

j = 17

**(6 17)(17 6) = (6 6): 0**

**CYCLE P = [6 <u>7</u> 4 <u>8</u> 16 <u>11</u> 5 <u>13</u> 18 <u>3</u> 1 <u>12</u> 17 <u>10</u> 6]: (6 4 16 5 18 1 17): 0**

(9 17)(17 5) = (9 5): <u>-11</u>;

**(9 17)(17 9) = (9 9): -8**

**CYCLE P = [9 <u>15</u> 20 <u>8</u> 16 <u>4</u> 7 <u>6</u> 10 <u>12</u> 17 <u>2</u> 9]: (9 20 16 7 10 17): -8**

(20 17)(17 5) = (20 5): <u>-13</u>; (20 17)(17 9) = (20 9): <u>-10</u>

j = 18

**(17 18)(18 17) = (17 17): 7**

**CYCLE P = [17 <u>10</u> 6 <u>7</u> 4 <u>16</u> 8 <u>20</u> 15 <u>9</u> 2 <u>11</u> 5 <u>13</u> 18 <u>12</u> 17]: (17 6 4 8 15 2 5 18): 7**

j = 19

(4 19)(19 3) = (4 3): <u>9</u>; (6 19)(19 20) = (6 20): 8; (13 19)(19 3) = (13 3): <u>10</u>;

j = 20

(8 20)(20 19) = (8 19): <u>-6</u>

**ACCEPTABLE CYCLES**

**1. DC 9; 0; (2 5 18 10 12 3 4 8 15). 2. DC 1,9; 7; (17 6 4 8 15 2 5 18).**

**3. DC 1, 8, 9; -12; (10 17 5 9 20 16 7). 4. DC 1, 8, 9; -10; (17 6 4 8 15).**

**5. DC 2, 8, 9; -3; (4 8 10 17 5 15 1). 6. DC 2, 9, 10; -3; (10 17 5 18 1 4 8).**

**7. DC 2, 9, 10; 0; (6 4 16 5 18 1 17). 8. DC 5, 7, 9; 2; (15 2 5 18 1 4 8).**

**9. DC 1, 4, 8, 9; -8; (9 20 16 7 10 17). 10. DC 1, 4, 8, 9; (9 20 16 7 10 12).**

**11. C 1, 4, 5, 7, 8; 4; (10 17 5 18 1). 12. C 2, 4, 6, 10; 1; (9 20 16 5).**

**2-CIRCUIT UNLINKED PATHS**

**1. DC 1, 8, 9; [*17* 6 4 8 15 2 11 *12*]. 2. DC 2, 6, 9, 10; -7; [*7* 10 17 5 18 1 *4*].**

**3. DC 2, 6, 7, 9, 10; -5; [*6* 4 5 18 1 *10*]. 4. DC 2, 3, 6, 9, 10; 5; [*5* 18 10 12 3 *11*].**

**5. C 1, 4, 5, 7, 8; 5; [*9* 20 16 5 *2*].**

$P_{60}$

| | 1 | 2 | 3 | 4 | 5 | 6 | 7 | 8 | 9 | 10 | 11 | 12 | 13 | 14 | 15 | 16 | 17 | 18 | 19 | 20 | |
|---|---|---|---|---|---|---|---|---|---|---|---|---|---|---|---|---|---|---|---|---|---|
| 1 | | 15 | 12 | 1 | 17 | | | 4 | 17 | 8 | 2 | 10 | | | 8 | 20 | 10 | | 20 | 9 | 1 |
| 2 | 18 | | 12 | 12 | 2 | 17 | | 4 | | 18 | | 10 | 5 | 15 | 8 | 4 | 10 | 5 | 1 | | 2 |
| 3 | | 15 | | 3 | 17 | | 16 | 4 | 5 | 13 | 2 | 10 | 3 | 15 | 8 | 20 | 10 | | 20 | 9 | 3 |
| 4 | 18 | 15 | 19 | | 17 | | 8 | 4 | 5 | 8 | 9 | 10 | 5 | | 18 | | 10 | 5 | 20 | 9 | 4 |
| 5 | 18 | 15 | 12 | 3 | | | 16 | 4 | 17 | 18 | | 10 | | 15 | 8 | 20 | 10 | 5 | 1 | 9 | 5 |
| 6 | 18 | 5 | 6 | 6 | 16 | | 8 | 4 | 5 | 8 | 13 | 11 | 6 | 15 | 8 | 4 | 1 | 5 | 1 | 19 | 6 |
| 7 | 18 | 5 | 12 | | 17 | | | 11 | 5 | 7 | 9 | 10 | 5 | 15 | 9 | 20 | 10 | 5 | 20 | 9 | 7 |
| 8 | 18 | 5 | 12 | | 17 | | 8 | | 5 | 7 | 9 | 10 | 5 | 15 | 9 | | 10 | 5 | 20 | 9 | 8 |
| 9 | 7 | 15 | 12 | | 17 | 17 | 16 | 11 | | 7 | 9 | 10 | 3 | 15 | 9 | 20 | 10 | 3 | 20 | 9 | 9 |
| 10 | 18 | 15 | 7 | 1 | 17 | | 16 | 4 | 5 | | 3 | 10 | 3 | 15 | 8 | 20 | 10 | 5 | 20 | 9 | 10 |
| 11 | | 15 | | 12 | | | 8 | 11 | 12 | 8 | | 10 | | | 8 | | 10 | 5 | 20 | 9 | 11 |
| 12 | 13 | 15 | 6 | 6 | 12 | 12 | | 4 | 5 | 13 | 3 | | 6 | 15 | 8 | 4 | | 6 | 1 | 9 | 12 |
| 13 | 13 | 5 | 19 | 1 | 17 | | | 4 | 5 | 13 | 2 | | | 15 | 9 | 4 | 10 | | 1 | 9 | 13 |
| 14 | | 15 | 14 | 3 | 17 | | | 4 | 5 | 8 | 2 | 10 | 5 | | 14 | 4 | 10 | 5 | | 9 | 14 |
| 15 | 13 | 15 | 6 | 6 | 2 | 17 | 8 | 4 | 17 | 8 | 2 | 11 | 6 | | | 4 | 11 | 6 | | 11 | 15 |
| 16 | 18 | 5 | 12 | | 17 | | 16 | | 5 | 7 | 3 | 10 | 5 | 15 | 8 | | 10 | 5 | 20 | 9 | 16 |
| 17 | 7 | 15 | 7 | 6 | 2 | 17 | 16 | 4 | 11 | | 2 | 2 | 6 | 15 | 8 | 4 | | 6 | 1 | 9 | 17 |
| 18 | 5 | 8 | 12 | 3 | 17 | | 16 | 4 | 5 | 18 | 3 | 10 | | 15 | 9 | 20 | 10 | | 1 | 9 | 18 |
| 19 | | | | | | | | | | | | | | | | | | | | | 19 |
| 20 | 7 | 20 | 12 | 1 | 17 | 17 | 16 | 11 | 17 | 7 | 13 | 10 | 3 | 13 | | 20 | 10 | 3 | | | 20 |
| | 1 | 2 | 3 | 4 | 5 | 6 | 7 | 8 | 9 | 10 | 11 | 12 | 13 | 14 | 15 | 16 | 17 | 18 | 19 | 20 | |

$\sigma_1^{-1}M^-(60)$

| | 3 | 9 | 1 | 7 | 11 | 10 | 4 | 16 | 2 | 6 | 5 | 17 | 18 | 19 | 20 | 8 | 12 | 13 | 14 | 15 | |
|---|---|---|---|---|---|---|---|---|---|---|---|---|---|---|---|---|---|---|---|---|---|
| | 1 | 2 | 3 | 4 | 5 | 6 | 7 | 8 | 9 | 10 | 11 | 12 | 13 | 14 | 15 | 16 | 17 | 18 | 19 | 20 | |
| 1 | **0** | -2 | ∞ | 2 | -5 | | | | -2 | 2 | -2 | -4 | | | -3 | -1 | -6 | | 3 | 0 | 1 |
| 2 | 4 | **0** | 1 | 4 | 2 | 9 | | 0 | ∞ | 2 | 0 | -4 | 3 | 5 | -1 | 2 | -6 | 3 | | | 2 |
| 3 | ∞ | -3 | **0** | 1 | -4 | | -6 | | -5 | 1 | -3 | -5 | 2 | 2 | -6 | -6 | -5 | 3 | 0 | -5 | 3 |
| 4 | -5 | 9 | 9 | **0** | -7 | | ∞ | -4 | -8 | 0 | -4 | -6 | 0 | | 8 | -2 | -8 | -6 | -3 | -6 | 4 |
| 5 | 2 | -4 | -1 | 0 | **0** | | -1 | -4 | -1 | 0 | ∞ | -4 | 7 | 1 | -5 | 0 | -8 | 1 | 11 | 1 | 5 |
| 6 | -7 | -12 | -4 | -8 | -11 | **0** | | -12 | -12 | ∞ | -5 | -7 | -8 | -7 | -13 | -10 | -7 | -10 | 2 | 8 | 6 |
| 7 | -9 | -8 | -5 | ∞ | -11 | | **0** | -7 | -12 | -4 | -8 | -10 | -4 | | -11 | -11 | -12 | -10 | -7 | -10 | 7 |
| 8 | -8 | -7 | -4 | | -10 | | 1 | **0** | -9 | -3 | -5 | -9 | -3 | 5 | -8 | ∞ | -11 | -9 | -6 | -9 | 8 |
| 9 | -3 | ∞ | -2 | | -11 | 6 | 0 | 5 | **0** | -4 | 4 | -10 | -3 | 7 | 1 | 1 | -12 | -2 | 5 | 2 | 9 |
| 10 | -11 | -7 | -10 | -3 | -7 | ∞ | -8 | -7 | -8 | **0** | -4 | -6 | 0 | -2 | -8 | -7 | -8 | -6 | -8 | -6 | 10 |
| 11 | | | | | ∞ | | | 1 | 1 | 5 | **0** | -1 | | | 0 | | -3 | | 6 | 3 | 11 |
| 12 | -7 | 2 | -5 | -9 | 0 | -1 | | 2 | -1 | -8 | 1 | **0** | -9 | 7 | 1 | -4 | ∞ | -7 | 2 | 1 | 12 |
| 13 | 2 | -3 | 10 | 4 | -6 | | | | -7 | 1 | 0 | -5 | **0** | 5 | -6 | 2 | -7 | ∞ | -2 | -5 | 13 |
| 14 | | 4 | 8 | 9 | 2 | | | 5 | 1 | 9 | 4 | 3 | 9 | **0** | 2 | 7 | 1 | 3 | ∞ | 3 | 14 |
| 15 | -3 | 1 | -1 | -5 | 3 | 3 | 3 | -9 | | -5 | 3 | -3 | -5 | 6 | **0** | -7 | 2 | -3 | 6 | ∞ | 15 |
| 16 | -10 | -9 | -6 | | -12 | | -1 | ∞ | -13 | -5 | 0 | -11 | -5 | -6 | -12 | **0** | -13 | -11 | -8 | -11 | 16 |
| 17 | -6 | -11 | -12 | -7 | -9 | 1 | | -11 | -9 | | -11 | ∞ | -8 | -6 | -12 | -9 | **0** | -8 | 0 | -3 | 17 |
| 18 | 1 | 9 | -2 | 1 | -8 | | -9 | -3 | -5 | -1 | 4 | -7 | ∞ | -2 | -8 | -8 | -9 | **0** | 0 | -3 | 18 |
| 19 | | | | | | | | | | | | | | ∞ | | | | | **0** | 6 | 19 |
| 20 | -5 | 0 | -7 | 1 | -13 | 9 | -2 | | -10 | -6 | -2 | -12 | -5 | 11 | ∞ | -1 | -14 | -3 | 3 | **0** | 20 |
| | 1 | 2 | 3 | 4 | 5 | 6 | 7 | 8 | 9 | 10 | 11 | 12 | 13 | 14 | 15 | 16 | 17 | 18 | 19 | 20 | |

j = 1

**(10 1)(1 10) = (10 10): -7**

**CYCLE P = [10 12 17 1 5 2 9 15 20 8 16 4 7 3 1 6 10]: (10 17 5 9 20 16 7 1): -7**

j = 2

**(5 2)(2 5) = (5 5): -2**

**CYCLE P = [5 13 18 6 10 17 12 1 3 7 4 16 8 20 15 9 2 11 5]: (5 18 10 12 3 4 8 15 2): -2**

j = 3

**(4 3)(3 4) = (4 4): 10**

**CYCLE P = [4 16 8 6 10 12 17 11 5 2 9 15 20 14 19 1 3 7 4]: (4 8 10 17 5 9 20 19 3): 10**

(10 3)(3 13) = (10 13): -8; (10 3)(3 18) = (10 18): -7; (17 3)(3 13) = (17 13): -10;

(17 3)(3 18) = (17 18): -9; (20 3)(3 18) = (20 18): -4

j = 5

(9 5)(5 1) = (9 1): -6;

**(9 5)(5 9) = (9 9): -12**

**CYCLE P = [9 15 20 8 16 4 7 6 10 12 17 11 5 2 9]: (9 20 16 7 10 17 5): -12**

(9 5)(5 13) = (9 13): -4; (9 5)(5 18) = (9 18): -10; (20 5)(5 1) = (20 1): -8; (20 5)(5 2) = (20 2): -10;

(20 5)(5 9) = (20 9): -14; (20 5)(5 13) = (20 13): -6; (20 5)(5 18) = (20 18): -12

j = 7

**(3 7)(7 3) = (3 3): -8**

**CYCLE P = [3 18 13 6 10 12 17 11 5 2 9 15 20 8 16 4 7 1 3]: (3 13 10 17 5 9 20 16 7): -8**

(18 7)(7 1) = (18 1): -12; (18 7)(7 3) = (18 3): -11;

j = 9

**(20 9)(9 20) = (20 20): -12**

**CYCLE P = [20 8 16 4 7 6 10 12 17 11 5 2 9 15 20]: (20 16 7 10 17 5 9): -12**

j = 13

**(10 13)(13 10) = (10 10): -7**

**CYCLE P = [10 12 17 11 5 2 9 15 20 8 16 4 7 1 3 18 13 6 10]: (10 17 5 9 20 16 7 3 13): -7**

j = 18

(9 18)(18 1) = (9 1): -9;

**(10 18)(18 10) = (10 10): -8**

**CYCLE P = [10 <u>12</u> 17 <u>11</u> 5 <u>2</u> 9 <u>15</u> 20 <u>8</u> 16 <u>4</u> 7 <u>1</u> 3 <u>13</u> 18 <u>6</u> 10]: (10 17 5 9 20 16 7 3 18): -8**

**(17 18)(18 17) = (17 17): 6**

**CYCLE P = [17 <u>10</u> 6 <u>7</u> 4 <u>8</u> 16 <u>4</u> 7 <u>1</u> 3 <u>13</u> 18 <u>12</u> 17]: (17 6 4 16 7 3 18): 6**

(20 18)(18 1) = (20 1): <u>-11</u>

**ACCEPTABLE CYCLES**

**1. DC 9; -8; (10 17 5 9 20 16 7 3 18). 2. DC 9; -8; (3 13 10 17 5 9 20 16 7).**

**3. DC 9; -2; (5 18 10 12 3 4 8 15 2). 4. DC 9; 10; (4 8 10 17 5 9 20 19 3).**

**5. DC 8, 9; -7; (10 17 5 9 20 16 7 1). 6. DC 1, 5, 8; -12; (9 20 16 7 10 17 5).**

**2-CIRCUIT UNLINKED PATHS**

**1. DC 9; -4; [*5* 18 10 12 3 4 8 15 2 *11*]. 2. DC 1, 8, 9; 2; [*20* 16 7 10 17 5 9 *15*].**

**3. DC 1, 2, 9; 2; [*20* 16 7 10 17 5 18 *15*]. 4. DC 1, 2, 8, 9; 5; [*20* 16 7 10 17 5 18 *15*].**

**5. DC 2, 3 5, 6, 9; 1; [*3* 13 10 17 5 *1*].**

$P_{80}$

| | 1 | 2 | 3 | 4 | 5 | 6 | 7 | 8 | 9 | 10 | 11 | 12 | 13 | 14 | 15 | 16 | 17 | 18 | 19 | 20 | |
|---|---|---|---|---|---|---|---|---|---|---|---|---|---|---|---|---|---|---|---|---|---|
| 1 | | 15 | 12 | 1 | 17 | | | 4 | 17 | 8 | 2 | 10 | | | 8 | 20 | 10 | | 20 | 9 | 1 |
| 2 | 18 | | 12 | 12 | 2 | 17 | | 4 | | 18 | | 10 | 5 | 15 | 8 | 4 | 10 | 5 | 1 | | 2 |
| 3 | | 15 | | 3 | 17 | | 16 | 4 | 5 | 13 | 2 | 10 | 3 | 15 | 8 | 20 | 10 | | 20 | 9 | 3 |
| 4 | 18 | 15 | 19 | | 17 | | 8 | 4 | 5 | 8 | 9 | 10 | 5 | | 18 | | 10 | 5 | 20 | 9 | 4 |
| 5 | 18 | 15 | 12 | 3 | | | 16 | 4 | 17 | 18 | | 10 | | 15 | 8 | 4 | 10 | 5 | 1 | 9 | 5 |
| 6 | 18 | 5 | 6 | 6 | 16 | | 8 | 4 | 5 | 8 | 13 | 11 | 6 | 15 | 8 | 4 | 1 | 5 | 1 | 19 | 6 |
| 7 | 18 | 5 | 12 | | 17 | | | 11 | 5 | 7 | 9 | 10 | 5 | 15 | 9 | 20 | 10 | 5 | 20 | 9 | 7 |
| 8 | 18 | 5 | 12 | | 17 | | 8 | | 5 | 7 | 9 | 10 | 5 | 15 | 9 | | 10 | 5 | 20 | 9 | 8 |
| 9 | 18 | 15 | 12 | | 17 | 17 | 16 | 11 | | 7 | 18 | 10 | 5 | 15 | 9 | 20 | 10 | 3 | 20 | 9 | 9 |
| 10 | 18 | 15 | 7 | 1 | 17 | | 16 | 4 | 5 | | 3 | 10 | 3 | 15 | 8 | 20 | 10 | 3 | 20 | 9 | 10 |
| 11 | | 15 | | 12 | | | 8 | 11 | 12 | 8 | | 10 | | | 8 | | 10 | 5 | 20 | 9 | 11 |
| 12 | 13 | 15 | 6 | 6 | 12 | 12 | | 4 | 5 | 13 | 3 | | 6 | 15 | 8 | 4 | | 6 | 1 | 9 | 12 |
| 13 | 13 | 5 | 19 | 1 | 17 | | | 4 | 5 | 13 | 2 | | | 15 | 9 | 4 | 10 | | 1 | 9 | 13 |
| 14 | | 15 | 14 | 3 | 17 | | | 4 | 5 | 8 | 2 | 10 | 5 | | 14 | 4 | 10 | 5 | | 9 | 14 |
| 15 | 13 | 15 | 7 | 6 | 2 | 17 | 8 | 4 | 17 | 8 | 2 | 11 | 6 | | | 4 | 11 | 6 | | 11 | 15 |
| 16 | 18 | 5 | 12 | | 17 | | 16 | | 5 | 7 | 3 | 10 | 5 | 15 | 8 | | 10 | 5 | 20 | 9 | 16 |
| 17 | 7 | 15 | 7 | 6 | 2 | 17 | 16 | 4 | 11 | | 2 | 2 | 3 | 15 | 8 | 4 | | 3 | 1 | 9 | 17 |
| 18 | 7 | 8 | 7 | 3 | 17 | | 16 | 4 | 5 | 18 | 3 | 10 | | 15 | 9 | 20 | 10 | | 1 | 9 | 18 |
| 19 | | | | | | | | | | | | | | | | | | | | | 19 |
| 20 | 18 | 5 | 12 | 1 | 17 | 17 | 16 | 11 | 5 | 7 | 13 | 10 | 5 | 13 | | 20 | 10 | 5 | | | 20 |
| | 1 | 2 | 3 | 4 | 5 | 6 | 7 | 8 | 9 | 10 | 11 | 12 | 13 | 14 | 15 | 16 | 17 | 18 | 19 | 20 | |

$\sigma_1^{-1} M^{-}(80)$

| | 3 | 9 | 1 | 7 | 11 | 10 | 4 | 16 | 2 | 6 | 5 | 17 | 18 | 19 | 20 | 8 | 12 | 13 | 14 | 15 | |
|---|---|---|---|---|---|---|---|---|---|---|---|---|---|---|---|---|---|---|---|---|---|
| | 1 | 2 | 3 | 4 | 5 | 6 | 7 | 8 | 9 | 10 | 11 | 12 | 13 | 14 | 15 | 16 | 17 | 18 | 19 | 20 | |
| 1 | **0** | -2 | ∞ | 2 | -5 | | | | -2 | 2 | -2 | -4 | | | -3 | -1 | -6 | | 3 | 0 | 1 |
| 2 | 4 | **0** | 1 | 4 | 2 | 9 | | 0 | ∞ | 2 | 0 | -4 | 3 | 5 | -1 | 2 | -6 | 3 | | | 2 |
| 3 | ∞ | -3 | **0** | 1 | -4 | | -6 | | -5 | 1 | -3 | -5 | 2 | 2 | -6 | -6 | -5 | 3 | 0 | -5 | 3 |
| 4 | -5 | 9 | 9 | **0** | -7 | | ∞ | -4 | -8 | 0 | -4 | -6 | 0 | | 8 | -2 | -8 | -6 | -3 | -6 | 4 |
| 5 | 2 | -4 | -1 | 0 | **0** | | -1 | -4 | -1 | 0 | ∞ | -4 | 7 | 1 | -5 | -2 | -8 | 1 | 11 | 1 | 5 |
| 6 | -7 | -12 | -4 | -8 | -11 | **0** | | -12 | -12 | ∞ | -5 | -7 | -8 | -7 | -13 | -10 | -7 | -10 | 2 | 8 | 6 |
| 7 | -9 | -8 | -5 | ∞ | -11 | | **0** | -7 | -12 | -4 | -8 | -10 | -4 | | -11 | -11 | -12 | -10 | -7 | -10 | 7 |
| 8 | -8 | -7 | -4 | | -10 | | 1 | **0** | -9 | -3 | -5 | -9 | -3 | 5 | -8 | ∞ | -11 | -9 | -6 | -9 | 8 |
| 9 | -9 | ∞ | -2 | | -11 | 6 | 0 | 5 | **0** | -4 | -3 | -10 | -4 | 7 | 1 | 1 | -12 | -10 | 5 | 2 | 9 |
| 10 | -11 | -7 | -10 | -3 | -7 | ∞ | -8 | -7 | -8 | **0** | -4 | -6 | -8 | -2 | -8 | -7 | -8 | -7 | -8 | -6 | 10 |
| 11 | | | | | ∞ | | | 1 | 1 | 5 | **0** | -1 | | | 0 | | -3 | | 6 | 3 | 11 |
| 12 | -7 | 2 | -5 | -9 | 0 | -1 | | 2 | -1 | -8 | 1 | **0** | -9 | 7 | 1 | -4 | ∞ | -7 | 2 | 1 | 12 |
| 13 | 2 | -3 | 10 | 4 | -6 | | | | -7 | 1 | 0 | -5 | **0** | 5 | -6 | 2 | -7 | ∞ | -2 | -5 | 13 |
| 14 | | 4 | 8 | 9 | 2 | | | 5 | 1 | 9 | 4 | 3 | 9 | **0** | 2 | 7 | 1 | 3 | ∞ | 3 | 14 |
| 15 | -3 | 1 | -1 | -5 | 3 | 3 | -8 | -9 | | -5 | 3 | -3 | -5 | 6 | **0** | -7 | 2 | -3 | 6 | ∞ | 15 |
| 16 | -10 | -9 | -6 | | -12 | | -1 | ∞ | -13 | -5 | 0 | -11 | -5 | -6 | -12 | **0** | -13 | -11 | -8 | -11 | 16 |
| 17 | -6 | -11 | -12 | -7 | -9 | 1 | | -11 | -9 | | -11 | ∞ | -10 | -6 | -12 | -9 | **0** | -9 | 0 | -3 | 17 |
| 18 | -12 | 9 | -11 | 1 | -8 | | -9 | -3 | -5 | -1 | 4 | -7 | ∞ | -2 | -8 | -8 | -9 | **0** | 0 | -3 | 18 |
| 19 | | | | | | | | | | | | | | ∞ | | | | | **0** | 6 | 19 |
| 20 | -11 | -10 | -7 | -3 | -13 | 9 | -2 | | -14 | -6 | -2 | -12 | -6 | 11 | ∞ | -1 | -14 | -12 | 3 | **0** | 20 |
| | 1 | 2 | 3 | 4 | 5 | 6 | 7 | 8 | 9 | 10 | 11 | 12 | 13 | 14 | 15 | 16 | 17 | 18 | 19 | 20 | |

We obtain precisely two acceptable cycles by extending two paths:

**(18 3)(3 18) = (18 18): -8**

**CYCLE P = [18 6 10 12 17 11 5 2 9 15 20 8 16 4 7 1 3 13 18]: (18 10 17 5 9 20 16 7 3): -8**

**(9 11)(11 9) = (9 9): -2**

**CYCLE P = [9 15 20 8 16 4 7 6 10 17 12 1 3 13 18 5 11 2 9]: (9 20 16 7 10 12 3 18): -2**

**ACCEPTABLE CYCLES**

**1. DC 9; -8; (18 10 17 5 9 20 16 7 3). 2. DC 9; -2; (9 20 16 7 10 12 3 18).**

Before proceeding further we combine all sets of acceptable cycles obtained and use our formula to limit the number of searches necessary to obtain a tour of smaller value than $|T_1|$.

**1. C 2, 4; 5; (11 9): 5. 2. C 9, 10; 9; (14 15). 3. C 2, 4, 10; 9; (5 9 15). 4. C 2, 4, 7; 9; (9 11 17).**

**5. C 2, 4, 10; 9; (20 2 11). 6. C 3, 5, 6, 7; 2; (4 8 10 12). 7. C 1, 3, 5, 6, 7; 0; (3 4 8 10 12).**

**8. DC 2, 8, 9, 10; 1; (5 1 4 8 10 12).**

**1. DC 8; 8; (14 3 4 8 10 17 5 9 15). 2. DC 8, 9; 2; (12 6 3 4 8 15 2 11).**

**3. DC 1, 8, 9; -12; (16 7 10 17 5 9 20). 4. DC 2, 9, 10; 1; (13 1 4 8 10 17 5).**

**5. DC 5, 7, 9; 2; (2 5 18 1 4 8 15). 6. DC 1, 8, 9, 10; -4; (8 7 10 17 9 11).**

**7. DC 2, 5, 7, 9; -4; (15 5 18 1 4 8). 8. DC 2, 8, 9, 10; 0; (1 4 8 10 17 5).**

**9. DC 1, 5, 6, 8, 9; 5; (7 10 17 5 9 20). 10. C 2, 4, 5, 6, 7; 10; (17 6 16 5 2).**

**11. C 2, 3, 6, 7, 10; 7; (12 4 8 15 2); 12. C 4, 5, 7, 8; -7; (18 10 17 5). 13. C 1, 5, 7, 8; 0; (17 16 13 1).**

**14. C 2, 4, 6, 10; 1; (5 9 20 16). 15. C 4, 5, 7, 8; 2; (12 5 13 10). 16. C 3, 5, 7, 6; (17 6 16 7).**

**17. C 3, 5, 7; 6; (17 6 4).**

**1. DC 9; 0; (2 5 18 10 12 3 4 8 15). 2. DC 1,9; 7; (17 6 4 8 15 2 5 18).**

**3. DC 1, 8, 9; -12; (10 17 5 9 20 16 7). 4. DC 1, 8, 9; -10; (17 6 4 8 15).**

**5. DC 2, 8, 9; -3; (4 8 10 17 5 15 1). 6. DC 2, 9, 10; -3; (10 17 5 18 1 4 8).**

**7. DC 2, 9, 10; 0; (6 4 16 5 18 1 17). 8. DC 5, 7, 9; 2; (15 2 5 18 1 4 8).**

**9. DC 1, 4, 8, 9; -8; (9 20 16 7 10 17). 10. DC 1, 4, 8, 9; (9 20 16 7 10 12).**

**11. C 1, 4, 5, 7, 8; 4; (10 17 5 18 1). 12. C 2, 4, 6, 10; 1; (9 20 16 5).**

**1. DC 9; -8; (10 17 5 9 20 16 7 3 18). 2. DC 9; -8; (3 13 10 17 5 9 20 16 7).**

**3. DC 9; -2; (5 18 10 12 3 4 8 15 2). 4. DC 9; 10; (4 8 10 17 5 9 20 19 3).**

**5. DC 8, 9; -7; (10 17 5 9 20 16 7 1). 6. DC 1, 5, 8; -12; (9 20 16 7 10 17 5).**

**1. DC 9; -8; (18 10 17 5 9 20 16 7 3). 2. DC 9; -2; (9 20 16 7 10 12 3 18).**

Combining them according to the number of points in a cycle as well as its value, we obtain:

**IX.**

**(1) DC 9; -8; (10 17 5 9 20 16 7 3 18). (2) DC 9; -7; (3 13 10 17 5 9 20 16 7).**

**(3) DC 9; -2; (5 18 10 12 3 4 8 15 2). (4) DC 9; 0; (2 5 18 10 12 3 4 8 15).**

**(5) DC 8; 8; (14 3 4 8 10 17 5 9 15). (6) DC 9; 10; (4 8 10 17 5 9 20 19).**

**-----------------------------------------------------------------------------------------------------------------**

**VIII**

**(1) DC 8, 9; -7; (10 17 5 9 20 16 7 1). (2) DC 4, 9; -2; (9 20 16 7 10 12 3 18).**

**(3) DC 8,9; 2; (12 6 3 4 8 15 2 11). (4) DC 8, 9; 7; (17 6 4 8 15 2 5 18).**

**------------------------------------------------------------------------------------------------------------**

**VII**

**(1) DC 1, 8, 9; -12; (9 20 16 7 10 17 5). (2) DC 2, 8, 9; -3; (4 8 10 17 5 15 1).**

**(3) DC 2, 9, 10; 0; (6 4 16 5 18 1 17). (4) DC 2, 9, 10; 1; (13 1 4 8 10 17 5).**

**(5) DC 5, 7, 9; 2; (2 5 18 1 4 8 15).**

**--------------------------------------------------------------------------------------------------------------**

**VI**

**(1) DC 1, 4, 8, 9; -8; (9 20 16 7 10 17). (2) DC 1, 8, 9, 10; -4; (8 7 10 17 9 11).**

**(3) DC 1, 4, 8, 9; -3; (9 20 16 7 10 12). (4) DC 2, 8, 9, 10; 0; (1 4 8 10 17 5).**

**(5) DC 2, 8, 9, 10; 1; (5 1 4 8 10 12). (6) DC 2, 5, 7, 9; 8; (15 5 18 1 4 8).**

**--------------------------------------------------------------------------------------------------------------**

**V**

**(1) C 1, ,3 , 5, 6, 7; 0; (3 4 8 10 12). (2) C 1, 4, 5, 7, 8; 4; (10 17 5 18 1).**

**(3) C 1, 5, 6, 8, 9; 4; (10 17 5 18 1). (4) C 2, 3, 6, 7, 10; 7; (12 4 8 15 2).**

**(5) C 2, 4, 5, 6, 7; 10; (17 6 16 5 2).**

**-----------------------------------------------------------------------------------------------------**

**IV**

**(1) C 4, 5, 7, 8; -7; (18 10 17 5). (2) C 1, 5, 7, 8; 0; (17 16 13 1). (3) C 2, 4, 6, 10; 1; (9 20 16 5).**

**(4) C 3, 5, 6, 7; 2; (4 8 10 12). (5) C 4, 5, 7, 8; 2; (12 5 13 10. (6) C 3, 5, 6, 7; 2; (17 6 16 7).**

---

**III**

**(1) C 3, 5, 7; 6; (17 6 4). (2) C 2, 4, 10; 9; (5 9 15). (3) C 2, 4, 7; 9; (9 11 17).**

---

**II**

**(1) C 2, 4,; 5; (11 9). (2) C 9, 10; 9; (14 15).**

---

We now deal with obtaining tours. $|T_1| = 68$. We thus wish to obtain tours of value no greater than $|T_1|$ - 1 = 67. Since we are using $\sigma_1^{-1}M^-$ to obtain acceptable and 2-circuit cycles, we've obtained cycles of value no greater than $|T_1| - |\sigma_1| - 1 = 11$. Furthermore, the sum of the value of *all* of the cycles in a tour can be no greater than 11.

We first obtain all tours when a point from precisely one 2-cycle has not been included in the cycle, i.e., **IX**.

**IX (1)** Earlier, we obtained a formula for the number of points necessary to obtain a tour from acceptable and 2-circuit cycles. Let $t$ be the number of 2-circuit cycles and $a$ the number of acceptable cycles. Then the number of points needed to obtain a tour is precisely $\frac{n}{2} + 3t + a - 1$. Since we are currently dealing with only acceptable cycles, the formula simplifies to $\frac{n}{2} + a - 1$. $\frac{n}{2} = 10$, while $a - 1 = 1$. Thus, we require 11 points. Therefore, since in all but one case in **IX,** we are missing a point from the 2-cycle (14 19), we patch the circuit obtained from (14 19) to each of **IX** (1),(2),(3),(4). We then obtain the following tours:

(1) (10 12 17 11 5 2 9 15 19 14 20 8 16 4 7 1 3 13 18 6): 57;

(2) (3 18 13 6 10 12 17 11 5 2 9 15 19 14 20 8 16 4 7 1): 58;

(3) (5 13 18 6 10 17 12 1 3 7 4 16 8 20 14 19 15 9 2 11): 63;

(4) (2 11 5 13 18 6 10 17 12 1 3 7 4 16 8 20 14 19 15 9): 65.

In order to better understand our method of patching circuits, we demonstrate it for **IX** (1).

$$[10\ \underline{12}\ 17\ \underline{11}\ 5\ \underline{2}\ 9\ \underline{15}\ 20\ \underline{8}\ 16\ \underline{4}\ 7\ \underline{1}\ 3\ \underline{13}\ 18\ \underline{6}\ 10]$$

$$\downarrow \quad \nwarrow$$

$$19 \rightarrow 14$$

From this point on, we are interested only in tours of value no greater than 56. Thus, the value in $\sigma_1^{-1}M^-$ of the sum of all cycles in a tour must be no greater than 0.

**VIII**

From our formula, if a tour is constructed from two cycles, then the number of points in the two cycles is 8 + 3 = 11. If it contains three cycles, we have 8 + 2 + 2 = 12. No 3-cycle contains points from both 2-cycles **8** and **9**. We have only two 2-cycles the sum of whose values is 15.

**VII**

Our possibilities here are 7 + 4, 7 + 3 + 2, 7 + 3 + 3.

7 + 4. We don't have any 4-cycle that contains a number from **9.**

7 + 3 + 2. No 3-cycle contains a number from **9.** However, the 2-cycle (14 15) does contain a number from **9**. However, its value is 9. No value of a 3-cycle is less than 6. –12 + 6 + 9 > 0. Thus, 7 + 3 + 2 yields no tour. 7 + 3 + 3 doesn't yield a tour since no 3-cycle contains a number from **9**.

**VI**

Our possibilities are 6 + 5, 6 + 4 + 2, 6 + 3 + 2 + 2, 6 + 2 + 2 + 2 + 2.

Only one 5-cycle contains a number from **9**, namely (6).. However, its value is 8, while no 5-cycle has a negative value. 6 + 4 + 2. No 4-cycle contains a number from **9**. (14 15) has value of 9. The smallest value of 6-cycle is –8. Only one 4-cycle has a negative value: (1) C 4, 5, 7, 8 (18 10 17 9). But numbers from 2-cycles 5 and 7 are both contained in **VI** (1). 6 + 3 + 3 + 2. The sum of the values of any two 3-cycles is at least 15. The smallest value of a 2-cycle is 6. 6 + 2 + 2 + 2 + 2. We have only two 2-cycles.

**V**

Our possibilities are 5 + 5 + 2, 5 + 4 + 3, 5 + 4 + 2 + 2, 5 + 3 + 3 + 2, 5 + 2 + 2 + 2 + 2

5 + 5 + 2. The sum of the values of any two 5-cycles is positive. The sum of each 2-cycle is positive.

5 + 4 + 3. No 4-cycle or 3-cycle contains a number from **9**. 5 + 4 + 2 + 2. The only negative value of a 2-cycle, 4-cycle, or 5-cycle is **IV** (1) whose value is –7. The sum of the values of the two 2-cycles is 15. 5 + 3 + 2 + 2.
**V** (1) has a value of 0. All 3-cycles and 2-cycles have positive values. 5 + 2 + 2 + 2 + 2. There are only two 2-cycles.
Thus, our smallest-valued tour **IX** (1): 57.

## 2-CIRCUIT UNLINKED PATHS

**1. C 5, 8; -7; [*6* 18 *10*]. 2. C 5, 7; -7; [*10* 17 *6*]. 3. C 3, 6; -3; [*4* 8 *7*]. 4. C 1, 3; -1; [*7* 1 *4*].**
**5. C 2, 4; 1; [*2* 5 *9*]. 6. C 2, 10; 2; [*9* 15 *2*]. 7. C 4, 7; 2; [*11* 17 *5*]. 8. C 1, 8; 3; [*3* 13 *1*].**
**9. C 2, 4; 3; [*5* 2 *11*]; 10. C 2, 4; 3; [*5* 9 *11*]. 11. C 3, 5, 6; -8; [*6* 4 8 *10*].**
**12. C 3, 5, 7; -2; [*7* 10 12 *4*]. 13. C 2, 4; 10; [*15* 2 11 *20*]. 14. C 1, 3, 5, 6, 7; 1; [*1* 4 8 10 12 *3*].**

## 2-CIRCUIT UNLINKED PATHS

**1. DC 1, 4, 8, 9; -5; [*12* 6 4 8 15 2 *17*]. 2. C 2, 4, 5, 7, 8; 6; [*2* 5 13 10 17 *9*].**

**3. C 2, 3, 5, 6, 7; 6; [*8* 7 10 17 9 *16*]. 4. C 2, 4, 5, 6, 7; 7; [*17* 6 16 5 2 *12*].**

**5. C 1, 5, 7, 8; -1; [*12* 6 13 1 *17*]. 6. C 1, 3, 5, 7; -4; [*7* 10 12 3 *4*]. 7. C 4, 5, 7 ,8; -3; [*17* 6 13 11 *12*];**

**8. C 4, 5, 7, 8; -1; [*18* 10 17 5 *13*]. 9. C 1, 5, 7, 8; 0; [*18* 10 12 3 *13*]. 10. C 2, 4, 6, 10; 0; [*20* 16 5 9 *15*].**

## 2-CIRCUIT UNLINKED PATHS

**1. DC 1, 8, 9; [*17* 6 4 8 15 2 11 *12*]. 2. DC 2, 6, 9, 10; -7; [*7* 10 17 5 18 1 *4*].**

**3. DC 2, 6, 7, 9, 10; -5; [*6* 4 5 18 1 *10*]. 4. DC 2, 3, 6, 9, 10; 5; [*5* 18 10 12 3 *11*].**

**5. C 1, 4, 5, 7, 8; 5; [*9* 20 16 5 *2*].**

## 2-CIRCUIT UNLINKED PATHS

**1. DC 9; -4; [*5* 18 10 12 3 4 8 15 2 *11*]. 2. DC 1, 8, 9; 2; [*20* 16 7 10 17 5 9 *15*].**

**3. DC 1, 2, 9; 2; [*20* 16 7 10 17 5 18 *15*]. 4. DC 1, 2, 8, 9; 5; [*20* 16 7 10 17 5 18 *15*].**

**5. DC 2, 3 5, 6, 9; 1; [*3* 13 10 17 5 *1*].**

We now present the set of 2-circuit paths of smallest value for each row number heading a path.

**1. C 1, 3, 5, 6, 7; 1; [*1* 4 8 10 12 *3*].**

**2. C 2, 4; 1; [*2* 5 *9*].**

**3. DC 2, 3, 5, 6, 9; 1; [*3* 13 10 17 5 *1*].**

**4. C 3, 6; -3; [*4* 8 *7*].**

**5. DC 9; -4; [*5* 18 10 12 3 4 8 15 2 *11*].**

**6. C 3, 5, 6; -8; [*6* 4 8 *10*].**

**7. DC 2, 6, 9, 10; -7; [*7* 10 17 5 18 1 *4*].**

**8. C 2, 3, 5, 6, 7; 6; [*8* 7 10 17 9 *16*].**

**9. C 2, 10; 2; [*9* 15 *2*].**

**10. C 5, 7; -7; [*10* 17 *6*].**

**11. C 4, 7; 2; [*11* 17 *5*].**

**12. DC 1, 4, 8, 9; -5; [*12* 6 4 8 15 2 *17*].**

**15. C 2, 4; 10; [*15* 2 11 *20*].**

**17. DC 1, 8, 9; -10; [*17* 6 4 8 15 2 11 *12*].**

**18. C 4, 5, 7, 8; -1; [*18* 10 17 5 *13*]**

**20. C 2, 4, 6, 10; 0; [*20* 16 5 9 *15*].**

## UNLINKED 2-CIRCUIT PATHS

(1) If $(a\ b)$ represents an acceptable path and we wish to use the arc $(b\ c)$ extend it to an unlinked 2-circuit path, $c$ must be the other point in the same 2-cycle as $a$.

(2) When backtracking to check if a new entry yields an acceptable or 2-circuit cycle, we must keep the following in mind. If we come to an entry in some $P_i$, where we have replaced an entry from a *previous* iteration with a new one, as we proceed from the last entry towards the new entry, we must check each entry that we reach to see if it belongs to the same 2-cycle of $\sigma_1$ as the new entry. Our reason for doing this is that such an old entry may have belonged to a different path obtained previously.

We now proceed to obtain 2-circuit cycles using the 2-circuit paths we obtained earlier. For easier comprehension, given a fixed column, linked paths are written in italics as they are obtained.

j = 1

(3 1)(1 19) = (3 19): 10

j = 2

(9 2)(2 11) = (9 11): 2; (9 2)(2 12) = (9 12): 7; (9 2)(2 17) = (9 17): 10

j = 3

(1 3)(3 13) = (1 13): 3; (1 3)(3 18) = (1 18): 4

j = 4

(7 4)(4 8) = (7 8): -11; (7 4)(4 15) = (7 15): 5

j = 5

(9 5)(5 13) = (9 13): 11; (9 5)(5 18) = (9 18): 3; (11 5)(5 1) = (11 1): 7;

(11 5)(5 2) = (11 2): 5; (11 5)(5 9) = (11 9): 1; (11 5)(5 13) = (11 13): 9;

(11 5)(5 18) = (11 18): 3

j = 6

(10 6)(6 3) = (10 3): -11; (10 6)(6 4) = (10 4): -15; (10 6)(6 13) = (10 13): -15;

(10 6)(6 16) = (10 16): -10; (10 6)(6 18) = (10 18): -13

j = 7

(4 7)(7 1) = (4 1): -6; (4 7)(7 3) = (4 3): -5; (4 7)(7 10) = (4 10): -7; (4 7)(7 17) = (4 17): 2

j = 8

(7 8)(8 15) = (7 15): 4

j = 9

(2 9)(9 15) = (2 15): 2; (2 9)(9 20) = (2 20): 3; (11 9)(9 15) = (11 15): 2;

(11 9)(9 20) = (11 20): 3

j = 10

(4 10)(10 12) = (4 12): -13; (4 10)(10 17) = (4 17): -15; (6 10)(10 12) = (6 12): -14;

(6 10)(10 17) = (6 17): -16

j = 11

(9 11)(11 8) = (9 8): 3; (9 11)(11 12) = (9 12): 3; (9 11)(11 17) = ( 17): 3

[**1** 4 8 10 12 **3**]          [**5** 18 10 12 3 4 8 15 2 **11**]

*(1 11)(11 17) = (1 17): 8; (5 11)(11 9) = (5 9): -3*

j = 12

(4 12)(12 3) = (4 3): -8; (4 12)(12 5) = (4 5): -13; (4 12)(12 9) = (4 9): -6;

(4 12)(12 13) = (4 13): 2; (9 12)(12 3) = (9 3): 8; (9 12)(12 4) = (9 4) 11;

(9 12)(12 6) = (9 6): 2; (17 12)(12 3) = (17 3): -5; (17 12)(12 13) = (17 13): 5

[**17** 6 4 8 15 2 11 **12**]

*(17 12)(12 5) = (17 5): -10; (17 12)(12 7) = (17 7): 6; (17 12)(12 9) = (17 9): -3*

j = 13

(1 13)(13 11) = (1 11): 6; (4 13)(13 11) = (4 11): 5; (10 13)(13 1) = (10 1): -13;

(10 13)(13 11) = (10 11): -12; (10 13)(13 14) = (10 14): 6; (18 13)(13 1) = (18 1): 1

j = 14

(10 14)(14 15) = (10 15): 9

j = 15

(2 15)(15 14) = (2 14): 8; (7 15)(15 2) = (7 2): 5; (7 15)(15 14) = (7 14): 10;

(10 15)(15 2) = (10 2): 10; (11 15)(15 14) = (11 14): 8; (20 15)(15 14) = (20 14): 6

j = 16

(8 16)(16 5) = (8 5): 7; (10 16)(16 5) = (10 5): -9; (10 16)(16 7) = (10 7): -11

j = 17

(4 17)(17 5) = (4 5): <u>-14</u>; (4 17)(17 9) = (4 9): <u>-11</u>; (6 17)(17 5) = (6 5): <u>-15</u>;

(6 17)(17 9) = (6 9): <u>-12</u>; (6 17)(17 11) = (6 11): <u>2</u>; (12 17)(17 5) = (12 5): <u>-4</u>

j = 18

(9 18)(18 1) = (9 1): <u>6</u>; (9 18)(18 10) = (9 10): <u>4</u>; (10 18)(18 8) = (10 8): <u>3</u>;

(10 18)(18 15) = (10 15): <u>1</u>; (11 18)(18 1) = (11 1): <u>4</u>; (11 18)(18 10) = (11 10): <u>2</u>

j = 20

(2 20)(20 16) = (2 16): <u>2</u>; (2 20)(20 19) = (2 19): <u>6</u>; (11 20)(20 19) = (11 19): <u>6</u>;

(15 20)(20 16) = (15 16): <u>9</u>

## 2-CIRCUIT UNLINKED CYCLES

**1. DC 2, 9, 10; -10; (*7* 10 17 5 18 1 *4* 8) 2. DC 2, 4, 9, 10; 5; (*1* 4 8 10 12 *3* 13)**

**3. C 3, 5, 6, 7; -5; (*4* 8 *7* 10 12 5) 4. C 3, 5, 6, 7; -15; (*6* 4 8 *10* 17) 5. C 2, 4, 7, 10; 7; (*9* 15 *2* 11 17)**

**6. C 2, 4, 7, 10; 10; (*9* 15 *2* 11 12) 7. C 5, 7, 8; -14; (*10* 17 *6* 13) 8. C C 2, 4, 10; 3; (*2* 5 *9* 15)**

**9. C 2, 4, 10; 3; (*2* 5 *9* 20) 10. C 2, 4, 10; 3; (*9* 15 *2* 5) 11. C 2, 4, 10; 3; (*9* 15 *2* 11)**

**12. C 2, 4, 7; 5; (*11* 17 *5* 9) 13. C 4, 7, 8; 10; (*11* 17 *5* 18)**

## 2-CIRCUIT LINKED PATHS

**1. DC 9; -3; [5** 18 10 12 3 4 8 15 *2* **11** *9***]. 2. DC 9; -3; [5** 18 10 *12* 3 4 8 15 2 **11** *17***].**

**3. DC 1, 8, 9; -10; [17** 6 4 8 15 2 *11* **12** *5***]. 4. DC 1, 8, 9; -3; [17** 6 4 8 15 *2* 11 **12** *9***].**

**5. DC 1, 8, 9; 6; [17** 6 *4* 8 15 2 11 **12** *7***]. 6. DC 1, 4, 8, 9; -1; [12** 6 4 8 15 *2* **17** *9***].**

**7. DC 2, 8, 9, 10; 8; [1** 4 8 10 *12* **3** 11 *17***]. 8. C 2, 4, 6 10; 1; [20** 16 5 *9* **15** *2***].**

**9. C 2, 4; 5; [2** *5* **9** *11***].**

## 2-CIRCUIT CYCLES OBTAINED FROM LINKED PATHS

**1. DC 9; -2; (5** 18 10 *12* 3 4 8 15 2 **11** *17***). 2. DC 9; 11; [17** 6 *4* 8 15 2 11 **12** *7***].**

$P_{100\ 2\text{–}CIRCUIT\ UNLINKED}$

| | 1 | 2 | 3 | 4 | 5 | 6 | 7 | 8 | 9 | 10 | 11 | 12 | 13 | 14 | 15 | 16 | 17 | 18 | 19 | 20 | |
|---|---|---|---|---|---|---|---|---|---|---|---|---|---|---|---|---|---|---|---|---|---|
| 1 | | | **12** | **1** | | | | **4** | 11 | **8** | 13 | **10** | 3 | | | | | **3** | | | 1 |
| 2 | | | | | **2** | | | | **5** | | | | | 15 | 9 | 20 | | | 20 | 9 | 2 |
| 3 | **5** | | | | **17** | | | | | **13** | | | **3** | | | | **10** | | 1 | | 3 |
| 4 | 7 | | 12 | | 17 | | **8** | **4** | 12 | 7 | 13 | 10 | 12 | | | | 10 | | | | 4 |
| 5 | | **15** | **12** | **3** | | | | **4** | | **18** | **2** | **10** | | | **8** | | | **5** | | | 5 |
| 6 | | | | **6** | 17 | | | **4** | | **8** | 17 | 10 | | | | | 10 | | | | 6 |
| 7 | **18** | 15 | | **1** | **17** | | | 4 | | **7** | | | | 15 | 8 | | **10** | **5** | | | 7 |
| 8 | | | | | | | **8** | | **17** | **7** | | | | | | **9** | **10** | | | | 8 |
| 9 | 18 | **15** | 12 | 12 | 2 | 12 | | 11 | | 18 | 2 | 11 | 5 | | **9** | | 11 | 5 | | | 9 |
| 10 | 13 | 15 | | | | **17** | | 18 | | | 13 | | 6 | 13 | 18 | 6 | **10** | 6 | | | 10 |
| 11 | 5 | | | | **17** | | | | 5 | 18 | | | 5 | 15 | 9 | | **11** | 5 | 20 | 9 | 11 |
| 12 | | **15** | | **6** | 17 | **12** | | **4** | | | | | | | **8** | | **2** | | | | 12 |
| 13 | | | | | | | | | | | | | | | | | | | | | 13 |
| 14 | | | | | | | | | | | | | | | | | | | | | 14 |
| 15 | | **15** | | | | | | | | | **2** | | | | | | | | | **11** | 15 |
| 16 | | | | | | | | | | | | | | | | | | | | | 16 |
| 17 | | **15** | 12 | **6** | | **17** | | **4** | | | **2** | **11** | 12 | | **8** | | | | | | 17 |
| 18 | 13 | | | | **17** | | | | | **18** | | | **5** | | | | **10** | | | | 18 |
| 19 | | | | | | | | | | | | | | | | | | | | | 19 |
| 20 | | | | | **16** | | | | **5** | | | | | 15 | **9** | **20** | | | | | 20 |
| | 1 | 2 | 3 | 4 | 5 | 6 | 7 | 8 | 9 | 10 | 11 | 12 | 13 | 14 | 15 | 16 | 17 | 18 | 19 | 20 | |

$\sigma_1^{-1} M^{-}(100)$ *2 – CIRCUIT UNLINKED*

| | 1 | 2 | 3 | 4 | 5 | 6 | 7 | 8 | 9 | 10 | 11 | 12 | 13 | 14 | 15 | 16 | 17 | 18 | 19 | 20 | |
|---|---|---|---|---|---|---|---|---|---|---|---|---|---|---|---|---|---|---|---|---|---|
| 1 | | | **1** | **2** | | | | **-2** | 8 | **2** | 6 | **-4** | *3* | | | | | *4* | | | 1 |
| 2 | | | | | **2** | | | 6 | **1** | | *5* | | | 8 | *2* | 2 | | | 6 | *3* | 2 |
| 3 | **1** | | | | **-4** | | | | | **3** | | | **2** | | | | **-5** | | *10* | | 3 |
| 4 | -6 | | -8 | | -14 | | **-3** | **-4** | -11 | *-7* | 5 | *-13* | *2* | | | | *-15* | | | | 4 |
| 5 | | **-4** | **-1** | **0** | | | | **-4** | | **0** | **-4** | **-6** | | | **-5** | | | **1** | | | 5 |
| 6 | | | | **-8** | -15 | | | **-12** | -12 | **-8** | 2 | *-14* | | | | | *-16* | | | | 6 |
| 7 | **-9** | 5 | | **-7** | **-11** | | | *-11* | | **-4** | | | | 10 | *4* | | **-12** | **-10** | | | 7 |
| 8 | | | | | 7 | | **1** | | **-7** | **-3** | | | | | | **6** | **-11** | | | | 8 |
| 9 | 6 | **2** | 8 | 11 | *4* | 2 | | 3 | | 4 | *2* | *3* | *11* | | **1** | *3* | *3* | *5* | | | 9 |
| 10 | -13 | 10 | -11 | -15 | -9 | **-7** | -11 | 3 | | | -12 | | *-15* | *6* | 1 | *-10* | **-8** | *-13* | | | 10 |
| 11 | 4 | | | | **2** | | | | 1 | 2 | | | *9* | 8 | *2* | | **1** | *3* | 6 | *3* | 11 |
| 12 | | **-13** | | **-9** | -4 | **-1** | | **-13** | | | | | | | **-14** | | **-5** | | | | 12 |
| 13 | | | | | | | | | | | | | | | | | | | | | 13 |
| 14 | | | | | | | | | | | | | | | | | | | | | 14 |
| 15 | | **1** | | | | | | | | | **1** | | | | | 9 | | | | **10** | 15 |
| 16 | | | | | | | | | | | | | | | | | | | | | 16 |
| 17 | 7 | **-11** | -5 | **-7** | | **1** | | **-11** | | | **-11** | **-10** | **5** | | **-12** | | | | | | 17 |
| 18 | | | | | **-8** | | | | | **-1** | | | **-1** | | | | **-9** | | | | 18 |
| 19 | | | | | | | | | | | | | | | | | | | | | 19 |
| 20 | | | | | **0** | | | | **-1** | | | | | 6 | **0** | **-1** | | | | | 20 |
| | 1 | 2 | 3 | 4 | 5 | 6 | 7 | 8 | 9 | 10 | 11 | 12 | 13 | 14 | 15 | 16 | 17 | 18 | 19 | 20 | |

j = 1

(4 1)(1 19) = (4 19): 3; (9 1)(1 4) = (9 4): 8; (10 1)(1 19) = (10 19): -4; (11 1)(1 4) = (11 4): 6

j = 3

(4 3)(3 13) = (4 13): -4; (4 3)(3 18) = (4 18): -5; (9 3)(3 7) = (9 7): 10; (17 3)(3 13) = (17 13): -1;

(17 3)(3 18) = (17 18): -2

j = 4

(9 4)(4 16) = (9 16): 6; (10 4)(4 2) = (10 2): 6; (10 4)(4 8) = (10 8): -19; (10 4)(4 9) = (10 9): 10;

(10 4)(4 15) = (10 15): -3; (10 4)(4 16) = (10 16): -17

j = 5

(4 5)(5 1) = (4 1): -9; (4 5)(5 2) = (4 2): -11; (4 5)(5 9) = (4 9): -15; (4 5)(5 13) = (4 13): -7;

(6 5)(5 1) = (6 1): -10; (6 5)(5 2) = (6 2): -12; (6 5)(5 9) = (6 9): -16; (6 5)(5 13) = (6 13): 3;

(10 5)(5 2) = (10 2): -6; (10 5)(5 9) = (10 9): -10; (12 5)(5 1) = (12 1): 1; (12 5)(5 13) = (12 13): 3

[**17** 6 7 8 15 2 11 **12** 5]

*(17 5)(5 1) = (17 1): -5; (17 5)(5 13) = (17 13): -3*

j = 6

(9 6)(6 3) = (9 3): -2; (9 6)(6 4) = (9 4): -6; (9 6)(6 13) = (9 13): -6; (9 6)(6 16) = (9 16): -1

j = 7

[**17** 6 7 8 15 2 11 **12** 7]

*(17 7)(7 3) = (17 3): 4; (17 7)(7 10) = (17 10): 2*

j = 8

(2 8)(8 7) = (2 7): 3; (2 8)(8 10) = (2 10): 10; (9 8)(8 7) = (9 7): 4; (10 8)(8 2) = (10 2): -7;

(10 8)(8 15) = (10 15): -20

j = 9

(1 9)(9 15) = (1 15): 9; (1 9)(9 20) = (1 20): 10; (4 9)(9 15) = (4 15): -14; (4 9)(9 20) = (4 20): -13;

(6 9)(9 15) = (6 15): -15; (6 9)(9 20) = (6 20): -14; (10 9)(9 20) = (10 20): -8

[**12** 6 4 8 15 2 **17** 9]

*(12 9)(9 11) = (12 11): 3*

j = 10

(2 10)(10 12) = (2 12): 4; (2 10)(10 17) = (2 17): 2; (9 10)(10 12) = (9 12): -2;

(9 10)(10 17) = (9 17): -4

j = 11

(1 11)(11 9) = (1 9): 7; (10 11)(11 9) = (10 9): -11

j = 12

[**11** 17 **5** 18 10 12] [**10** 17 **6** 13 11 12]

*(11 12)(12 3) = (11 3): 1; (11 12)(12 4) = (11 4): 4; (10 12)(12 3) = (10 3): -6;*

*(10 12)(12 4) = (10 4): -3; (10 12)(12 9) = (10 9): -4*

j = 13

(6 13)(13 11) = (6 11): -5; (9 13)913 1) = (9 1): -4

j = 15

(4 15)(15 14) = (4 14): -8; (6 15)(15 14) = (6 14): -9; (6 15)(15 19) = (6 19): 8;

(10 15)(15 2) = (10 2): -19; (10 15)(15 5) = (10 5): -11; (10 15)(15 14) = (10 14): -14

j = 16

(2 16)(16 7) = (2 7): 1; (15 16)(16 7) = (15 7): 8

j = 17

[**1** 4 8 10 13 **3** 13 11 17]

*(1 17)(17 9) = (1 9): 11*

j = 18

(4 18)(18 11) = (4 11): 2

j = 20

(4 20)(20 13) = (4 13): 1; (6 20)(20 19) = (6 19): -11

## 2-CIRCUIT UNLINKED CYCLES

**1. DC 1, 3, 4, 9; 0; (9** 15 **2** 5 18 10 17**). 2. DC 1, 3, 6, 9; 5; (9** 15 **2** 5 18 10 12**).**

**3. C 3, 5, 6, 7; -15; (10** 17 **6** 4 8**). 4. C 1, 5, 7, 8; -9; (10** 17 **6** 13 1**). 5. C 1, 3, 6; -4; (4** 8 **7** 1**).**

## 2-CIRCUIT LINKED PATHS

**1. DC 9; -3; [17** *6* 4 8 15 2 11 **12** 3 18 *10***]. 2. DC 9; -2; [17** *6* 4 8 15 2 11 **12** 3 13 *10***].**

**3. DC 9; 5; [7** 10 17 *5* 18 1 **4** 8 15 2 *11***]. 4. DC 9; 9; [1** 4 *8* 10 12 **3** 13 11 9 20 *16***].**

**5. DC 9; 10; [7** 10 *17* 5 18 1 **4** 8 15 2 *12***]. 6. DC 1, 8, 9; -15; [4** *8* **10** 17 5 9 20 *16***].**

**7. DC 1, 8, 9; -14; [4** *8* **7** 10 17 5 9 20 *16***]. 8. DC 1, 8, 9; [12** 6 4 8 15 *2* **17** 5 *9***].**

**9. DC 1, 8, 9; 1; [6** *4* 8 **10** 17 5 9 20 *7***]. 10. DC 2, 8, 9, 10; 7; [1** 4 8 10 *12* **3** 11 *17***].**

**11. DC 1, 8, 9, 10; 10; [8** 7 10 17 *9* **16** 5 *2***]. 12. C 4, 5, 7, 8; -11; [10** *17* **6** 13 11 *12***].**

**13. C 4, 5, 6, 7; 9; [10** *17* **6** 16 5 *12***].**

## 2-CIRCUIT CYCLES OBTAINED FROM LINKED PATHS

**1. DC 9; -11; (17** *6* 4 8 15 2 11 **12** 3 18 *10***). 2. DC 9; -10; [17** *6* 4 8 15 2 11 **12** 3 13 *10***].**

$P_{120\ 2-CIRCUIT\ UNLINKED}$

| | 1 | 2 | 3 | 4 | 5 | 6 | 7 | 8 | 9 | 10 | 11 | 12 | 13 | 14 | 15 | 16 | 17 | 18 | 19 | 20 | |
|---|---|---|---|---|---|---|---|---|---|---|---|---|---|---|---|---|---|---|---|---|---|
| 1 | | | **12** | **1** | | | | **4** | 11 | **8** | 13 | **10** | 3 | | 9 | | | **3** | | 9 | 1 |
| 2 | | | 6 | 6 | **2** | 12 | 16 | 11 | **5** | 8 | | 10 | 6 | 15 | 9 | 20 | 10 | 6 | 20 | 9 | 2 |
| 3 | **5** | | | | **17** | | | | | **13** | | | **3** | | | | **10** | | 1 | | 3 |
| 4 | 5 | 5 | 12 | | 17 | | **8** | **4** | 5 | 7 | 18 | 10 | 20 | 15 | 9 | | 10 | 3 | 1 | 9 | 4 |
| 5 | | **15** | **12** | **3** | | | | **4** | | **18** | **2** | **10** | | | **8** | | | **5** | | | 5 |
| 6 | 5 | 5 | | **6** | 17 | | | **4** | 5 | **8** | 13 | 10 | 5 | 15 | 9 | | 10 | | 15 | 9 | 6 |
| 7 | **18** | 15 | | **1** | **17** | | | 4 | | **7** | | | | 15 | 8 | | **10** | **5** | | | 7 |
| 8 | | | | | | | **8** | | **17** | **7** | | | | | | **9** | **10** | | | | 8 |
| 9 | 13 | **15** | 6 | 6 | 2 | 12 | 8 | 11 | | 18 | 2 | 11 | 6 | | **9** | 6 | 10 | 5 | | | 9 |
| 10 | 13 | 15 | 6 | 6 | 15 | **17** | | 4 | 11 | | 13 | | 6 | 15 | 8 | 6 | **10** | 6 | 1 | 9 | 10 |
| 11 | 5 | | | 1 | **17** | | | | 5 | 18 | | | 5 | 15 | 9 | | **11** | 5 | 20 | 9 | 11 |
| 12 | 5 | **15** | | **6** | 17 | **12** | | **4** | | | | | | | **8** | | **2** | | | | 12 |
| 13 | | | | | | | | | | | | | | | | | | | | | 13 |
| 14 | | | | | | | | | | | | | | | | | | | | | 14 |
| 15 | | **15** | | | | | 8 | | | | **2** | | | | | 20 | | | | **11** | 15 |
| 16 | | | | | | | | | | | | | | | | | | | | | 16 |
| 17 | | **15** | 12 | **6** | | **17** | | **4** | | 18 | **2** | **11** | 3 | | **8** | | | 3 | | | 17 |
| 18 | 13 | | | | **17** | | | | | **18** | | | **5** | | | | **10** | | | | 18 |
| 19 | | | | | | | | | | | | | | | | | | | | | 19 |
| 20 | | | | | **16** | | | | **5** | | | | | 15 | **9** | **20** | | | | | 20 |
| | 1 | 2 | 3 | 4 | 5 | 6 | 7 | 8 | 9 | 10 | 11 | 12 | 13 | 14 | 15 | 16 | 17 | 18 | 19 | 20 | |

$\sigma_1^{-1}M^-(120)$ *2 – CIRCUIT UNLINKED*

| | 1 | 2 | 3 | 4 | 5 | 6 | 7 | 8 | 9 | 10 | 11 | 12 | 13 | 14 | 15 | 16 | 17 | 18 | 19 | 20 | |
|---|---|---|---|---|---|---|---|---|---|---|---|---|---|---|---|---|---|---|---|---|---|
| 1 | | | **1** | **2** | | | | **-2** | *7* | **2** | *6* | **-4** | *3* | | *9* | | | *4* | | *10* | 1 |
| 2 | | | | | **2** | | 1 | | **1** | | | | | *8* | *2* | *2* | | *-1* | *6* | *3* | 2 |
| 3 | **1** | | | | **-4** | | | | | **3** | | | **2** | | | | **-5** | | *10* | | 3 |
| 4 | -9 | -13 | *-8* | | *-14* | | **-3** | **-4** | *-15* | *-7* | *2* | *-13* | *1* | -8 | *-14* | | *-15* | *-5* | *3* | *-13* | 4 |
| 5 | | **-4** | **-1** | **0** | | | | **-4** | | **0** | **-4** | **-6** | | | **-5** | | | **1** | | | 5 |
| 6 | -10 | -12 | | **-8** | *-15* | | | **-12** | *-16* | **-8** | *-5* | *-14* | *-8* | *-9* | *-15* | | *-16* | | *8* | *-14* | 6 |
| 7 | **-9** | *5* | | **-7** | **-11** | | | *-11* | | **-4** | | | | *10* | *4* | | **-12** | **-10** | | | 7 |
| 8 | | | | | *7* | | **1** | | **-7** | **-3** | | | | | | **6** | **-11** | | | | 8 |
| 9 | -4 | **2** | -2 | -6 | *4* | *2* | *4* | *3* | | *4* | *2* | *3* | *11* | | **1** | *-1* | *-4* | *5* | | | 9 |
| 10 | *-13* | -19 | *-11* | *-15* | -11 | **-7** | *-11* | *-19* | *-10* | | *-12* | | *-15* | -14 | *-20* | *-17* | **-8** | *-13* | *-4* | *-8* | 10 |
| 11 | *4* | | | *6* | **2** | | | | *1* | *2* | | | *9* | *8* | *2* | | **1** | *3* | *6* | *3* | 11 |
| 12 | 1 | **-13** | | **-9** | *-4* | **-1** | | **-13** | | | | | | | **-14** | | **-5** | | | | 12 |
| 13 | | | | | | | | | | | | | | | | | | | | | 13 |
| 14 | | | | | | | | | | | | | | | | | | | | | 14 |
| 15 | | **1** | | | | | | | | | **1** | | | | | *9* | | | | **10** | 15 |
| 16 | | | | | | | | | | | | | | | | | | | | | 16 |
| 17 | *7* | **-11** | *-5* | **-7** | | **1** | | **-11** | | | **-11** | **-10** | ***-1*** | | **-12** | | | *-2* | | | 17 |
| 18 | | | | | **-8** | | | | | **-1** | | | **-1** | | | | **-9** | | | | 18 |
| 19 | | | | | | | | | | | | | | | | | | | | | 19 |
| 20 | | | | | **0** | | | | **-1** | | | | | *6* | **0** | **-1** | | | | | 20 |
| | 1 | 2 | 3 | 4 | 5 | 6 | 7 | 8 | 9 | 10 | 11 | 12 | 13 | 14 | 15 | 16 | 17 | 18 | 19 | 20 | |

j = 1

**(4 1)(1 4) = (4 4): -7**

**CYCLE P = [4 $\underline{16}$ 8 $\underline{4}$ 7 $\underline{6}$ 10 $\underline{12}$ 17 $\underline{11}$ 5 $\underline{3}$ 1 $\underline{7}$ 4]**

$P_1$ **= (4 $\underline{16}$ 8),** $P_2$ **= (7 $\underline{6}$ 10 $\underline{12}$ 17 $\underline{11}$ 5 $\underline{3}$ 1 $\underline{7}$ 4)**

(4 1)(1 19) = (4 19): 0; (6 1)(1 19) = (6 19): -1; (9 1)(1 19) = (9 19): 5; (12 1)(1 19) = (12 19): 10

j = 2

(4 2)(2 11) = (4 11): -3; (10 2)(2 11) = (10 11): -19

j = 3

(9 3)(3 13) = (9 13): 0; (9 3)(3 18) = (9 18): 1

j = 4

(9 4)(4 8) = (9 8): -10; (9 4)(4 16) = (9 16): -8

j = 5

(10 5)(5 9) = (10 9): -12

j = 7

(2 7)(7 1) = (2 1): $\underline{-2}$; (2 7)(7 3) = (2 3): $\underline{-1}$; (2 7)(7 10) = (2 10): -3; (2 7)(7 17) = (2 17): 6

j = 10

(2 10)(10 12) = (2 12): -9; (2 10)(10 17) = (2 17): -11

j = 11

(10 11)(11 19) = (10 19): -7

j = 12

(2 12)(12 3) = (2 3): $\underline{-4}$; (2 12)(12 13) = (2 13): 6

j = 14

(6 14)(14 3) = (6 3): $\underline{-1}$

$P_{140\ 2\text{–CIRCUIT UNLINKED}}$

| | 1 | 2 | 3 | 4 | 5 | 6 | 7 | 8 | 9 | 10 | 11 | 12 | 13 | 14 | 15 | 16 | 17 | 18 | 19 | 20 | |
|---|---|---|---|---|---|---|---|---|---|---|---|---|---|---|---|---|---|---|---|---|---|
| 1 | | | **12** | **1** | | | | **4** | 11 | **8** | 13 | **10** | 3 | | 9 | | | **3** | | 9 | 1 |
| 2 | 7 | | 12 | 6 | **2** | 12 | 16 | 11 | **5** | 7 | | 10 | 12 | 15 | 9 | 20 | 10 | 6 | 20 | 9 | 2 |
| 3 | **5** | | | | **17** | | | | | **13** | | | **3** | | | | **10** | | 1 | | 3 |
| 4 | 5 | 5 | 12 | | 17 | | **8** | **4** | 5 | 7 | 2 | 10 | 20 | 15 | 9 | | 10 | 3 | 1 | 9 | 4 |
| 5 | | **15** | **12** | **3** | | | | **4** | | **18** | **2** | **10** | | | **8** | | | **5** | | | 5 |
| 6 | 5 | 5 | 14 | **6** | 17 | | | **4** | 5 | **8** | 13 | 10 | 5 | 15 | 9 | | 10 | | 1 | 9 | 6 |
| 7 | **18** | 15 | | **1** | **17** | | | 4 | | **7** | | | | 15 | 8 | | **10** | **5** | | | 7 |
| 8 | | | | | | | **8** | | **17** | **7** | | | | | | **9** | **10** | | | | 8 |
| 9 | 13 | **15** | 6 | 6 | 2 | 12 | 8 | 4 | | 18 | 2 | 11 | 3 | | **9** | 4 | 10 | 3 | 1 | | 9 |
| 10 | 13 | 15 | 6 | 6 | 15 | **17** | | 4 | 5 | | 2 | | 6 | 15 | 8 | 6 | **10** | 6 | 11 | 9 | 10 |
| 11 | 5 | | | 1 | **17** | | | | 5 | 18 | | | 5 | 15 | 9 | | **11** | 5 | 20 | 9 | 11 |
| 12 | 5 | **15** | | **6** | 17 | **12** | | **4** | | | | | | | **8** | | **2** | | 1 | | 12 |
| 13 | | | | | | | | | | | | | | | | | | | | | 13 |
| 14 | | | | | | | | | | | | | | | | | | | | | 14 |
| 15 | | **15** | | | | | 8 | | | | **2** | | | | | 20 | | | | **11** | 15 |
| 16 | | | | | | | | | | | | | | | | | | | | | 16 |
| 17 | | **15** | 12 | **6** | | **17** | | **4** | | 18 | **2** | **11** | 3 | | **8** | | | 3 | | | 17 |
| 18 | 13 | | | | **17** | | | | | **18** | | | **5** | | | | **10** | | | | 18 |
| 19 | | | | | | | | | | | | | | | | | | | | | 19 |
| 20 | | | | | **16** | | | | **5** | | | | | 15 | **9** | **20** | | | | | 20 |
| | 1 | 2 | 3 | 4 | 5 | 6 | 7 | 8 | 9 | 10 | 11 | 12 | 13 | 14 | 15 | 16 | 17 | 18 | 19 | 20 | |

$\sigma_1^{-1}M^{-}(140)\ 2-CIRCUIT\ UNLINKED$

| | 1 | 2 | 3 | 4 | 5 | 6 | 7 | 8 | 9 | 10 | 11 | 12 | 13 | 14 | 15 | 16 | 17 | 18 | 19 | 20 | |
|---|---|---|---|---|---|---|---|---|---|---|---|---|---|---|---|---|---|---|---|---|---|
| 1 | | | **1** | **2** | | | | **-2** | *7* | **2** | *6* | **-4** | *3* | | *9* | | | *4* | | *10* | 1 |
| 2 | -2 | | -4 | | **2** | | *1* | | **1** | *-3* | | *-9* | 6 | *8* | *2* | *2* | *-11* | *-1* | *6* | *3* | 2 |
| 3 | **1** | | | | **-4** | | | | | **3** | | | **2** | | | | **-5** | | *10* | | 3 |
| 4 | *-9* | *-13* | *-8* | | *-14* | | **-3** | **-4** | *-15* | *-7* | *-3* | *-13* | *1* | *-8* | *-14* | | *-15* | *-5* | *0* | *-13* | 4 |
| 5 | | **-4** | **-1** | **0** | | | | **-4** | | **0** | **-4** | **-6** | | | **-5** | | | **1** | | | 5 |
| 6 | *-10* | *-12* | -1 | **-8** | *-15* | | | **-12** | *-16* | **-8** | *-5* | *-14* | *-8* | *-9* | *-15* | | *-16* | | *-1* | *-14* | 6 |
| 7 | **-9** | *5* | | **-7** | **-11** | | | *-11* | | **-4** | | | | *10* | *4* | | **-12** | **-10** | | | 7 |
| 8 | | | | | *7* | | **1** | | **-7** | **-3** | | | | | | **6** | **-11** | | | | 8 |
| 9 | *-4* | **2** | *-2* | *-6* | *4* | *2* | *4* | *-10* | | *4* | *2* | *3* | *3* | | **1** | *-8* | *-4* | *1* | *5* | | 9 |
| 10 | *-13* | *-19* | *-11* | *-15* | *-11* | **-7** | *-11* | *-19* | *-12* | | *-19* | | *-15* | *-14* | *-20* | *-17* | **-8** | *-13* | *-7* | *-8* | 10 |
| 11 | *4* | | | *6* | **2** | | | | *1* | *2* | | | *9* | *8* | *2* | | **1** | *3* | *6* | *3* | 11 |
| 12 | *1* | **-13** | | **-9** | *-4* | **-1** | | **-13** | | | | | | | **-14** | | **-5** | | *10* | | 12 |
| 13 | | | | | | | | | | | | | | | | | | | | | 13 |
| 14 | | | | | | | | | | | | | | | | | | | | | 14 |
| 15 | | **1** | | | | | | | | | **1** | | | | | *9* | | | | **10** | 15 |
| 16 | | | | | | | | | | | | | | | | | | | | | 16 |
| 17 | *7* | **-11** | *-5* | **-7** | | **1** | | **-11** | | | **-11** | **-10** | ***-1*** | | **-12** | | | *-2* | | | 17 |
| 18 | | | | | **-8** | | | | | **-1** | | | **-1** | | | | **-9** | | | | 18 |
| 19 | | | | | | | | | | | | | | | | | | | | | 19 |
| 20 | | | | | **0** | | | | **-1** | | | | | *6* | **0** | **-1** | | | | | 20 |
| | 1 | 2 | 3 | 4 | 5 | 6 | 7 | 8 | 9 | 10 | 11 | 12 | 13 | 14 | 15 | 16 | 17 | 18 | 19 | 20 | |

We have only three paths that can be extended: (2 1): -2; (2 3): -4; (6 3): -1. In order to save space, we examine each of these possibilities using rooted trees, $\sigma_1^{-1}M^{-}(140)\ 2-CIRCUIT\ UNLINKED$ and $\sigma_1^{-1}M^{-}$.

(2  1)(1  10) doesn't yield a path smaller than one obtained earlier.  We now examine (2  3):

**2**
↓
*5*
↓
**9** → *20* → *16* → *7* → *12* → *3*
↓ -2
*13*

**6**
↓
*4*
↓
*8*
↓
**10** → *17* → *5* → *9* → *15* → *14* → *3*
↓ 2
*18*

We thus cannot obtain any further cycles.

**2-CIRCUIT UNLINKED CYCLES**

**1. DC 2, 14, 19; -7; (*4* 8 *7* 10 17 5 1).**

**2-CIRCUIT LINKED PATHS**

**1. DC 9; 11; [9** *15* **2** 11 12 6 3 13 1 19 *20*]. **2. DC 8, 9; 5; [12** *6* 4 8 15 2 **17** 5 1 *10*].

**3. DC 1, 8, 9; -18; [10** *17* **6** 4 8 15 2 11 *12*]. **4. DC 1, 4, 8, 9; -14; [10** *17* **6** 4 8 15 2 *12*].

**2-CIRCUIT LINKED CYCLE OBTAINED FROM LINKED PATH**

**1. DC 8, 9; -1; (12** *6* 4 8 15 2 **17** 5 1 *10*).

Before going on, it must be noted that we haven't fully followed the idea of the modified F-W algorithm. Namely, by choosing the smallest 2-circuit unlinked path in each row to check for unlinked 2-circuit cycles, we arbitrarily chose the smallest-valued unlinked path in each row. This doesn't guarantee that we will obtain the smallest possible unlinked 2-circuit cycle with a particular row value as a determining vertex. In order to do so, we would have had to apply the F-W algorithm to *each* possible unlinked 2-circuit path.  Let $n_i$ be the number of acceptable unlinked 2-circuit paths obtained in row $i$, while $MAX(UNLINKED) = \max\{n_i \mid i = 1, 2, ..., n\}$.  In order to obtain all such possibilities, we would have had to apply the F-W algorithm to $MAX(UNLINKED)$ cost matrices whose entries consist of unlinked 2-circuit paths. Similarly, in order to obtain *all* possibilities for linked 2-circuit cycles, we would have had to apply the algorithm to $MAX(LINKED)$ cost matrices. We now add 2-circuit cycles to the set of acceptable cycles obtained earlier. Unlike acceptable cycles, 2-circuit cycles have the property that certain numbers in 2-cycles do not yield edges that correspond to 2-cycles of $\sigma_1$. Thus, the two numbers in such a 2-cycle belong to *different* circuits. Such numbers are either printed in boldface or in italics. Assume that $C$ is a 2-circuit cycle. Let $p$ be the number of points in $C$, while $lp$ is the number of its linking points. Then, if $C$ is an unlinked cycle, $lp = p-2$; if $C$ is linked, then $lp = p-4$ implying if $C$ is unlinked, $p \geq 4$, linked, $p \geq 5$.

To recapitulate, a *linking point* of a cycle, $C$, is the only point of a 2-cycle of $\sigma_1$ contained in $C$. This distinguishes it from a *non-linking point* of $C$ which is one of the two points of a 2-cycle of $\sigma_1$ contained in $C$.

**DC** represents "doesn't contain", **C**, "contains". Each Roman numeral represents $\frac{n}{2}$ - (the number of linking point 2-cycles in $C$). A number in italics and before **C** or **DC** represents a 2-cycle that contains no linking points. When used in patching cycle **Q** to cycle **Q*,** they must always be 2-cycles that have no linking points in **Q***. The numbers in boldface directly after **C** or **DC** represent 2-cycles each of which contains a linking point. A *linking point 2-cycle* is a 2-cycle containing a linking point. We denote the patching of cycles $C_1$ and $C_2$ by $C_1 \oplus C_2$.

Each cycle in $S$ contains at most one point belonging to a 2-cycle, say $t$, where another cycle of $S$ contains a point in $t$. Furthermore, $t$ must belong to a *linking edge* (2-cycle) of the two cycles.

**Theorem 16** While constructing a tour, assume that we have patched $j+1$ cycles to an initial cycle $C_i$ to obtain a tree, $T_{j+1}$, formed from a set of linked cycles $C_{j+1}$ linked to $C_i$. Let $N_j$ represent the number of 2-cycles of $\sigma_1$ that contain no points in $C_i$, while $N_{j+1}$ is the number of 2-cycles in $\sigma_1$ that contain at least one point in $C_{j+1}$ but none in $C_i$. Then the respective number of remaining acceptable, unlinked and linked cycles needed to obtain a tour satisfies the inequality $a_{j+} + u_{j+} + 2l_{j+} \leq N_{j+1}$. Furthermore, the number of cycles needed to construct a tour is no greater than $a_{j+1} + u_{j+1} + l_{j+1} + N_1 - N_{j+1} - l_{j+1} + 1$ .The addition of 1 signifies the inclusion of $C_i$.

**Proof.** Let $C_1, C_2, \ldots, C_j$ be the sequence of patching of $j$ cycles added one at a time to a tree of cycles $T_j$ whose root is $C_i$. If $N_j$ represents the number of 2-cycles of $\sigma_1$ that contain no points in $C_j$, then $SN_j$ is the set of 2-cycles containing no points in $C_j$. We now note the following: if $C_a$ is an acceptable cycle that can be patched to $C_j$, then it contains one point in $C_j$ and *at least* one point in a 2-cycle of $SN_j$. Next, assume that $C_u$ is an unlinked 2-circuit cycle that can be patched to $C_j$. Then it has two points that are patched to $C_j$ and *at least* one 2-cycle and perhaps one or more points that belong to $SN_j$. Finally, if $C_l$ is a linked 2-circuit cycle, two pairs of its points are 2-cycles that belong to $SN_j$ and possibly more points also do so, while one point belongs to a 2-cycle in $C_j$. More generally, it follows that if $C_{j+1}$ is the patching following $C_j$ after patching one of the cycles described above, the following inequality holds: $a_{j+1} + u_{j+1} + 2l_{j+1} \leq N_{j+1}$. Here $a_{j+1}$, $u_{j+1}$ and $l_{j+1}$ are the respective number of acceptable cycles, unlinked and linked cycles used in the first $j+1$ patchings to form $T_{j+1}$. Denote by $a_{j+}, u_{j+}, l_{j+}$, the respective number of cycles remaining that can be used to construct a tour by patching. Then $a_{j+} + u_{j+} + 2l_{j+} \leq N_1 - N_{j+1}$. We have already patched $a_{j+1} + u_{j+1} + l_{j+1}$ cycles. We have only $N_1 - N_{j+1}$ linking point 2-cycles left that contain no points in $C_{j+1}$. It follows that the maximum number of cycles that can exist in a tour is no greater than $a_{j+1} + u_{j+1} + l_{j+1} + N_1 - N_{j+1} - l_{j+1} + 1$

*Note 1*. If $C_i$ is an unlinked cycle, then it yields *two* circuits that can't be linked. If a tour exists containing it, there must also exist a patching of cycles containing two linking points – one each to a linking point of each circuit. If we can find no satisfactory cycle of this nature, then no tour exists containing $C_i$.

*Note 2*. We always use the value of the smallest-valued tour obtained thus far as an upper bound for values of any new tour obtained.. When we obtain a tour of smaller value, we substitute its value for the previous one.

*Note 3*. If we obtain a patched set of cycles that contains at least one point in each 2-cycle of $\sigma_1$ but whose set of points doesn't satisfy the formula for the number of points in a set of cycles that yields a tour, $p = \frac{n}{2} + 3t + a - 1$, then we have obtained a derangement, i.e., a set of pair-wise disjoint circuits containing all $n$ points. Here $t$ is the subset of 2-circuit cycles, $a$, the subset of acceptable cycles. If the derangement has a smaller value than all tours encountered thus far, perhaps it can be used to obtain a smaller-valued tour as an upper bound.

*Note 4*. When backtracking to check if a new entry yields an acceptable or 2-circuit cycle, we must keep the following in mind. If we come to an entry in some $P_i$, where we have replaced an entry from a *previous* iteration with a new one, as we proceed from the last entry towards the new entry, we must check each entry that we reach to see if it belongs to the same 2-cycle of $\sigma_1$ as the new entry. Our reason for doing this is that such an old entry may have belonged to a different path obtained previously.

**1**. (1 3); **2**: (2 9); **3**: (4 7); **4**: (5 11); **5**: (6 10); **6**: (8 16); **7**: (12 17); **8**: (13 18) ; **9**: (14 19); **10**: (15 20).

**IX.**

**(1)** *5, 7* **DC 9; -11; (17** *6* 4 8 15 2 11 **12** 3 18 *10*). **(2) 5,** *7* **DC 9; -10; (17** *6* 4 8 15 2 11 **12** 3 13 *10*).

**(3) DC 9; -8; (10 17 5 9 20 16 7 3 18). (4) DC 9; -7; (3 13 10 17 5 9 20 16 7).**

**(5) DC 9; -2; (5 18 10 12 3 4 8 15 2). (6)** *4, 7* **DC 9; -2; (5** 18 10 *12* 3 4 8 15 2 **11** *17*).

**(7) DC 9; 0; (2 5 18 10 12 3 4 8 15). (8) DC 8; 8; (14 3 4 8 10 17 5 9 15).**

---

**VIII**

**(1) DC 8, 9; -7; (10 17 5 9 20 16 7 1). (2) DC 4, 9; -2; (9 20 16 7 10 12 3 18).**

**(3) DC 8, 9; -1;** *5, 7* **(12** *6* 4 8 15 2 **17** 5 1 *10***); (4) DC 8,9; 2; (12 6 3 4 8 15 2 11).**

**(5) DC 8, 9; 7; (17 6 4 8 15 2 5 18). (6) DC 1, 8; 10; (4 8 10 17 5 9 20 19).**

---

**VII**

**(1) DC 1, 8, 9; -12; (9 20 16 7 10 17 5). (2)** *3* **DC 2, 9, 10; -10; (***7* **10 17 5 18 1** *4* **8),**

**(3) DC 2, 8, 9; -3; (4 8 10 17 5 15 1). (4) DC 2, 9, 10; 0; (6 4 16 5 18 1 17).**

**(5) DC 2, 9, 10; 1; (13 1 4 8 10 17 5). (6) DC 5, 7, 9; 2; (2 5 18 1 4 8 15).**

**(7)** *3, 7* **DC 1, 8, 9; 11; (17** 6 *4* 8 15 2 11 **12** *7***).**

---

**VI**

**(1) DC 1, 4, 8, 9; -8; (9 20 16 7 10 17). (2)** *3* **DC 2, 8, 9, 10; -7; (***4* **8** *7* **10 17 5 1).**

**(3) DC 1, 8, 9, 10; -4; (8 7 10 17 9 11). (4) DC 1, 4, 8, 9; -3; (9 20 16 7 10 12).**

**(5) DC 2, 8, 9, 10; 0; (1 4 8 10 17 5). (6)** *2* **DC 1, 3, 4, 9; 0; (***9* **15** *2* **5 18 10 17);**

**(7) DC 2, 8, 9, 10; 1; (5 1 4 8 10 12). (8)** *1* **DC 2, 4, 9, 10; 5; (***1* **4 8 10 12** *3* **13).**

**(9)** *2* **DC 1, 3, 6, 9; 5; (***9* **15** *2* **5 18 10 12). (10) DC 2, 5, 7, 9; 8; (15 5 18 1 4 8).**

---

**V**

**(1) C 1, ,3 , 5, 6, 7; 0; (3 4 8 10 12). (2) C 1, 4, 5, 7, 8; 4; (10 17 5 18 1).**

**(3) C 1, 4, 5, 7, 8; 4; (10 17 5 18 1). (4) C 2, 3, 6, 7, 10; 7; (12 4 8 15 2).**

**(5) C 2, 4, 5, 6, 7; 10; (17 6 16 5 2).**

---

**IV**

**(1) C 4, 5, 7, 8; -7; (18 10 17 5). (2)** *3* **C 4, 5, 6, 7; -5; (***4* **8** *7* **10 12 5).**

**(3) C 1, 5, 7, 8; 0; (17 16 13 1). (4) C 2, 4, 6, 10; 1; (9 20 16 5).**

**(5) C 3, 5, 6, 7; 2; (4 8 10 12). (6) C 4, 5, 7, 8; 2; (12 5 13 10).**

**(7) C 3, 5, 6, 7; 2; (17 6 16 7).**

---

**III**

**(1)** *5* **C 3, 6, 7; -15; (***6* **4 8** *10* **17). (3)** *5* **C 1, 7, 8; -9; (***10* **17** *6* **13 1**).

**(4) C 3, 5, 7; 6; (17 6 4). (5) C 2, 4, 10; 9; (5 9 15). (6) C 2, 4, 7; 9; (9 11 17).**

---

**II**

**(1)** *5* **C 7, 8; -14; (***10* **17** *6* **13). (2)** *3* **C 1, 6; -4; (**4 **8** *7* **1). (3)** *2* **C 4, 10; 3; (***2* **5** *9* **15).**

**(4)** *2* **C 4, 10; 3; (***2* **5** *9* **20). (5)** *2* **C 4, 10; 3; (***9* **15** *2* **5). (6 )** *2* **C 4, 10; 3; (***9* **15** *2* **11).**

**(7)** *4* **C 2, 7; 5; (***11* 17 *5* **9). (8) C 2, 4,; 5; (11 9). (9) C 9, 10; 9; (14 15).**

**(10)** *4* **C 7, 8; 10; (***11* **7** *5* **18).**

---

**XI (1)**: *5, 7* **DC 9; -11**; (**17** *6* 4 8 15 2 11 **12** 3 18 *10*).

**II(2)**: **C 9, 10; 9; (14 15). 10**: (15 20).

**XI(1)** is a linked 2-circuit cycle. The two circuits we obtain from it are $P_1$ = (17 10 6 7 4 16 8 20 15 9 2 5 11) and $P_2$ = (12 1 3 13 18 6 10). The edges [10 6] = [6 10] may be deleted from the respective circuits to obtain the circuit (17 10 12 1 3 13 18 6 7 4 16 8 20 15 9 2 5 11). **II(2)** is an acceptable cycle that yields the circuit $Q_1$ = (14 20 15 19). Both $P_1 \oplus P_2$ and $Q_1$ contain the edge [20 15]. Deleting this edge from the two circuits and patching yields the tour $T_{P_1 \oplus P_2 \oplus Q_1}$

(17 10 12 1 3 13 18 6 7 4 16 8 20 14 19 15 9 2 5 11): -2. The perfect matching $\sigma_1$ has a value of 56. Therefore, $|T_{P_1 \oplus P_2 \oplus Q_1}|$ = 54.

---

**IX(2)**: **5,** *7* **DC 9; -10; (17** *6* 4 8 15 2 11 **12** 3 13 *10***).**

**II(2): C 9, 10: (14 15).**

**IX(2)** is a linked 2-circuit cycle. The first circuit obtained is $P_1$**.** The second is $P_{21}$ = (12 1 3 18 13 6 10).

Again deleting [10 6] = [6 10], we obtain

$P_1 \oplus P_{21}$ = (17 10 12 1 3 18 13 6 7 4 16 8 20 15 9 2 5 11). Patching $Q_1$ to $P_1 \oplus P_{21}$, we obtain the tour

$T_{P_1 \oplus P_{21} \oplus Q_1}$ = (17 10 12 1 3 18 13 6 7 4 16 8 20 14 19 15 9 2 5 11): -1.

It follows that $|T_{P_1 \oplus P_{21} \oplus Q_1}|$ = 55.

**IX(6)** is a linked 2-circuit cycle that yields a tour whose value is 63**.**

**VIII(3):** *5, 7* **DC 8, 9; -1; (12** *6* 4 8 15 2 **17** 5 1 *10*)**. VIII(3)** is a linked 2-cycle circuit. We shall frequently be using the formula given in theorem 1.21 as well as that for the number of points in all cycles. Furthermore, every 2-circuit cycle must contain at least four points.

If we consider only cycles other than **VIII(3)**, $a + l + 2u \le N$. In this case, $N$ = 2, $l \ge 0$. Since **II(2)** is the only cycle containing a point in **9**, $a \ge 1$. Thus, $1 \le a + l + 2u \le 2$. It follows that if $u > 0$, $a + l + 2u \ge 3$. Thus, $u$ = 0. We have the following possibilities: 1, 0; 2, 0. $a = 1$: We can't have a 2-circuit cycle since there is only one 2-cycle of $\sigma_1$ without a point in **VIII(3)** ⊕ **II(2): 8.** Thus, $l$ **=** 0. Our only possibility is $a$ = 2. Our only 2-cycle other than **II(2)** is **II(1)**. But **II(1)** doesn't have a point in **8**. Thus, we can't obtain a tour.

---

**VII(2):** *3* **DC 2, 9, 10; -10; (***7* **10 17 5 18 1** *4* **8). VII(2)** is an unlinked 2-circuit. Thus, $p$ = 8, $lp = 6$. Assuming we don't include **VII(2)**, $a + l + 2u \le 3$. Since $a \ge 1$, **II(2)** is the only cycle containing a point in **9**; it also contains a point in **10**. Thus, it can't be directly linked to **VII(2)**. If $l \ge 1$ or $u \ge 1$, there must exist a 2-circuit cycle containing at least five points. This can't occur. Any cycle that linked **VII(2**) to **II(2)** must contain linking points in **10** and **VII(2)** with at most one other point in **2**. Thus, it moves at most three points.. Therefore, it can't be a 2-circuit cycle. Our only possibilities are $2 \le a \le 3$. Since **VII(2)** is a 2-circuit cycle, our formula yields $8 + 2 + x = 10 + 3 + a - 1 \Rightarrow x = 2 + a$. Thus, if $a = 2$, x = 4 implying that a 4-cycle must link **VII(2)** to **II(2)**. **|VII(2)| + {II(2)|** = -1. |**IV(4)|** = -7. It doesn't contain a point in **10**. No other 4-cycle would yield a tour of value less than 54. If $a$ = 3, x = 5. min |**V**| = 0. Thus, **|VII(2)| + {II(2)|** + min |**V**| = -1. Thus, no tour can be obtained from **VII(2).**

---

**VII(7):** *3, 7* **DC 1, 8, 9; 11; (17** 6 *4* 8 15 2 11 **12** *7*)**. VII(7)** is a linked 2-circuit cycle. $p = 9, lp = 5$.

Assuming we don't include **VII(7),** we must include **II(2).** As in **VII(2)**, we require a cycle that has a linking point in **10**. Since $l = u = 0$, it follows that $2 \le a \le 3$. If $a = 2$, we must have a 3-cycle with linking points in **2**, **10** and a linking point of **VII(7).** If $a = 3$, we have two 2-cycles. min |**II**| = min |**III**| = 6. Thus, either

**|VII(7)| + |II(2)| + min** |**III**| = 26 or |**VII(7)|** + |**II(2)|** + 2(min |**II**|) = 32. No tour is possible.

---

**VI(2)**: *3* **DC 2, 8, 9, 10; -7; (***4* **8** *7* **10 17 5 1)**. **VI(2)** is an unlinked 2-circuit cycle with $p = 7, lp = 5$ . **II(2)** must be included. If $l = 1$, a 2-circuit 5-cycle must contain points in **2, 8, 10** as well as a linking point in **VI(2)**. No such cycle exists. An unlinked 4- cycle must satisfy the above conditions. Thus, $l = w = 0$ . It follows that excluding **II(2),** $1 \le a \le 3$ . Using our formula for the number of points in a tour, 7 + 2 + x = 10 + 3 + $a$ - 1 $\Rightarrow x = 3 + a$ . If $a$ = 1, we have a 4-cycle containing points in **2, 8, 10,** no point in **3** and *one* linking point in **VI(2)**. No acceptable 4-cycle satisfies these conditions. $a$ = 2: a 3-cycle and a 2-cycle |**VI(2)**| + |**II(2)**| + min |**II**| + min |**III**| = 14. $a$ = 3: Since x = 6, we must have three 2-cycles. We have only one left. No tour is possible.

---

**VI(6):** *2* **DC 1, 3, 4, 9; 0; (***9* **15** *2* **5 18 10 17).** This is an unlinked cycle. $p = 7, lp = 5$ **II(2)** must belong to any set of cycles that yield a tour. **II(2)** can be linked to **VI(6)**. |**VI(6)**| + |**II(2)**| = 9. Excluding **VI(6)** and **II(2)**, Suppose $l = 1$. There must exist a linked 5-cycle containing one point in each of **1, 3, 4** as well as precisely *one* linking point in **VI(6)**. Not possible. No unlinked 4-cycle satisfies these conditions. $l = w = 0$ . Excluding **II(2)**, $1 \le a \le 3 \Rightarrow x = 3 + a$ . $a = 1$: An acceptable 4-cycle containing 1, 3, 4. None exists. $a = 2$: A 2-cycle and a 3-cycle. The only remaining 2-cycle is **II(1) C 2, 4; 6; (11 9)**. Although it contains a point in **4**, its point in **2** can't be linked to **VI(6)**. In fact, *no* 3-cycle containing a point in **2** can used in patching. $a = 3$: Since x = 6, we would require three 2-cycles. Not possible. No tour exists.

---

**VI(9):** *2* **DC 1, 3, 6, 9; 5; (***9* **15** *2* **5 18 10 12).** Unlinked cycle. $p = 7, lp = 5$ . **II(2)** can be patched to **VI(9)**. Excluding **VI(9)** and **II(2),** $a + l + 2u \le N, N = 3$. As previously, $l = w = 0$ . Thus, $1 \le a \le 3$ . $a = 1$: 9 + x = 10 + 3 + 1 $\Rightarrow$ x = 5. A 5-cycle. |**VI(9)**| + |**II(2)**| = 14**.** Thus, our 5-cycle must have a value no greater than –17. None do. $a = 2$: x = 6. 2 + 4 or 3 + 3. |**VI(9)**| + |**II(2)**| + min |**II**| = 20. No 4- cycle has a value no greater than –23. Next, two 3-cycles. min |**III**| = 6. By the same reasoning, we can't patch two 3-cycles to **VI(9)** and **II(2)**. $a = 3$: x = 7. The only possibility is 2 + 2 + 3. This would imply that we had three 2-cycles. No tour is possible by patching.

---

Before going on,

(1) We assume that the initial cycle , $C_i$, has a value that is no greater than that of the current $T_{smallest}$ . $S$ is the set of linking point 2-cycles of $C_i$ . As a reminder, all of the 2-cycles of an acceptable cycle are linking points, while all but one in unlinked cycles and all but two in linked cycles have the same property.

(2) Our first branches are constructed by linking to $C_i$ .

(3) When the number of points in $C_i$ is comparatively small, we use the formula that yields the largest number of cycles possible: $a + l + 2u \leq N + 1$ to construct a tour. In the current case, $N = 4$. As we proceed along a branch, we subtract from 1 from $N + 1$ if a new cycle linked to a branch is acceptable or linked, 2, if it is unlinked. We no longer continue linking to a branch when we've obtained a non-positive value from our subtractions.

(4) In general, given a cycle that is a candidate for patching, say $C_i$, none of its non-linking point 2-cycles belongs to $S$ or any of the cycles comprising a branch $B_j$, while at most one of its linking point 2-cycles belongs to $B_j$.

(5) Given a branch $B_j$ that can still be extended, check to see that every 2-cycle not in $C_i$ has one of its points contained in at least one of the cycles found in (4). Start with the 2-cycle that has points in the fewest cycles. If a 2-cycle of $\sigma_1$ exists that has no points in any branch in (4), go to (10).

(6) Patch those cycles in (2) that can be linked to $C_i$.

*Note.* Unlinked cycles can yield a derangement consisting of edges. In such a case, we would have to link circuits rather than cycles. If a tour exists, we always would be able to patch trees formed from sets of patched circuits in the forest constructed until we obtain a tour.

(7) As we go along, we use the formula $p = \frac{n}{2} + 3t + a - 1$ to keep track of the number of points in each branch. Here $t$ represents the number of 2-cycles in a branch. If a branch contains at least one point in each 2-cycle but doesn't satisfy this formula, we can obtain a derangement but not a tour from it..

(8) If no branch of the tree contains at least one point in each 2-cycle of $\sigma_1$, go to (10).

(9) If the sum of the values of cycles in each branch is not smaller than $|T_{smallest}|$, we go to (10). Otherwise, we construct a tour from the smallest-valued tour obtained from a branch and rename it $T_{smallest}$.

(10) Continue with the algorithm. Our next step is to place the smallest-valued linked paths in each row of an $n \ X \ n$ path matrix $P_{linked}$. We then follow the same procedure on $M$ that we used when obtaining unlinked cycles.

**V(1)** ***3*** **DC 1, 2, 8, 9, 10; -5.**

(2) **III(i) (i = 3, 4, 5 ,6)** *2* **C 2, 4, 10; 3.**

**III(9) C 2, 4, 10; 9.**

**II(1) C 2, 4, 5; 5.**

**II(2) C 9, 10; 9.**

(3) None of the tours contains a point in **1** or **8**. No tour possible.

**IV(1) 5 DC 1, 2, 4, 8, 9, 10; -15.**

(2) **IV(5) C 2, 4, 6, 10; 1.**

**III(i) (i = 3, 4, 5, 6)** *2* **C 2, 4, 10; 3.**

**II(1) C 2, 4; 5.**

**II(2) C 9, 10; 9.**

(3) None of the cycles in (2) contains a point in **1** or **8**. No tour possible.

**IV(3)** *5* **DC 2, 3, 4, 6, 9, 10; -9.**

(2) **IV(5) C 2, 4, 6, 10; 1.**

**III(i) (i = 3, 4, 5, 6)** *2* **C 2, 4, 10; 3.**

**III(7)** *4* **C 2, 4, 7; 5.**

**III(9) C 2, 4, 10; 9.**

**III(10) C 2, 4, 7; 9.**

(3) We now construct a tree of linked cycles with **IV(3)** as its root.

$$IV(2)\ 5\ DC\ \ 2,3,4,6,9,10$$

$\swarrow$ $\downarrow$ $\searrow$

$III(2)\ 3,5\,DC\ 2,4,9,10$ $\quad$ $III(7)\ 4,5\ DC\ 3,6,9,10$ $\quad$ $III(10)\ 5\ DC\ 3,6,9,10$

$\downarrow$ $\downarrow$ $\downarrow$

$IV(5)\ 2,3,5\ C\ 10$ $\quad$ $III(2)\ 3,4,5\ DC\ 9,10.$ $\quad$ $III(2)\ 3,5\ DC\ 9,10$

$\downarrow$

$II(2)\ C\ 9,10$

We first note that **IV(2)**, **III(2)** and **III(7)** are all unlinked cycles, while IV(5) and **II(2)** are acceptable cycles. Since **IV(2)** has only three linking point 2-cycles, from the formula $a + l + 2u \le N+1 = 3 + 1 = 4$, the maximum number of cycles that a tour can contain is 4. But $a = 2, u = 2.$ Therefore, $a + l + 2u = 6$. Now consider the formula for the number of points that a set of cycles yielding a tour can have: $p = \frac{n}{2} + 3t + a - 1$. Thus,

$p = \frac{n}{2} + 3t + a - 1 = 10 + 6 + 1 = 17$. Since the sum of the points in the cycles **IV(3)**, **III(2)**, **IV(5)**, **II(2)** is

5 + 4 + 4 + 2 = 15, we obviously have not obtained a tour. However, by patching together circuits that can be patched, we can construct a set of disjoint circuits of edges that contains 20 points, i.e., a derangement the value of whose circuits is 53. Thus, by relaxing the constraints imposed by the formulas, we could obtain a minimally-valued derangement. This situation only occurs when $S$ contains at least one unlinked cycle: The two circuits obtained from linked cycle can be patched to form one circuit, while an acceptable cycle always yields one circuit. We now obtain a derangement from the branch $IV(2) \oplus III(2) \oplus IV(5) \oplus II(2)$.

**IV(2)** = (*10* **17** *6* **13** **1**) yields the circuits of edges **A(1)** = (10 12 17), **A(2)** = (6 18 13 3 1).

**III(2)** = (*4* 8 *7* **1**): **B(1)** = (4 16 8), **B(2)** = (7 3 1).

**IV(5)** = (9 20 16 5): **C** = (9 15 20 8 16 11 5 2).

**II(2)** = (14 15): **D** = (14 20 15 19).

**A(2)** ⊕ **B(2)** = (6 18 13 3 7 1), **C** ⊕ **B(1) =** (9 15 20 8 4 16 11 5 **2**)

**C** ⊕ **B(1)** ⊕ **D** = (9 15 19 14 20 8 4 16 11 5 2)

Our derangement of circuits is thus **(10 12 17)(6 18 13 3 7 1)(9 15 19 14 20 8 4 16 11 5 2): 53**

We now place in each row the 2-circuit linked path of smallest value. Our purpose is to obtain linked 2-circuit cycles.

**2-CIRCUIT LINKED PATHS**

**1. DC 2, 8, 9, 10; 8; [1** 4 8 10 *12* **3** 11 *17***]. 2. C 2, 4; 5; [2** 5 **9** 11**].**

**3. DC 9; -3; [5** 18 10 12 3 4 8 15 *2* **11** *9***]. 4. DC 9; -3; [5** 18 10 *12* 3 4 8 15 2 **11** *17***].**

**5. DC 9; 11; [9** *15* **2** 11 12 6 3 13 1 19 *20***]. 6. DC 1, 8, 9; -18; [10** *17* **6** 4 8 15 2 11 *12***].**

**7. DC 1, 4, 8, 9; -14; [10** *17* **6** 4 8 15 2 *12***]. 8. DC 1, 4, 8, 9; -1; [12** 6 4 8 15 *2* **17** *9***].**

**9. DC 8, 9; 5; [12** *6* 4 8 15 2 **17** 5 1 *10***]. 10. DC 1, 8, 9; -10; [17** 6 4 8 15 2 *11* **12** *5***].**

**11. DC 1, 8, 9; -3; [17** 6 4 8 15 *2* 11 **12** *9***]. 12. DC 1, 8, 9; 6; [17** 6 *4* 8 15 2 11 **12** *7***].**

**13. C 2, 4, 6 10; 1; [20** 16 5 *9* **15** *2***].**

**3., 4., 5.** cannot be extended further. The only cycle containing a point in **9** is (14 15) whose value is 9. Thus, we can't obtain a tour whose value is no greater than –3 by respectively patching (14 15) to each of them. We thus can delete them from consideration. We now construct linked paths by placing the linked path of smallest value in each row.

$P_{160}$(*LINKED CIRCUITS*)

| | 1 | 2 | 3 | 4 | 5 | 6 | 7 | 8 | 9 | 10 | 11 | 12 | 13 | 14 | 15 | 16 | 17 | 18 | 19 | 20 | |
|---|---|---|---|---|---|---|---|---|---|---|---|---|---|---|---|---|---|---|---|---|---|
| 1 | | | **12** | **1** | | | | **4** | | **8** | **3** | **10** | | | | | **11** | | | | 1 |
| 2 | | | | | **2** | | | | **5** | | **9** | | | | | | | | | | 2 |
| 3 | | | | | | | | | | | | | | | | | | | | | 3 |
| 4 | | | | | | | | | | | | | | | | | | | | | 4 |
| 5 | | | | | | | | | | | | | | | | | | | | | 5 |
| 6 | | | | | | | | | | | | | | | | | | | | | 6 |
| 7 | | | | | | | | | | | | | | | | | | | | | 7 |
| 8 | | | | | | | | | | | | | | | | | | | | | 8 |
| 9 | | | | | | | | | | | | | | | | | | | | | 9 |
| 10 | | **15** | | **6** | | **17** | | **4** | | | **2** | **11** | | | **8** | | **10** | | | | 10 |
| 11 | | | | | | | | | | | | | | | | | | | | | 11 |
| 12 | | **15** | | **6** | | **12** | | **4** | **17** | | | | | | **8** | | **2** | | | | 12 |
| 13 | | | | | | | | | | | | | | | | | | | | | 13 |
| 14 | | | | | | | | | | | | | | | | | | | | | 14 |
| 15 | | | | | | | | | | | | | | | | | | | | | 15 |
| 16 | | | | | | | | | | | | | | | | | | | | | 16 |
| 17 | | **15** | | **6** | **12** | **17** | | **4** | | | **2** | **11** | | | **8** | | | | | | 17 |
| 18 | | | | | | | | | | | | | | | | | | | | | 18 |
| 19 | | | | | | | | | | | | | | | | | | | | | 19 |
| 20 | | **15** | | | **16** | | | | **5** | | | | | | **9** | **20** | | | | | 20 |
| | 1 | 2 | 3 | 4 | 5 | 6 | 7 | 8 | 9 | 10 | 11 | 12 | 13 | 14 | 15 | 16 | 17 | 18 | 19 | 20 | |

\

$\sigma_1^{-1}M^-(160)$ (LINKED CIRCUITS

| | 3 | 9 | 1 | 7 | 11 | 10 | 4 | 16 | 2 | 6 | 5 | 17 | 18 | 19 | 20 | 8 | 12 | 13 | 14 | 15 | |
|---|---|---|---|---|---|---|---|---|---|---|---|---|---|---|---|---|---|---|---|---|---|
| | 1 | 2 | 3 | 4 | 5 | 6 | 7 | 8 | 9 | 10 | 11 | 12 | 13 | 14 | 15 | 16 | 17 | 18 | 19 | 20 | |
| 1 | **0** | | **X** | **X** | | | | **X** | | **X** | **X** | **X** | | | | | **X** | | | | 1 |
| 2 | | **0** | | | **X** | | | | **X** | | **X** | | | | | | | | | | 2 |
| 3 | ∞ | | **0** | | | | | | | | | | | | | | | | | | 3 |
| 4 | | | | **0** | | | ∞ | | | | | | | | | | | | | | 4 |
| 5 | | | | | **0** | | | | | | ∞ | | | | | | | | | | 5 |
| 6 | | | | | | **0** | | | | ∞ | | | | | | | | | | | 6 |
| 7 | | | | ∞ | | | **0** | | | | | | | | | | | | | | 7 |
| 8 | | | | | | | | **0** | | | | | | | | ∞ | | | | | 8 |
| 9 | | ∞ | | | | | | | **0** | | | | | | | | | | | | 9 |
| 10 | | **X** | | **X** | | **X** | | **X** | | **0** | **X** | **X** | | | **X** | | **X** | | | | 10 |
| 11 | | | | | ∞ | | | | | | **0** | | | | | | | | | | 11 |
| 12 | | **X** | | **X** | | **X** | | **X** | **X** | | | **0** | | | **X** | | **X** | | | | 12 |
| 13 | | | | | | | | | | | | | **0** | | | | | ∞ | | | 13 |
| 14 | | | | | | | | | | | | | | **0** | | | | | ∞ | | 14 |
| 15 | | | | | | | | | | | | | | | **0** | | | | | ∞ | 15 |
| 16 | | | | | | | | ∞ | | | | | | | | **0** | | | | | 16 |
| 17 | | **X** | | **X** | **X** | **X** | | **X** | | | **X** | **X** | | | **X** | | **0** | | | | 17 |
| 18 | | | | | | | | | | | | | ∞ | | | | | **0** | | | 18 |
| 19 | | | | | | | | | | | | | | ∞ | | | | | **0** | | 19 |
| 20 | | **X** | | | **X** | | | | **X** | | | | | | **X** | **X** | | | | **0** | 20 |
| | 1 | 2 | 3 | 4 | 5 | 6 | 7 | 8 | 9 | 10 | 11 | 12 | 13 | 14 | 15 | 16 | 17 | 18 | 19 | 20 | |

j = 2

(20 2)(2 12) = (20 12): 6

j = 5

(17 5)(5 1) = (17 1): <u>-5</u>; (17 5)(5 13) = (17 13): -3; (17 5)(5 18) = -9

j = 9

(12 9)(9 11) = (12 11): 3

j = 11

(2 11)(11 8) = (2 8): 6; (2 11)(11 12) = (2 12): 6; (2 11)(11 17) = (2 17): 6

j = 12

(2 12)(12 3) = (2 3): 11

j = 18

(17 18)(18 3) = (17 3): -6

$P_{180}$ (*LINKED CIRCUITS*)

| | 1 | 2 | 3 | 4 | 5 | 6 | 7 | 8 | 9 | 10 | 11 | 12 | 13 | 14 | 15 | 16 | 17 | 18 | 19 | 20 | |
|---|---|---|---|---|---|---|---|---|---|---|---|---|---|---|---|---|---|---|---|---|---|
| 1 | | | **12** | **1** | | | | **4** | | **8** | **3** | **10** | | | | | **11** | | | | 1 |
| 2 | | | 12 | | **2** | 12 | | 11 | **5** | | **9** | 11 | | | | | 11 | | | | 2 |
| 3 | | | | | | | | | | | | | | | | | | | | | 3 |
| 4 | | | | | | | | | | | | | | | | | | | | | 4 |
| 5 | | | | | | | | | | | | | | | | | | | | | 5 |
| 6 | | | | | | | | | | | | | | | | | | | | | 6 |
| 7 | | | | | | | | | | | | | | | | | | | | | 7 |
| 8 | | | | | | | | | | | | | | | | | | | | | 8 |
| 9 | | | | | | | | | | | | | | | | | | | | | 9 |
| 10 | | **15** | 12 | **6** | | **17** | | **4** | | | **2** | **11** | | | **8** | | **10** | | | | 10 |
| 11 | | | | | | | | | | | | | | | | | | | | | 11 |
| 12 | | **15** | | **6** | | **12** | | **4** | **17** | | 9 | | | | **8** | | **2** | | | | 12 |
| 13 | | | | | | | | | | | | | | | | | | | | | 13 |
| 14 | | | | | | | | | | | | | | | | | | | | | 14 |
| 15 | | | | | | | | | | | | | | | | | | | | | 15 |
| 16 | | | | | | | | | | | | | | | | | | | | | 16 |
| 17 | 5 | **15** | 18 | **6** | **12** | **17** | | **4** | | | **2** | **11** | 5 | | **8** | | | 5 | | | 17 |
| 18 | | | | | | | | | | | | | | | | | | | | | 18 |
| 19 | | | | | | | | | | | | | | | | | | | | | 19 |
| 20 | | **15** | 12 | | **16** | 12 | | | **5** | | | 2 | | | **9** | **20** | | | | | 20 |
| | 1 | 2 | 3 | 4 | 5 | 6 | 7 | 8 | 9 | 10 | 11 | 12 | 13 | 14 | 15 | 16 | 17 | 18 | 19 | 20 | |

$\sigma_1^{-1}M^{-}(180)$ (LINKED CIRCUITS

| | 3 | 9 | 1 | 7 | 11 | 10 | 4 | 16 | 2 | 6 | 5 | 17 | 18 | 19 | 20 | 8 | 12 | 13 | 14 | 15 | |
|---|---|---|---|---|---|---|---|---|---|---|---|---|---|---|---|---|---|---|---|---|---|
| | 1 | 2 | 3 | 4 | 5 | 6 | 7 | 8 | 9 | 10 | 11 | 12 | 13 | 14 | 15 | 16 | 17 | 18 | 19 | 20 | |
| 1 | **0** | | **X** | **X** | | | | **X** | | **X** | **X** | **X** | | | | | **X** | | | | 1 |
| 2 | | **0** | 11 | | **X** | | | 6 | **X** | | **X** | *6* | | **22** | | | 6 | | | | 2 |
| 3 | ∞ | | **0** | | | | | | | | | | | | | | | | | | 3 |
| 4 | | | | **0** | | | ∞ | | | | | | | | | | | | | | 4 |
| 5 | | | | | **0** | | | | | | ∞ | | | | | | | | | | 5 |
| 6 | | | | | | **0** | | | | ∞ | | | | | | | | | | | 6 |
| 7 | | | | ∞ | | | **0** | | | | | | | | | | | | | | 7 |
| 8 | | | | | | | | **0** | | | | | | | | ∞ | | | | | 8 |
| 9 | | ∞ | | | | | | | **0** | | | | | | | | | | | | 9 |
| 10 | | **X** | -13 | **X** | | **X** | | **X** | | **0** | **X** | **X** | | | **X** | | **X** | | | | 10 |
| 11 | | | | | ∞ | | | | | | **0** | | | | | | | | | | 11 |
| 12 | | **X** | | **X** | | **X** | | **X** | **X** | | *3* | **0** | | | **X** | | **X** | | | | 12 |
| 13 | | | | | | | | | | | | | **0** | | | | | ∞ | | | 13 |
| 14 | | | | | | | | | | | | | | **0** | | | | | ∞ | | 14 |
| 15 | | | | | | | | | | | | | | | **0** | | | | | ∞ | 15 |
| 16 | | | | | | | | ∞ | | | | | | | | **0** | | | | | 16 |
| 17 | -5 | **X** | -6 | **X** | **X** | **X** | | **X** | | | **X** | **X** | *-3* | | **X** | | **0** | *-9* | | | 17 |
| 18 | | | | | | | | | | | | | ∞ | | | | | **0** | | | 18 |
| 19 | | | | | | | | | | | | | | ∞ | | | | | **0** | | 19 |
| 20 | | **X** | 11 | | **X** | 5 | | | **X** | | | *6* | | | **X** | **X** | | | | **0** | 20 |
| | 1 | 2 | 3 | 4 | 5 | 6 | 7 | 8 | 9 | 10 | 11 | 12 | 13 | 14 | 15 | 16 | 17 | 18 | 19 | 20 | |

j = 3

(10 3)(3 13) = (10 13): -11; (10 3)(3 18) = (10 18): -10

j = 6

(2 6)(6 3) = (2 3): 1; (2 6)(6 4) = (2 4): -3; (2 6)(6 13) = (2 13): -3; (2 6)(6 16) = (2 16): 2;

(2 6)(6 18) = (2 18): -1; (20 6)(6 3) = (20 3): 1; (20 6)(6 4) = (20 4): -3; (20 6)(6 13) = (20 13): -3;

(20 6)(6 18) = (20 18): -1

j = 8

(2 8)(8 7) = (2 7): 7; (2 8)(8 10) = (2 10): 10; (2 8)(8 15) = (2 15): 5

j = 10

(2 10)(10 12) = (2 12): 4; (2 10)(10 17) = (2 17): 2

j = 13

(2 13)(13 1) = (2 1): -1;

**(10 13)(13 10) = (10 10): -10**

**CYCLE P = [10 12 17 10 6 7 4 16 8 20 15 9 2 5 11 17 12 1 3 18 13 6 10]**

$P_1$ **= (10 12 17),** $P_2$ **= (6 7 4 16 8 20 15 9 2 5 11 17 12 1 3 18 13)**

**(10 17 6 4 8 15 2 11 12 3 13)**

(20 13)(13 1) = (20 1): -1

j = 15

**(2 15)(15 2) = (2 2): 6**

**CYCLE P = [2 11 5 2 9 5 11 16 8 20 15 9 2]**

$P_1$ **= (2 11 5),** $P_2$ **= (9 5 11 16 8 20 15)**

**(2 5 9 11 8 15)**

(2 15)(15 14) = (2 14): 11

j = 16

(2 16)(16 7) = (2 7): 1

j = 18

**(10 18)(18 10) = (10 10): -11**

**CYCLE P = [10 12 17 10 6 7 4 16 8 20 15 9 2 5 11 17 12 1 3 13 18 6 10]**

$P_1$ **= (10 12 17),** $P_2$ **= (6 7 4 16 8 20 15 9 2 5 11 17 12 1 3 13 18)**

**(10 17 6 4 8 15 2 11 12 3 18)**

$P_{200}(LINKED\ CIRCUITS)$

| | 1 | 2 | 3 | 4 | 5 | 6 | 7 | 8 | 9 | 10 | 11 | 12 | 13 | 14 | 15 | 16 | 17 | 18 | 19 | 20 | |
|---|---|---|---|---|---|---|---|---|---|---|---|---|---|---|---|---|---|---|---|---|---|
| 1 | | | **12** | **1** | | | | **4** | | **8** | **3** | **10** | | | | | **11** | | | | 1 |
| 2 | 13 | | 6 | 6 | **2** | 12 | 16 | 11 | **5** | 8 | **9** | 10 | 6 | 15 | 8 | 6 | 10 | 6 | | | 2 |
| 3 | | | | | | | | | | | | | | | | | | | | | 3 |
| 4 | | | | | | | | | | | | | | | | | | | | | 4 |
| 5 | | | | | | | | | | | | | | | | | | | | | 5 |
| 6 | | | | | | | | | | | | | | | | | | | | | 6 |
| 7 | | | | | | | | | | | | | | | | | | | | | 7 |
| 8 | | | | | | | | | | | | | | | | | | | | | 8 |
| 9 | | | | | | | | | | | | | | | | | | | | | 9 |
| 10 | | **15** | 12 | **6** | | **17** | | **4** | | | **2** | **11** | 3 | | **8** | | **10** | 3 | | | 10 |
| 11 | | | | | | | | | | | | | | | | | | | | | 11 |
| 12 | | **15** | | **6** | | **12** | | **4** | **17** | | 9 | | | | **8** | | **2** | | | | 12 |
| 13 | | | | | | | | | | | | | | | | | | | | | 13 |
| 14 | | | | | | | | | | | | | | | | | | | | | 14 |
| 15 | | | | | | | | | | | | | | | | | | | | | 15 |
| 16 | | | | | | | | | | | | | | | | | | | | | 16 |
| 17 | 5 | **15** | 18 | **6** | **12** | **17** | | **4** | | | **2** | **11** | 5 | | **8** | | | 5 | | | 17 |
| 18 | | | | | | | | | | | | | | | | | | | | | 18 |
| 19 | | | | | | | | | | | | | | | | | | | | | 19 |
| 20 | 13 | **15** | 6 | 6 | **16** | 12 | | | **5** | | | 2 | 6 | | **9** | **20** | | 6 | | | 20 |
| | 1 | 2 | 3 | 4 | 5 | 6 | 7 | 8 | 9 | 10 | 11 | 12 | 13 | 14 | 15 | 16 | 17 | 18 | 19 | 20 | |

$\sigma_1^{-1}M^-(200)$ LINKED CIRCUITS

| | 3 | 9 | 1 | 7 | 11 | 10 | 4 | 16 | 2 | 6 | 5 | 17 | 18 | 19 | 20 | 8 | 12 | 13 | 14 | 15 | |
|---|---|---|---|---|---|---|---|---|---|---|---|---|---|---|---|---|---|---|---|---|---|
| | 1 | 2 | 3 | 4 | 5 | 6 | 7 | 8 | 9 | 10 | 11 | 12 | 13 | 14 | 15 | 16 | 17 | 18 | 19 | 20 | |
| 1 | **0** | | **X** | **X** | | | | **X** | | **X** | **X** | **X** | | | | | **X** | | | | 1 |
| 2 | -1 | **0** | 1 | -3 | **X** | 5 | 1 | *6* | **X** | *10* | **X** | *4* | *-3* | 11 | *5* | *2* | *2* | *-1* | | | 2 |
| 3 | ∞ | | **0** | | | | | | | | | | | | | | | | | | 3 |
| 4 | | | | **0** | | | ∞ | | | | | | | | | | | | | | 4 |
| 5 | | | | | **0** | | | | | | ∞ | | | | | | | | | | 5 |
| 6 | | | | | | **0** | | | | ∞ | | | | | | | | | | | 6 |
| 7 | | | | ∞ | | | **0** | | | | | | | | | | | | | | 7 |
| 8 | | | | | | | | **0** | | | | | | | | ∞ | | | | | 8 |
| 9 | | ∞ | | | | | | | **0** | | | | | | | | | | | | 9 |
| 10 | | **X** | *-13* | **X** | | **X** | | **X** | | **0** | **X** | **X** | *-11* | | **X** | | **X** | *-10* | | | 10 |
| 11 | | | | | ∞ | | | | | | **0** | | | | | | | | | | 11 |
| 12 | | **X** | | **X** | | **X** | | **X** | **X** | | *3* | **0** | | | **X** | | **X** | | | | 12 |
| 13 | | | | | | | | | | | | | **0** | | | | | ∞ | | | 13 |
| 14 | | | | | | | | | | | | | | **0** | | | | | ∞ | | 14 |
| 15 | | | | | | | | | | | | | | | **0** | | | | | ∞ | 15 |
| 16 | | | | | | | | ∞ | | | | | | | | **0** | | | | | 16 |
| 17 | *-5* | **X** | *-6* | **X** | **X** | **X** | | **X** | | | **X** | **X** | *-3* | | **X** | | **0** | *-9* | | | 17 |
| 18 | | | | | | | | | | | | | ∞ | | | | | **0** | | | 18 |
| 19 | | | | | | | | | | | | | | ∞ | | | | | **0** | | 19 |
| 20 | -1 | **X** | 1 | -3 | **X** | *5* | | | **X** | | | *6* | *-3* | | **X** | **X** | | *-1* | | **0** | 20 |
| | 1 | 2 | 3 | 4 | 5 | 6 | 7 | 8 | 9 | 10 | 11 | 12 | 13 | 14 | 15 | 16 | 17 | 18 | 19 | 20 | |

We are unable to extend any path further. We have obtained the following linked circuits.

**IX**

**(1)** ***5, 7*** **DC 9; -11; (10 17 6 4 8 15 2 11 12 3 18). (2)** ***5, 7*** **DC 9; -10; (10 17 6 4 8 15 2 11 12 3 13).**

**IV**

**(1)** ***2, 4*** **DC 1, 3, 5, 7, 8, 9 ; 6; (2 5 9 11 15).**

**IX(1)** and **IX(2)** are the same cycles as **IX(1)** and **IX(2)** obtained earlier. We thus will use that table of acceptable, unlinked and linked cycles given earlier together with **IV(1)** to see if we can obtain a tour of value no greater than –3.

In obtaining tours, we consider all tours of value no greater than 56. Since $|\sigma_1| = 56$, we are interested in sets of cycles the total sum of whose values is no greater than 0. We first rank the occurrence of points of 2-cycles in cycles. Those that occur least often are given the ranking i; those that secondly least often, we rank ii, etc. . In general, we attempt to first patch cycles that contain points of lowest ranking. Our reasoning here is that those 2-cycles having no representation in large cycles have this property because they occur less often in cycles than others do. We assume that initial cycles contain the greatest number of points and values no greater than zero. When patching, we indicate the sets of unlinked 2-cycles remaining (boldface), unlinked 2-cycles (italics). As we proceed along a branch, we give the sum of the values of the cycles patched so far at each node of the tree. The 2-cycles that occur in fewest cycles are **9**, **8**, **1, 10, 2** in that order.

**IX(1)** 5, 7 **DC 9; -11**

*l IX(1) 5, 7 DC 9; −11* ⊕ *a II(9) 5, 7, 10; −2* *t = 1, a = 1*. The total number of points we have is 11 + 2 = 13. From our formula, if a tour exists $p = \frac{n}{2} + 3t + a - 1 = 10 + 3 = 13$. We thus have a tour. Furthermore, since its value is –2, this tour becomes a new upper bound.

We new present it as a simple illustration of our method of patching.

$\left[\mathbf{17}\ \underline{10}\ \mathit{6}\ \underline{7}\ 4\ \underline{16}\ 8\ \underline{20}\ 15\ \underline{9}\ 2\ \underline{5}\ 11\ \underline{17}\ \mathbf{12}\ \underline{1}\ 3\ \underline{13}\ 18\ \underline{\mathit{6}}\ 10\ \underline{12}\ \mathbf{17}\right],\ \left[14\ \underline{20}\ 15\ \underline{19}\ 14\right]$

$\left[17\ 10\ \ \mathit{6}\ 7\ 4\ 16\ 8\ 20\ 15\ 9\ 2\ 5\ 11\right],\ \left[12\ 1\ 3\ 13\ 18\ \mathit{6}\ 10\right]$

```
[17  10    6  7  4  16  8  20    15  9  2  5  11]
    ↓   ↑                 ↓    ↑
    12  18                14 → 19
    ↓   ↑
    1   ↑
    ↓   ↑
    3→13
```

We thus obtain $T_{54}$ = [17 10 12 1 3 13 18 6 7 4 16 8 20 14 19 15 9 2 5 11].

*l* **IX(2)** *5, 7* **DC 9; -10**

*l* **IX(2)** *5,7* **DC 9; -10** ⊕ *a* **II(9)** *5, 7, 10;* **-1**

This yields a tour of value 55.

*a* **IX(3) DC 9; -8**

*a* **IX(3) DC 9; -8** ⊕ *a* **II(9)** *10;* **1**

A tour of value 57.

*a* **IX(4) DC 9; -7**

*a* **IX(4) DC 9; -7** ⊕ *a* **II(9)** *10;* **2**

A tour of value 58

*a* **IX(5) DC 9; -2**

*a* **IX(5) DC 9; -2** ⊕ *a* **II(9)** *10*; **7**

A tour of value 63.

*l* **IX(6)** *4, 7* **DC 9**

*l* **IX(6)** *4, 7* **DC 9; -2** ⊕ *a* **II(9)** *10*; **7**

A tour of value 63.

*a* **IX(7) DC 9; 0**

*a* **IX(7) DC 9; 0** ⊕ *a* **II(9)** *10*; **9**

A tour of value 65.

**IX(8) DC 8; 8**

There exists no 2-cycle containing a point in **8.**

*a* **VIII(1) DC 8, 9; -7**

No 3-cycle contains a point in both **8** and **9** and a linking point in some other 2-cycle.

No 2-cycle contains a point in **8**.

*a* **VIII(2) DC 4, 9; -2**

*a* **II(8) C 2, 4; 5** can be linked to **VIII(2)** yielding *a* **VIII(2) DC 4, 9; -2** ⊕ *a* **II(8)** *2* **DC 9; 3.**
We can't patch **II(9)** because the resulting tour would be of value 68, the same as that of $T_1$.

*a* **VIII(3) DC 8, 9; -1**

No 3-cycle contains a point in both **8** and **9** and a linking point in some other 2-cycle.

No 2-cycle contains a point in **8**.

*a* **VII(1) DC 1, 8, 9; -12**

*a* **VII(1) DC 1, 8, 9; -12** ⊕ *a* **II(9)** *10* **DC 1, 8; -3**

No 4-cycle contains a point in each of **1**, **8**, **9**. No 3-cycle contains a point in both (a) **1** and **8,** (b) **1 and 10,** (c) **8** and **10**. No 2-cycle contains a point in **1** or **8**.

*u* **VII(2)** *3* **DC 2, 9, 10; -10**

No 4-cycle contains points in **2**, **9** and **10**.

*a* **III(4) C 2, 4, 10; 9** can be linked to both **VII(2)** and **II(9).**

*u* **VII(2)** *3* **DC 2, 9, 10; -10** ⊕ *a* **III(4)** *3, 4* **DC 9; -1**

*a* **II(9)** can now be patched: *u* **VII(2)** *3* **DC 2, 9, 10; -10** ⊕ *a* **III(4)** *3, 4* **DC 9; -1** ⊕ *a* **II(9)** 3, 4, 10; **8**

**|VII(2)| + |III(4)| + |II(9)| = 64.** $p = 8 + 3 + 2 = 13$. From our formula for a tour,

$p = \frac{n}{2} + 3t + a - 1 = 10 + 3 + 1 = 14$. Thus, we have obtained a derangement whose value is 64.

On the other hand, for each i = 3, 4, 5, 6, **II(i) 2 C 4, 10; 3** can also be patched to both **VII(2)** and **II(9).** In this case we obtain *u* **VII(2)** *3* **DC 2, 9, 10; -10** ⊕ *u* **II(i)** *2, 3, 4* **DC 9; -7** ⊕ *a* II(9) *2, 3, 4, 10*; **2**.

From our formula for the number of points yielding a tour, we see that we have obtained a derangement whose value is 58.

*a* **VII(3) DC 2, 8, 9; -3**

We can patch **II(8)** and **II(9)** to **VII(3)**.

*a* **VII(3) DC 2, 8, 9; -3** ⊕ *a* **II(8)** *4* **DC 8, 9; 2** ⊕ **II(9)** *4, 10* **DC 8; 11**.

No 2-cycle contains a point in **8**.

*a* **VI(1) DC 1, 4, 8, 9; -8**

The only cycles containing a point in **9** are **VIII(6)** and **II(9)**. No 4-cycle contains a point in **1**, **4** and **8**. Only one 4-cycle contains a point in **1: II(2).** However, it also contains a non-linking 2-cycle **3**. No 2-cycle or 3-cycle contains a point in **1**.

*u* **VI(2) 3 DC 2, 8, 9, 10; -7**

We can link **III(4), II(i)**, i = 1, 2, ..., 6 to **VI(2)**. We can then link **II(9)** to each of those patchings. However, we still cannot link a point in 8 to construct a derangement or a tour.

*a* **VI(3) DC 1, 8, 9, 10; -4**

Other than **II(9**), no 2-cycle, 3-cycle or 4-cycle contains a point in **1** or **10** that doesn't include at least two points that are linking points in **VI(3)**. Thus, we can't obtain a derangement or a tour.

*a* **VI(4) DC 1, 4, 8, 9; -3.** A patching cycle can go through at most five points.

We first check to see if we can patch **9** to **VI(4)**.

**II(9)** is the only cycle containing a point of **9** that can be patched to **VI(4)**. We obtain

*a* **VI(4) DC 1, 4, 8, 9; -3** ⊕ *a 10* **II(9) DC 1, 4, 8**; **6.**

We now see if we can patch a cycle containing a point of **8**.

**II(10)** can be patched. However, we keep in mind that the new sum of values is greater than 11.

*a* **VI(4) DC 1, 4, 8, 9; -3** ⊕ *a 10* **II(9) DC 1, 4, 8**; **6** ⊕ *u 4, 7, 10* **II(10) DC 1; 16**

No 2-cycle contains a point in **1**.

No derangement or tour.

*a* **IV(1) DC 1, 2, 3, 6, 9, 10; -7.** We now search for cycles that are non-negative if they contain *more than* four points together with all cycles that contain *fewer* than five points. No cycle with a point in **9** can be patched to **IV(1).** No cycle containing a point in **1** can be patched to **IV(1)**. However, a cycle, **IV(4)**, containing points in **2**, **6** and **10** can be patched to **IV(1)**: *a* **IV(1) DC 1, 2, 3, 6, 9, 10; -7** ⊕ *a 4* **IV(4) DC 1, 3, 9; -6**

**II(9)** containing a point in **9** can now be patched:

*a* **IV(1) DC 1, 2, 3, 6, 9, 10; -7** ⊕ *a 4* **IV(4) DC 1, 3, 9; -6** ⊕ *u 4, 10* **II(2) DC 1, 3; 3**.

Finally, **II(2)** containing points in **1** and **3** can be patched.

*a* **IV(1) DC 1, 2, 3, 6, 9, 10; -7**⊕ *a 4* **IV(4) DC 1, 3, 9; -6** ⊕ *u 4, 10* **II(2) DC 1, 3; -10** ⊕ *a 4, 6, 10* **II(9); -1**

The total number of points is 14. If it is a tour, since there is one 2-circuit cycle and three acceptable cycles, our formula gives $p = \frac{n}{2} + 3t + a - 1 = 10 + 3 + 3 = 16.$ We thus have a derangement whose value is 55.

*u* **IV(2)** *3* **DC 1, 2, 8, 9, 10; -5**

No cycle containing a point in **9** or **8** can be patched to **IV(2).**

No cycle containing a point in **1** can be patched. The smallest-valued cycles containing precisely points in **2**, **4** and **10** are **II(i)**, i = 3, 4, 5, 6. We patch **II(3)** to **IV(2)**: *u* **IV(2)** *3* **DC 1, 2, 8, 9, 10; -5** ⊕ *u 3, 4* **II(3) DC 1, 8, 9; -2**.

**II(9)** containing a point in **9** can be patched:

*u* **IV(2)** *3* **DC 1, 2, 8, 9, 10; -5** ⊕ *u 3, 4* **DC 1, 8, 9; -2** ⊕ *a 3, 4, 10* **II(9) DC 1, 8; 7**.

No cycle containing points in **8** can be patched.

No tour or derangement.

*u* **III(1)** *5* **DC 1, 2, 4, 8, 9, 10; -15**

No point in **9** can be patched. No point in **8** can be patched. No point in **1** can be patched.

**IV(4)** can be patched.

*u* **III(1)** *5* **DC 1, 2, 4, 8, 9, 10; -15** ⊕ *a* **IV(4)** *5, 6* **DC 1, 8, 9; -14**

No 4-cycle can be patched. The only 3-cycle containing a point in **8** can be patched. **II(9)** can be patched:

*u* **III(1)** *5* **DC 1, 2, 4, 8, 9, 10; -15** ⊕ *a* **IV(4)** *5, 6* **DC 1, 8, 9; -14** ⊕ *a 5, 6, 10* **DC 1, 8; -5**

**III(2)** and **II(10)** both contain a point in **8** but they can't be patched. Thus, we can't get a derangement or tour.

**II(7)** can be patched to **III(1)**.

*u* **III(1)** *5* **DC 1, 2, 4, 8, 9, 10; -15** ⊕ *u* **II(5)** *4, 5, 7* **DC 1, 8, 9, 10; -10**

**IV(4)** can't be patched. Thus, no 4-cycle can be patched. From our previous case, no point in **8** can be patched.

**III(5)** can also be patched to **III(1).**

*u* **III(1)** *5* **DC 1, 2, 4, 8, 9, 10; -15** ⊕ *a* **III(5)** *5, 7* **DC 1, 8, 9, 10; -6**

No 2-cycle 3-cycle, 4-cycle or 5-cycle can be patched. No derangement or tour.

*u* **III(2)** *5* **DC 2, 3, 4, 6, 9, 10; -9**

No point in **9** can be patched to **III(2)**. **IV(5)** can be patched to **III(2)**.

*u* **III(2)** *5* **DC 2, 3, 4, 6, 9, 10; -9** ⊕ *a 5, 7* **IV(5) DC 2, 4, 9; -7**

**II(9)** can be patched:

*u* **III(2)** *5* **DC 2, 3, 4, 6, 9, 10; -9** ⊕ *a 5, 7* **IV(5) DC 2, 4, 9; -7** ⊕ *a 5, 7, 10* **II(9) DC 2, 4; 2**

No 2-cycle or 3-cycle can be patched that contains a point in **2** or **3**.

No derangement or tour.

*u* **II(1) 5 C 7, 8; -14** = *u* **II(1)** *5* **DC 1, 2, 3, 4, 6, 9, 10; -14**

Each of the following can be patched to **II(1)**:

**(i)** *u* **II(7)** *4* **C 2, 7; 5. (ii)** *a* **III(5) C 2, 4, 7; 9. (iii)** *a* **V(4) C 2, 3, 6, 7, 10; 7.**

**(iv)** *a* **VII(6) C 1, 2, 3, 4, 6, 8, 10; 2**.

**II(1)(i)** *u* **II(1)** *5* **DC 1, 2, 3, 4, 6, 9, 10; -14** ⊕ *u* **II(7)** *4, 5, 7* **DC 1, 2, 3, 6, 9, 10; -9**

**II(1)(ii)** *u* **II(1)** *5* **DC 1, 2, 3, 4, 6, 9, 10; -14** ⊕ *a* **III(5)** *5, 7* **DC 1, 3, 6, 9, 10; -5**

**II(1)(iii)** *u* **II(1)** *5* **DC 1, 2, 3, 4, 6, 9, 10; -14** ⊕ *a* **V(4)** *5, 7* **DC 1, 4, 9; -7**

**II(1)(iv)** *u* **II(1)** *5* **DC 1, 2, 3, 4, 6, 9, 10; -14** ⊕ *a* **VII(6) 5, 7, 8 DC 9; -12**

**II(1)(i)**: No 2-cycle, 3-cycle, 4-cycle, 5-cycle or 6-cycle can be patched. No derangement or tour.

**II(1)(ii)**: Same as in **II(1)(i)**.

**II(1)(iii)**: Only **II(8)** and **II(9)** can be patched.

**(a) II(8)**: *u* **II(1)** *5* **DC 1, 2, 3, 4, 6, 9, 10; -14** ⊕ *a* **V(4)** *5, 7* **DC 1, 4, 9; -7**⊕ *a* **II(8)** *5, 7, 2* **DC 1, 9; -2**

**(b) II(9)**: *u* **II(1)** *5* **DC 1, 2, 3, 4, 6, 9, 10; -14** ⊕ *a* **V(4)** *5, 7* **DC 1, 4, 9; -7** ⊕ *a* **II(9)** 5, 7, 10 **DC 1, 4; 2**

We can patch **II(9)** to **(a)** and **II(8)** to **(b).** In both cases, we obtain **DC 1.** No 2-cycle contains a point in **1**.

No derangement or tour.

**II(1)(iv)**: **II(1)(iv)** *u* **II(1)** *5* **DC 1, 2, 3, 4, 6, 9, 10; -14** ⊕ *a* **VII(6)** *5, 7, 8* **DC 9; -12**

We can patch **II(9)** to **II(1)(iv)** to obtain either a derangement or a tour. We obtain

**II(1)(iv)** *u* **II(1)** *5* **DC 1, 2, 3, 4, 6, 9, 10; -14** ⊕ *a* **VII(6)** *5, 7, 8* **DC 9; -12** ⊕ **II(9)** *5, 7, 8, 10*; **-3**

We have a total of thirteen points in **II(1), VII(6)** and **II(9).** From our formula,

$p = \frac{n}{2} + 3t + a - 1 = 10 + 3 + 1 = 14$ . We thus have a derangement whose value is 53.

*u* **II(2)** *3* **C 1, 6** = *3* **DC 2, 4, 5, 7, 8, 9, 10; -4**

The following can be patched to **II(2)**:

**II(2)(i)**: *u* **II(2)** *3* **DC 2, 4, 5, 7, 8, 9, 10; -4** ⊕ *u* **VI(6)** *2, 3, 6* **DC 4, 9; -4**

We can patch **II(9)**.

*u* **II(2)** *3* **DC 2, 4, 5, 7, 8, 9, 10; -4** ⊕ *u* **VI(6)** *2, 3, 6* **DC 4, 9; -4** ⊕ *a* **II(9)** *2, 3, 6, 10* **DC 4; 5**

No 2-cycle can be patched. No derangement or tour.

**II(2)(ii)**: *u* **II(2)** *3* **DC 2, 4, 5, 7, 8, 9, 10; -4** ⊕ *a* **V(2)** *1, 3* **DC 2, 9, 10; 0**

We can patch **III(5)** or **II(i)** (i = 3, 4, 5, 6).

*u* **II(2)** *3* **DC 2, 4, 5, 7, 8, 9, 10; -4** ⊕ *a* **V(2)** *1, 3* **DC 2, 9, 10; 0** ⊕ *a* **III(5)** *1, 3, 4* **DC 9; 9**

We could patch **II(9)** to obtain a derangement or a tour. However, its value would be 74.

Each of the **II(i)**'s would have the same effect. We therefore choose **II(3)**.

*u* **II(2)** *3* **DC 2, 4, 5, 7, 8, 9, 10; -4** ⊕ *a* **V(2)** *1, 3* **DC 2, 9, 10; 0** ⊕ *u* **II(3)** *1, 3, 4* **DC 9; 3**

We can patch **II(9)** to obtain a derangement or tour whose value is 68.

**II(2)(iii)**: *u* **II(2)** *3* **DC 2, 4, 5, 7, 8, 9, 10; -4** ⊕ *a* **V(5)** *3, 6* **DC 8, 9; 6**

We can patch **II(9)** which would yield **DC 9; 15**. But we couldn't patch a cycle containing a point in **8**.

Thus, we can't obtain a derangement or a tour.

**II(2)(iv)**: *u* **II(2)** *3* **DC 2, 4, 5, 7, 8, 9, 10; -4** ⊕ *a* **IV(3)** *1, 3* **DC 2, 4, 9, 10; -4**

We can patch **III(6)** or **II(7)**.

*u* **II(2)** *3* **DC 2, 4, 5, 7, 8, 9, 10; -4** ⊕ *a* **IV(3)** *1, 3* **DC 2, 4, 9, 10; -4** ⊕ *a* **III(6)** *1, 3, 7* **DC 9, 10; 5**

**II(9)** can't be patched. Thus, we can't obtain a derangement or tour.

*u* **II(2)** *3* **DC 2, 4, 5, 7, 8, 9, 10; -4** ⊕ *a* **IV(3)** *1, 3* **DC 2, 4, 9, 10; -4** ⊕ *u* **II(7)** *1, 3, 4, 7* **DC 9, 10; 1**

**II(9)** can't be patched. Thus, we can't obtain a derangement or tour.

## REFERENCES

[1] Helsgaun, K., “An Effective Implementation of the Lin-Kernighan Traveling Salesman Heuristic”, DATALOGISKE SKRIFTER (Writings on Computer Science), No. 81, 1998, Roskilde University.

[2] Jonker, R. and Volgenant, T., “Transforming asymmetric into symmetric traveling salesman problems”, Oper. Res. Let., **2**, 161-163 (1983).